\font\eusm=eusm10                   % AMS_Tex Euler script font

                   % bold \Cal font

\font\eusms=eusm7                       % eusm scriptsize

\font\eusmss=eusm5                      % eusm tiny

\font\scriptsize=cmr7

\input amstex

\documentstyle{fic}
  \hsize=7.0truein
  \vsize=9.0truein
  \hoffset -0.1truein

\NoBlackBoxes

\define\Ac{{\Cal A}}                         % A caligraphic

\define\Afr{{\frak A}}                       % A fraktur

\define\ah{{\hat a}}                         % a hat

\define\bh{\hat b}                           % b hat

\define\Bof{B}                               % for B of H

\define\Cc{{\Cal C}}                         % C caligraphic

\define\clspan{\overline\lspan}              % closed linear span

\define\Cpx{\bold C}                         % complex numbers symbol

\define\diag{\text{\rm diag}}                % diag

\define\dif{\text{\it d}}                    % d for integration

\define\eqdef{{\;\overset\text{def}\to=\;}}     % eq def

\define\freeF{F}                             % free group F

\define\freeprodi{\mathchoice                % free product (index set I)
     {\operatornamewithlimits{\ast}
      _{\iota\in I}}
     {\raise.5ex\hbox{$\dsize\operatornamewithlimits{\ast}
      _{\sssize\iota\in I}$}\,}
     {\text{oops!}}{\text{oops!}}}

\define\hh{\hat h}                           % h hat

\define\Hil{{\mathchoice                     % Hilbert space 
     {\text{\eusm H}}
     {\text{\eusm H}}
     {\text{\eusms H}}
     {\text{\eusmss H}}}}

\define\Hilo{{\overset{\scriptsize o}        % Hilbert space with 'zero' on top
     \to\Hil}}

\define\htld{{\tilde h}}                     % h tilde

\define\Integers{\bold Z}                    % Integers

\define\kt{{\tilde k}}                       % k tilde

\define\Lambdao{{\Lambda\oup}}               % Lambda^o

\define\lspan{\text{\rm span}@,@,@,}         % span

\define\mfr{{\frak m}}                       % m fraktur

\define\MvN{{\Cal M}}                        % M caligraphic

\define\nm#1{||#1||}                         % norm symbol

\define\ntr{{\text{\rm ntr}}}                % ntr

\define\Naturals{{\bold N}}                  % Natural numbers symbol

\define\NvN{{\Cal N}}                        % N caligraphic

\define\onehat{\hat 1}                       % one hat

\define\oup{^{\text{\rm o}}}                 % ^o

\define\owedge{{                             % circle wedge
     \operatorname{\raise.5ex\hbox{\text{$
     \ssize{\,\bigcirc\llap{$\ssize\wedge\,$}\,}$}}}}}

\define\PvN{{\Cal P}}                        % P caligraphic

\define\QED{$\hfill$\qed\enddemo}            % QED

\define\QvN{{\Cal Q}}                        % Q caligraphic

\define\Real{{\bold R}}                      % Real numbers symbol

\define\Reals{{\bold R}}                     % Real numbers symbol

\define\restrict{\lower .3ex                 % restriction bar
     \hbox{\text{$|$}}}

\define\Sd{\text{\it Sd}}                    % Sd

\define\smd#1#2{\underset{#2}\to{#1}}          % summand with tracial weights

\define\smdp#1#2#3{\overset{#3}\to             % smd plus projection
     {\smd{#1}{#2}}}

\define\smp#1#2{\overset{#2}\to                % sm with only projection
     {#1}}                                     % (no tracial weight)

\define\tocdots                              % table of contents dots
  {\leaders\hbox to 1em{\hss.\hss}\hfill}    %  (see [Knuth], p. 223)

\define\Tr{\text{\rm Tr}}                    % Tr (trace)

\define\xh{\hat x}                           % x hat

  \newcount\notasecflag \notasecflag=0
  \def\notasec{\notasecflag=1}

  \newcount\secno \secno=0 \newcount\subsecno
  \def\newsec#1{\procno=0 \subsecno=0 \notasecflag=0
    \advance\secno by 1 \edef#1{\number\secno}
    \edef\currentsec{\number\secno}}
  \def\newsubsec#1{\procno=0 \advance\subsecno by 1 \edef#1{\number\subsecno}
    \edef\currentsec{\number\secno.\number\subsecno}}

  \newcount\appendixno \appendixno=0
  \def\newappendix#1{\procno=0 \notasecflag=0 \advance\appendixno by 1
    \ifnum\appendixno=1 \edef\appendixalpha{\hbox{A}}
      \else \ifnum\appendixno=2 \edef\appendixalpha{\hbox{B}} \fi
      \else \ifnum\appendixno=3 \edef\appendixalpha{\hbox{C}} \fi
      \else \ifnum\appendixno=4 \edef\appendixalpha{\hbox{D}} \fi
      \else \ifnum\appendixno=5 \edef\appendixalpha{\hbox{E}} \fi
      \else \ifnum\appendixno=6 \edef\appendixalpha{\hbox{F}} \fi
    \fi
    \edef#1{\appendixalpha}
    \edef\currentsec{\appendixalpha}}

  \newcount\procno \procno=0
  \def\newproc#1{\advance\procno by 1
   \ifnum\notasecflag=0 \edef#1{\currentsec.\number\procno}
   \else \edef#1{\number\procno}
   \fi}

  \newcount\tagno \tagno=0
  \def\newtag#1{\advance\tagno by 1 \edef#1{\number\tagno}}

\notasec
 \newproc{\MainOne}
  \newtag{\FPAoneAtwo}
  \newtag{\MzeroD}
  \newtag{\FPExamples}
 \newproc{\ExtrFreeProd}
 \newproc{\AlsoBofH}
 \newproc{\FPInfMany}
  \newtag{\FPMoreExamples}
 \newproc{\DefLambdao}
\newsec{\FreeSubcompl}
 \newproc{\FreelyCompl}
 \newproc{\FiniteStdEmbChar}
 \newproc{\StdEmbFacts}
 \newproc{\FreelySubcompl}
 \newproc{\FreelySubcomplWords}
 \newproc{\TwoProjections}
  \newtag{\twoprojections}
 \newproc{\FreelySubcomplChain}
  \newtag{\freelysubcomplchain}
 \newproc{\FreelySubcomplWell}
\newsec{\FreeAmalg} %2
 \newproc{\FreenessAmalg}
 \newproc{\FreenessAmalgCondExp}
 \newproc{\FreelySubcomplAmalg}
\newsec{\CondExp} %3
 \newproc{\CondExpFromDense}
 \newproc{\CondExpImplByProj}
  \newtag{\EhatIsP}
  \newtag{\EImplProj} %(4)
 \newproc{\FreeProdOfCondExp}
 \newproc{\FreeFactorCondExp}
\newsec{\ExtrAlmPerSt} %4
 \newproc{\AlmPer}
 \newproc{\AlmPerFreeProd}
 \newproc{\SpectralTh}
  \newtag{\spectraldecomp}
 \newproc{\Msg}
 \newproc{\CondExpAlmPer}
 \newproc{\SdProp}
 \newproc{\DefExtrAPSt}
 \newproc{\MustBeFactor}
 \newproc{\ExtrAPStDecomp}
 \newproc{\ExtrAPStCutDowns}
 \newproc{\WhenAlmPer}
  \newtag{\psiphi}
 \newproc{\InWhatSenseExtr}
\newsec{\FreeProdLemmas} %5
 \newproc{\FreeProdDirectSum}
 \newproc{\FreeProdDirectSumMatrixAlg}
 \newproc{\FreeProdDirectSumBoth}
 \newproc{\FreeProdMatUnits}
  \newtag{\hpgen}
  \newtag{\cupthree}
 \newproc{\FreeProdCutDown}
 \newproc{\FreelySubcomplCutDown}
\newsec{\FPwrtMAP} %6
 \newproc{\AddingOneComplimented}
  \newtag{\allNvNDef}
  \newtag{\NvNminus}
  \newtag{\NvNplus}
 \newproc{\AddingOneTorsion}
  \newtag{\vmodN}
  \newtag{\modNsum} %(11)
 \newproc{\AddingTwoComplimented}
  \newtag{\NvNtZeroOne}
  \newtag{\uOnea}
  \newtag{\uOneb}
  \newtag{\NvNtminus}
  \newtag{\NvNtplus}
 \newproc{\AddingTwoTorsion}
  \newtag{\TmodNsum} %(17)
 \newproc{\FPTypeIIIFactors} %6.5
\newsec{\FPfd} %7
 \newproc{\pAphi}
 \newproc{\PtSpecFd}
  \newtag{\TracialWt}
 \newproc{\FPExtrAPfd} %7.3
  \newtag{\Qfp}
  \newtag{\xAltProd}
  \newtag{\chiI} %(21)
 \newproc{\NotationFD}
 \newproc{\FPFD}
  \newtag{\DorNot}
  \newtag{\DoneDtwo}
  \newtag{\Mresult}
 \newproc{\MatrixUnitsMeet}
  \newtag{\QplusDcases}
  \newtag{\PTheta}
 \newproc{\StdEmbFPFD}
  \newtag{\NzeroInMzero}
\newsec{\Fullness} %8
 \newproc{\ExampleSd}
 \newproc{\UsuallyBarnett}
 \newproc{\MnMnFull}
 \newproc{\CentralizerPI}
\newsec{\Questions}
 \newproc{\Qisomlambda}
 \newproc{\Qcoreone}
 \newproc{\Qisomone}
 \newproc{\Qnonfull}
  
\def\Arveson{Arveson [1974]}
\def\Barnett{Barnett [1995]}
\def\ConnesZZThesis{Connes [1973]}
\def\ConnesZZAlmPer{Connes [1974]}
\def\DykemaZZFreeProdR{Dykema [1993(a)]}
\def\DykemaZZFreeDim{Dykema [1993(b)]}
\def\DykemaZZInterp{Dykema [1994(a)]}
\def\DykemaZZTM{Dykema [1994(b)]}
\def\DykemaZZAmalg{Dykema [1995(a)]}
\def\DykemaZZAlmPer{Dykema [1995(b)]}
\def\DykemaZZSNU{Dykema [SNU]}
\def\GolodetsNessonov{Golodets and Nessonov [1987]}
\def\MvNZZiv{Murray and von Neumann [1943]}
\def\PedersenTakesaki{Pedersen and Takesaki [1973]}
\def\PopaZZUnivJonesAlg{Popa [1993]}
\def\Pukansky{Pukansky [1953]}
\def\RadulescuZZFundGp{R\u adulescu [1992(a)]}
\def\RadulescuZZAmalgSubfactors{R\u adulescu [1994]}
\def\RadulescuZZOneParamGp{R\u adulescu [1992(b)]}
\def\RadulescuZZlambda{R\u adulescu [preprint]}
\def\ReedSimonZZI{Reed and Simon [1972]}

\def\TakesakiZZDualityCrossProd{Takesaki [1973]}
\def\Tomiyama{Tomiyama [1957]}
\def\VoiculescuZZSymetries{Voiculescu [1985]}
\def\VoiculsecuZZCircSemiCirc{Voiculsecu [1990]}
\def\VDNbook{Voiculescu, Dykema and Nica [1992]}

\topmatter

  \title Free products of finite dimensional and other von Neumann algebras
         with respect to non-tracial states.
  \endtitle

  \rightheadtext{Free products}

  \author Kenneth J\. Dykema   \\
   \address Department of Mathematics \\ University of California \\
            Berkeley CA 94720--3840 \\ ---------- \\
            present address: \\
            Department of Mathematics and Computer Science \\
            Odense Universitet, Campusvej 55 \\ DK--5230 Odense M, Denmark \\
            e-mail: dykema\@imada.ou.dk
   \endaddress
  \endauthor

  \leftheadtext{K.J\. Dykema}

  \thanks Supported by a National Science
          Foundation Postdoctoral Fellowship and by a Fields Institute
          Fellowship. \endthanks

  \abstract 
The von Neumann algebra free product of arbitary finite dimensional von
Neumann algebras with respect to arbitrary faithful states, at least
one of which is not a trace, is found to be a type~III factor possibly
direct sum a finite dimensional algebra.
The free product state on the type~III factor is what we
call an extremal
almost periodic state, and has centralizer isomorphic to $L(\freeF_\infty)$.
This allows further classification the type~III factor and provides another
construction of full type~III$_1$ factors having arbitrary $\Sd$~invariant of
Connes.
The free products considered in this paper are not limited to free products of
finite dimensional algebras, but can be of a quite general form.
  \endabstract

\endtopmatter

\document \baselineskip=18pt

\head Introduction. \endhead

In the context of Voiculescu's theory of freeness,
(see~\cite{\VDNbook}), the free product of operator algebras is a basic
operation, which is to freeness what the tensor product is to independence.
Free products of several sorts of von Neumann algebras have been investigated,
including free products of finite dimensional and hyperfinite von Neumann
algebras with respect to traces, (\cite{\VoiculsecuZZCircSemiCirc},
\cite{\DykemaZZFreeProdR}, \cite{\DykemaZZInterp}, \cite{\DykemaZZFreeDim}),
certain amalgamated free products with respect to traces,
(\cite{\PopaZZUnivJonesAlg}, \cite{\RadulescuZZAmalgSubfactors},
\cite{\DykemaZZAmalg}) and free products
of various von Neumann algebras with respect to non--tracial states,
(\cite{\RadulescuZZlambda}, \cite{\Barnett}, \cite{\DykemaZZTM}).
In~\cite{\RadulescuZZlambda}, the free product of the diffuse abelian von
Neumann algebra and $M_2(\Cpx)$ with respect to the state
$\Tr\left(\left(\smallmatrix\frac1{1+\lambda}&0\\0&\frac\lambda{1+\lambda}
\endsmallmatrix\right)\cdot\right)$ for $0<\lambda<1$ was shown to be a
type~III$_\lambda$ factor with core $L(\freeF_\infty)\otimes\Bof(\Hil)$.
In~\cite{\Barnett} and~\cite{\DykemaZZTM}, some general results about free
products of von Neumann algebras with respect to non--tracial states were
proved, but several natural questions were left unanswered.
For example, in the free product
$$ (\MvN,\phi)=(M_2(\Cpx),
\Tr\left(\left(\matrix\frac1{1+\lambda}&0\\0&\frac\lambda{1+\lambda}
\endmatrix\right)\cdot\right))*
(M_2(\Cpx),
\Tr\left(\left(\matrix\frac1{1+\mu}&0\\0&\frac\mu{1+\mu}
\endmatrix\right)\cdot\right)) $$
for $0<\lambda,\mu<1$,
$\MvN$ was by~\cite{\DykemaZZTM} known to be a factor only for certain values
of $\lambda$ and $\mu$, and for other values this remained unknown.
One of the consequences of our main result is that such a free product is
always a factor, and the centralizer of the free product state $\phi$ is
isomorphic to $L(\freeF_\infty)$.

In this paper, we investigate the free product of finite dimensional von
Neumann algebras (and others) with respect to arbitrary faithful states
that are not traces.
The main result can be summarized as follows:
\proclaim{Theorem \MainOne}
Let
$$ (\MvN,\phi)=(A_1,\phi_1)*(A_2,\phi_2) \tag{\FPAoneAtwo} $$
be the von Neumann algebra free product of finite dimensional algebras with
respect to faithful states, at least one of which is non--tracial.
Then
$$ \MvN=\MvN_0\oplus D\quad\text{or}\quad\MvN=\MvN_0, \tag{\MzeroD} $$
where $D$ is a finite dimensional algebra and where $\MvN_0$ is a type~III
factor.
The type~I part, $D$, can be found from knowledge of $\phi_1$ and $\phi_2$, as
described below.
The restriction, $\phi_0\eqdef\phi\restrict_{\MvN_0}$ to the type~III
part, is an almost periodic state (or functional) whose centralizer
is isomorphic to the II$_1$--factor $L(\freeF_\infty)$.
The point spectrum of the modular operator of $\phi_0$,
$\Delta_{\phi_0}$, is equal to the
subgroup of $\Real^*_+$ generated by the union of the point spectra of
$\Delta_{\phi_1}$ and of $\Delta_{\phi_2}$.
Thus, in Connes' classification of $\MvN_0$ as type~III$_\lambda$,
one can find $\lambda$ and always $0<\lambda\le1$.
\endproclaim
Note that in~(\FPAoneAtwo) when both $\phi_1$ and $\phi_2$ are
traces, $\MvN$ was found
in~\cite{\DykemaZZFreeDim}.
In that case, $\MvN$ was of a similar form, except that the non--(type~I) part,
$\MvN_0$, was then a II$_1$ factor related to free groups,
(see~\cite{\RadulescuZZAmalgSubfactors} and~\cite{\DykemaZZInterp}).
In fact, the description in~\cite{\DykemaZZFreeDim} of $D$
is a special case of the description that follows.

Let us now describe the type~I part, $D$, in~(\MzeroD).
(A more precise discussion is found in~\S\FPfd.)
Given a faithful state, $\psi$, on a finite dimensional von Neumann algebra,
$D=\bigoplus_{j=1}^K M_{n_j}(\Cpx)$, we write
$$ (D,\psi)=\bigoplus_{j=1}^K \smd{M_{n_j}(\Cpx)}
{\alpha_{j,1},\ldots,\alpha_{j,n_j}} $$
to mean that the restriction of $\psi$ to the $j$th summand of $D$ is given by
a diagonal density matrix with $\alpha_{j,1},\ldots,\alpha_{j,n_j}$ down the
diagonal.
If
$$ \aligned
(A_1,\phi_1)&=\bigoplus_{j=1}^{K_1}\smd{M_{n_j}(\Cpx)}
{\alpha_{j,1},\cdots,\alpha_{j,n_j}} \\
(A_2,\phi_2)&=\bigoplus_{j=1}^{K_2}\smd{M_{m_j}(\Cpx)}
{\beta_{j,1},\cdots,\beta_{j,m_j}}
\endaligned  $$
and if $\MvN$ is their free product as in~(\FPAoneAtwo), to determine whether
$\MvN$ is a factor and, when it is not, to find the type~I part of $\MvN$, we
first ``mate'' every simple summand of $A_1$ with every simple summand of $A_2$
and examine the offspring, if any.
When the $j$th summand of $A_1$, namely
$\smd{M_{n_j}(\Cpx)}{\alpha_{j,1},\cdots,\alpha_{j,n_j}}$, is mated with the
$k$th summand of $A_2$, namely
$\smd{M_{m_k}(\Cpx)}{\beta_{k,1},\cdots,\beta_{k,m_k}}$,
then there can be offspring only if at least one of $n_j$ and $m_k$ is
equal to $1$, {\it i.e\.} if at least one of the simple summands is just a
copy of $\Cpx$ in the corresponding $A_\iota$.
In that case, say if $n_j=1$, then there is offspring if and only if
$$ \sum_{i=1}^{m_k}\frac1{\beta_{k,i}}<\frac1{1-\alpha_{j,1}}, $$
and then the offspring is
$$ \smd{M_{m_k}(\Cpx)}{\gamma_1,\ldots,\gamma_{m_k}}\quad\text{where}\quad
\gamma_i=\beta_{k,i}\left(1-(1-\alpha_{j,1})\sum_{p=1}^{m_j}\frac1{\beta_{k,p}}
\right). $$
Analogous offspring results if $m_k=1$ instead of $n_j=1$.
A system of matrix units for the offspring is the meet, $\owedge$,
defined in~\cite{\DykemaZZAmalg},
of systems of matrix units for the simple summands that were mated,
$M_{n_j}(\Cpx)$ and $M_{m_k}(\Cpx)$, (see~\S7).
Thus the support projection of
this offspring lies under the support projection of $M_{n_j}(\Cpx)$ and the
support projection of $M_{m_k}(\Cpx)$.
The type~I part of $\MvN$ is equal to the direct sum of all offspring
produced in this way by all possible matings of the original simple summands.
If none of the matings were successful in producing offspring, then $\MvN$ is a
factor.

For example, we have
$$ \spreadlines{3ex} \align
\bigg(\smd\Cpx{\frac15}\oplus\smd\Cpx{\frac45}\bigg)*
 \smd{M_2(\Cpx)}{\frac23,\frac13}
 &=\MvN_{\text a}\oplus\smd{M_2(\Cpx)}{\frac1{15},\frac1{30}}
 \tag{\FPExamples a} \\
\bigg(\smd\Cpx{\frac14}\oplus\smd\Cpx{\frac34}\bigg)*
 \smd{M_2(\Cpx)}{\frac23,\frac13}
 &=\MvN_{\text b} \tag{\FPExamples b} \\
\smd{M_2(\Cpx)}{\frac14,\frac34}*
 \smd{M_2(\Cpx)}{\frac23,\frac13}
 &=\MvN_{\text c}
 \tag{\FPExamples c} \\
\bigg(\smd{M_2(\Cpx)}{\frac1{10(1+\sqrt8)},\frac{\sqrt8}{10(1+\sqrt8)}}
 \oplus\smd\Cpx{\frac9{10}}\bigg)
 *\bigg(\smd{M_2(\Cpx)}{\frac13,\frac16}\oplus\smd\Cpx{\frac12}\bigg)
 &=\MvN_{\text d}\oplus\smd{M_2(\Cpx)}{\frac1{30},\frac1{60}}
 \oplus\smd\Cpx{\frac25}
 \tag{\FPExamples d} \\
\bigg(\smd{M_2(\Cpx)}{\frac4{20},\frac1{20}}\oplus\smd\Cpx{\frac34}\bigg)
 *\bigg(\smd{M_2(\Cpx)}{\frac13,\frac16}\oplus\smd\Cpx{\frac12}\bigg)
 &=\MvN_{\text e}\oplus\smd\Cpx{\frac14}
 \tag{\FPExamples e} \\
\bigg(\smd\Cpx{\frac1{41}}\oplus\smd\Cpx{\frac{40}{41}}\bigg)
 *\bigg(\smd{M_3(\Cpx)}{\frac1{16},\frac1{16},\frac18}
 \oplus\smd{M_2(\Cpx)}{\frac1{40},\frac9{40}}
 \oplus\smd{M_2(\Cpx)}{\frac18,\frac18}
 \oplus\smd\Cpx{\frac14}\bigg)
 &=\MvN_{\text f}\oplus\smd{M_3(\Cpx)}{\frac1{656},\frac1{656},\frac1{328}}
 \oplus\smd{M_2(\Cpx)}{\frac{25}{328},\frac{25}{328}}
 \oplus\smd\Cpx{\frac{37}{164}}, \tag{\FPExamples f}
\endalign $$
where $\MvN_{\text a}$, $\MvN_{\text b}$, $\MvN_{\text c}$, $\MvN_{\text d}$,
$\MvN_{\text e}$ and $\MvN_{\text f}$ are type~III factors.

Suppose that $\MvN$ is not a factor, {\it i.e\.} that some matings result in
offspring, then one easily sees that there is a {\it dominant projection}, $p$,
which is by definition a minimal and central projection of $A_\iota$, for
$\iota=1$ or $\iota=2$, such that all of the offspring are the progeny of the
mating of $p$ with simple summands from $A_{\iota'}$ where $\iota'\neq\iota$.
For example, in~(\FPExamples a), (\FPExamples d) and~(\FPExamples f), the only
dominant projection is $0\oplus1$ in the algebra on the left side of the $*$,
whereas in~(\FPExamples f), both $0\oplus1$ on the left and $0\oplus1$ on the
right side of $*$ are dominant projections.  The dominant projection $p\in
A_\iota$, if it exists, must be the largest (with respect to $\phi_\iota$) of
the projections in $A_\iota$ that are both minimal and central.

We now turn to a discussion of the type~III factor, $\MvN_0$, appearing
in Theorem~\MainOne.
\cite{\ConnesZZThesis} defined a state, $\phi$, on a von Neumann algebra,
$\MvN$, to be almost periodic if
$L^2(\MvN,\phi)$ has a basis of eigenvectors for the modular operator,
$\Delta_\phi$, of $\phi$.
In that case, Arveson--Connes spectral theory (\cite{\Arveson},
\cite{\ConnesZZThesis}) provides a decomposition of
$\MvN$ as a direct sum of spectral subspaces.
If, in addition, the centralizer, $\MvN_\phi$, of $\phi$ is a factor, then, it
is easily seen, the
point spectrum of $\Delta_\phi$ is a subgroup of $\Real_+^*$ and every spectral
subspace is of the form $\MvN_\phi v$ or $v^*\MvN_\phi$ for an
isometry $v$.
In analogy with Connes' results about full factors, \cite{\ConnesZZAlmPer}, we
call an almost periodic state, $\phi$, having centralizer that is a factor an
{\it extremal almost periodic} state, and we define $\Sd(\phi)$ to be the point
spectrum of $\Delta_\phi$.
(At the end of~\S5, we show that term ``extremal'' is in some sense
appropriate.)
By Corollary~3.2.7 of~\cite{\ConnesZZThesis}, from the fact that $\MvN_\phi$ is
a factor we have also that Connes' invariant $S(\MvN)$ is equal to the closure
of $\Sd(\phi)$.
Hence in Theorem~\MainOne, $\phi_0$ is extremal
almost periodic and $\Sd(\phi_0)$ is found from the initial data, namely
$\phi_1$ and $\phi_2$.
Moreover, regarding Connes' classification of $\MvN_0$ as type~III$_\lambda$,
we have:
$\MvN_0$ is type III$_\lambda$ for $0<\lambda<1$ if
$\Sd(\phi_0)=\lambda^\Integers=\{\lambda^n\mid n\in\Integers\}$ and
otherwise $\Sd(\phi_0)$ is dense in $\Reals_+^*$, and
$\MvN_0$ must be type III$_1$.

For example,
\roster
\item"in (\FPExamples a)," $\MvN_{\text a}$ is a type III$_\lambda$ factor for
$\lambda=1/2$;
\item"in (\FPExamples b)," $\MvN_{\text b}$ is a type III$_\lambda$ factor for
$\lambda=1/2$;
\item"in (\FPExamples c)," $\MvN_{\text c}$ is a type III$_1$ factor and
the free product state, $\phi$, has $\Sd(\phi)$ equal to the (dense) subgroup
of $\Real_+^*$ generated by $1/3$ and $1/2$;
\item"in (\FPExamples d)," $\MvN_{\text d}$ is a type III$_\lambda$ factor for
$\lambda=1/\sqrt2$;
\item"in (\FPExamples e)," $\MvN_{\text e}$ is a type III$_\lambda$ factor for
$\lambda=1/2$;
\item"in (\FPExamples f)," $\MvN_{\text f}$ is a type III$_1$ factor and the
free product state, $\phi$, has $\Sd(\phi)$ equal to the (dense) subgroup of
$\Real_+^*$ generated by $1/2$ and $1/9$.
\endroster

Extremal almost periodic states play a prominent role in the proof of
Theorem~\MainOne.
An example of the sort of results proved and used in the course of this paper
is the following theorem.
\proclaim{Theorem \ExtrFreeProd}
Suppose for $\iota=1,2$ that $A_\iota$ is a separable type~III factor and
$\phi_\iota$ is an extremal almost periodic state on $A_\iota$ whose
centralizer is isomorphic to the hyperfinite II$_1$ factor, $R$, or to
$L(\freeF_\infty)$.
Let
$$ (\MvN,\phi)=(A_1,\phi_1)*(A_2,\phi_2). $$
Then $\MvN$ is a type~III factor, $\phi$ is extremal almost periodic and
$\Sd(\phi)$ is the subgroup of $\Reals_+^*$ generated by
$\Sd(\phi_1)\cup\Sd(\phi_2)$.
\endproclaim

By taking inductive limits, Theorem~\MainOne{} can be extended to more
general algebras.
\proclaim{Theorem \AlsoBofH}
Let
$$ (\MvN,\phi)=(A_1,\phi_1)*(A_2,\phi_2), $$
where $\phi_1$ and $\phi_2$ are
faithful states, at least one of which is not a trace, where
$A_1$ and $A_2$ are nontrivial algebras, each of which is one of the following:
\roster
\item"(i)" finite dimensional
\item"(ii)" $\Bof(\Hil)$ for seperable, infinite dimensional Hilbert space
           $\Hil$
\item"(iii)" a type~III factor on which $\phi_\iota$ is extremal almost
             periodic and has centralizer isomorphic to the hyperfinite II$_1$
             factor, $R$, or to $L(\freeF_\infty)$;
\item"(iv)" a diffuse hyperfinite von Neumann algebra on which $\phi_\iota$ is
            a trace
\item"(v)" an interpolated free group factor (of type~II$_1$),
            see \cite{\RadulescuZZAmalgSubfactors}, \cite{\DykemaZZInterp}
\item"(vi)" a possibly infinite direct sum of algebras from (i)-(v).
\endroster
Then
$$ \MvN=\MvN_0\oplus D\quad\text{or}\quad\MvN=\MvN_0, $$
with $\MvN_0$ and $D$ as described in Theorem~\MainOne, (and more explicitly
below).
\endproclaim
In Theorem~\AlsoBofH, one finds the finite dimensional part, $D$, of $\MvN$, if
it exists, by mating every matrix algebra direct summand of $A_1$ with every
matrix algebra summand of $A_2$, as described above, and $D$ is the sum of the
offspring of these matings.
\proclaim{Theorem \FPInfMany}
Let
$$ (\MvN,\phi)=\operatornamewithlimits*_{\iota\in I}(A_\iota,\phi_\iota), $$
where $I$ is finite or countably infinite, $|I|\ge2$, each $\phi_\iota$ is a
faithful state on $A_\iota$, at least one $\phi_\iota$ is not a trace and each
$A_\iota$ is of the form described in Theorem~\AlsoBofH.
Then
$$ \MvN=\MvN_0\oplus D\quad\text{or}\quad\MvN=\MvN_0, $$
with $\MvN_0$ and $D$ as described in Theorem~\MainOne, (and more explicitly
below).
\endproclaim
In this case, one finds the finite dimensional part, $D$, of $\MvN$, if it
exists, by, for every choice of a matrix algebra direct summand,
say $\smd{M_{n_\iota}(\Cpx)}{\alpha_{\iota,1},\ldots,\alpha_{\iota,n_\iota}}$
from each $A_\iota$, mating them all together, and letting $D$ be the sum of
the offspring resulting from all such matings.
A mating of one summand from each $A_\iota$ results in offspring only if all
but one of the $n_\iota$
equal $1$, say $n_\iota=1$ if $\iota\in I\backslash\{\iota'\}$, and
then offspring results if and only if
$$ \frac1{\sum_{\iota\in I\backslash\{\iota'\}}(1-\alpha_{\iota,1})}
>\sum_{j=1}^{n_\iota'}\frac1{\alpha_{\iota',j}} $$
and the offspring is
$$ \smd{M_{n_{\iota'}}(\Cpx)}{\gamma_1,\ldots,\gamma_{n_{\iota'}}}
\quad\text{where}\quad
\gamma_i=\alpha_{\iota',i}\left(
1-\biggl(\sum_{\iota\in I\backslash\{\iota'\}}(1-\alpha_{\iota,1})\biggr)
\sum_{p=1}^{n_{\iota'}}\frac1{\alpha_{\iota',p}}\right). $$
We also have, as in Theorem~\MainOne, that $\MvN_0$ is a type~III factor,
$\phi_0=\phi\restrict_{\MvN_0}$ is extremal almost periodic and $\Sd(\phi_0)$
is equal to the subgroup of $\Reals_+^*$ generated by the union over
$\iota\in I$ of the point spectra of $\Delta_{\phi_\iota}$.

For example,
$$ \gather
{\align (\Bof(\Hil),
 \Tr(\left(\smallmatrix1/2&*&&\\&1/4&&\\&&1/8&\\&&&\ddots\endsmallmatrix\right)
 \cdot))&*
 (\Bof(\Hil),
 \Tr(\left(\smallmatrix3/4&&&\\&3/16&&\\&&3/64&\\&&&\ddots
 \endsmallmatrix\right)
 \cdot))
 =\MvN_{\text g} \tag{\FPMoreExamples g} \\ \vspace{2ex}
\biggl((\Bof(\Hil),
 \Tr(\left(\smallmatrix1/2&&&\\&1/4&&\\&&1/8&\\&&&\ddots\endsmallmatrix\right)
 \cdot))\oplus\smd\Cpx{9/10}\biggr)&*\smd{M_2(\Cpx)}{\frac34,\frac14}
 =\MvN_{\text h}\oplus\smd{M_2(\Cpx)}{\frac7{20},\frac7{60}}
 \tag{\FPMoreExamples h}
 \endalign} \\ \vspace{3ex}
\operatornamewithlimits{\ast}_{n\ge1}
 \smd{M_2(\Cpx)}{\frac1{1+\root n\of2},\frac{\root n\of2}{1+\root n\of2}}
 =\MvN_{\text k},
 \tag{\FPMoreExamples k}
\endgather $$
where $\MvN_{\text g}$, $\MvN_{\text h}$
and $\MvN_{\text k}$ are type~III factors and
\roster
\item"in (\FPMoreExamples g)," $\MvN_{\text g}$ is a type III$_\lambda$ factor
for $\lambda=1/2$;
\item"in (\FPMoreExamples h)," $\MvN_{\text h}$ is a type III$_1$ factor and
the free product state restricted to $\MvN_{\text h}$, call it $\phi_0$,
has $\Sd(\phi_0)$ equal to the subgroup of
$\Real_+^*$ generated by $1/2$ and $1/3$;
\item"in (\FPMoreExamples k)," $\MvN_{\text k}$ is a type III$_1$ factor and
the free product state, $\phi$,
has $\Sd(\phi)$ equal to the subgroup of
$\Real_+^*$ generated by $\{\root n\of2\mid n\ge1\}$.
\endroster

The proofs of the above theorems involve what is loosely speaking a ``type~III
version'' of the following fact from group theory.
Let $F_2=\langle a,b\rangle$ be the free group with free generators~$a$
and~$b$.
Let $\pi:F_2\rightarrow\Integers$ be the group homomorphism such that
$\pi(a)=0$ and $\pi(b)=1$.
Then $\ker\pi\cong F_\infty$ and $\ker\pi$ is freely generated by
$(b^nab^{-n})_{n\in\Integers}$.
By ``type~III version'' we mean that trying to find the centralizer of certain
free product states is like trying to find $\ker\pi$ above, where instead of
$b$ we have a nonunitary isometry.
Of course, there are some technically more involved aspects, but the basic idea
is as easy as this.
We plan to write up the proof of a special case, namely $M_2(\Cpx)*M_2(\Cpx)$
with respect to arbitrary faithful states, in~\cite{\DykemaZZSNU}.
This proof will contain most of the essential ideas of the proof of the general
case, but should more transparent, and we would recommend it to those
interested in an introduction to the techniques.

Many results of this paper are proved using the technique of {\it free
etymology}.
This technique consists of proving the freeness of certain algebras by looking
at words in those algebras, and investigating the ``roots'' of those words.
To be more specific, suppose we want to show that a family
$(A_\iota)_{\iota\in I}$  of subalgebras of $A$ is free in $(A,\phi)$,
and each $A_\iota$ is defined in terms of
some other subalgebras $(B_j)_{j\in J}$ that we know are free.
Since we work so much with words, let us use a special notation, which we have
used before in~\cite{\DykemaZZFreeDim}, \cite{\DykemaZZInterp}
and~\cite{\DykemaZZAmalg}.
\proclaim{Definition \DefLambdao}\rm
Suppose $A$ is an algebra, $\phi:A\rightarrow\Cpx$ is linear.
For $B\subseteq A$ any subalgebra, we define $B\oup$ (pronounced $B$--bubble)
to be $B\oup=B\cap\ker\phi$.
Suppose $X_\iota\subseteq A$, are subsets, ($\iota\in I$).
The set of {\it reduced words} in $(X_\iota)_{\iota\in I}$ is denoted
$\Lambdao((X_\iota)_{\iota\in I})$, and is defined to be the set of
$a_1a_2\cdots a_n$ such that $n\ge1$, $a_j\in X_{\iota_j}$ and
$\iota_1\neq\iota_2\neq\cdots \iota_n$.
If $|I|=2$, we may call them {\it alternating products}.
\endproclaim
Now a free etymology proof of the freeness of $(A_\iota)_{\iota\in I}$ would go
as follows:
we need only show $\phi(x)=0$
for every $x\in\Lambdao((A_\iota\oup)_{\iota\in I})$;
after doing some work, one shows that $x$ is equal to, for example if $A$ is a
von Neumann algebra, a
strong--operator limit of sums of reduced words in $(B_j\oup)_{j\in J}$,
and hence one concludes that $\phi(x)=0$.

In~\S\Fullness, using a result of~\cite{\Barnett}, we show in many cases that
$\MvN_0$ in Theorems~\MainOne,~\AlsoBofH{} and~\FPInfMany{} is a full factor.
Thus, for these factors, Connes' invariant for full factors, $\Sd(\MvN_0)$, is
equal to $\Sd(\phi_0)$, and thereby, for instance by taking free products as in
example~(\FPMoreExamples k), we obtain for an arbitrary countable dense
subgroup, $\Gamma$, of
$\Real_+^*$, a free product factor that is a full type~III$_1$ factor,
has $\Gamma$ as its $\Sd$~invariant and has extremal almost periodic state with
centralizer isomorphic to $L(\freeF_\infty)$.
See Corollary~4.4 of~\cite{\ConnesZZAlmPer} and~\cite{\GolodetsNessonov} for
other constructions of full~III$_1$ factors with arbitrary $\Sd$~invariant.

\vfil
\noindent{\bf Acknowledgements.}
Most of this work was done while I was at the Fields Institute during
the special year in Operator Algebras, 1994/95.  
I would like to thank George Elliott for the organizing that marvelous year.

\vbox{
\noindent
\line{\bf Contents:\hfil}
\line{\S\FreeSubcompl.  Free Subcomplementation\tocdots p.7}
\line{\S\FreeAmalg.  Freeness with amalgamation over a subalgebra\tocdots p.11}
\line{\S\CondExp.  Conditional Expectations\tocdots p.11}
\line{\S\ExtrAlmPerSt.  Extremal almost periodic states\tocdots p.13}
\line{\S\FreeProdLemmas.  Some technical lemmas\tocdots p.16}
\line{\S\FPwrtMAP.  Free products of factors with respect to extremal
 almost periodic states.\tocdots p.19}
\line{\S\FPfd.  Free products of finite dimensional algebras\tocdots p.31}
\line{\S\Fullness.  Fullness of free product factors\tocdots p.44}
\line{\S\Questions.  Questions\tocdots p.45}
\line{References \tocdots p.46}
}

\head \FreeSubcompl.  Free Subcomplementation. \endhead

In~\cite{\DykemaZZFreeDim}, the free product of arbitrary finite dimensional or
AFD von Neumann algebras with respect to faithful traces was determined, and
found to be
an interpolated free group factor~(\cite{\RadulescuZZAmalgSubfactors},
\cite{\DykemaZZInterp}) possibly direct sum a finite dimensional algebra.
The free product of finite dimensional algebras was determined using the
technique 
that is now called {\it free etymology}, and passing to AFD algebras was
permitted by the use of the notion of a {\it standard embedding} of
interpolated free group factors.
The definition of a standard embedding of interpolated free group factors,
$L(\freeF_s)\hookrightarrow L(\freeF_t)$ is a generalization of the sort of
embedding $L(\freeF_n)\hookrightarrow L(\freeF_m)$, $n<m$ or $n=m=\infty$,
obtained by mapping the $n$ free generators of $\freeF_n$ to $n$ (not all)
of the $m$ free generators of $\freeF_m$.

The characterization of standard embeddings in Proposition~\FiniteStdEmbChar{}
is pleasing to the intuition and will be used in this paper.

\proclaim{Definition \FreelyCompl}\rm
Suppose $(\MvN,\phi)$ is a von Neumann algebra with normal state, and let
$1\in A\subseteq\MvN$ be a von Neumann subalgebra.
We say that $A$ is {\it freely complemented} in $(\MvN,\phi)$ if there is a von
Neumann subalgebra $B\subseteq\MvN$ such that $\MvN$ is generated by $A\cup B$
and $A$ and $B$ are free in $(\MvN,\phi)$.
To be more specific, we can say that $A$ is freely complemented by $B$.
\endproclaim

\proclaim{Proposition \FiniteStdEmbChar}
Suppose $\pi:L(\freeF_s)\hookrightarrow L(\freeF_t)$ is an injective, normal,
trace--preserving 
$*$--homomorphism.
Then $\pi$ is a standard embedding if and only if $\pi(L(\freeF_s))$ is freely
complemented in $L(\freeF_t)$ by an algebra isomorphic to
$L(\freeF_r)\oplus\Cpx$ or simply to $L(\freeF_r)$, some $1<r\le\infty$.
\endproclaim
\demo{Proof}
Let $A=\pi(L(\freeF_s))$ be freely complemented in $L(\freeF_t)$ by
$B\cong L(\freeF_r)\oplus\Cpx$ or $\cong L(\freeF_r)$.
Then by~4.4 and~1.2 of~\cite{\DykemaZZFreeDim}, $\pi$ is a standard embedding.
Since for any two standard embeddings
$\pi_1:L(\freeF_s)\hookrightarrow L(\freeF_t)$ and
$\pi_2:L(\freeF_s)\hookrightarrow L(\freeF_t)$, there is an automorphism
$\alpha$ of $L(\freeF_t)$ such that $\pi_1=\pi_2\circ\alpha$, (this is easily
seen from the definition of a standard embedding~\cite{\DykemaZZFreeDim}), and
since for any $s$ and $t$ we can find $B=L(\freeF_r)\oplus\Cpx$ and a trace on
it such that $L(\freeF_s)*B\cong L(\freeF_t)$, it follows that every standard
embedding is of this form.
\QED

The above characterization of standard embeddings allows one to easily prove
the following useful facts about standard embeddings, some of which first
appeared in~\cite{\DykemaZZFreeDim}.
\proclaim{Proposition \StdEmbFacts}
\roster
\item"(i)" The composition of standard embeddings is a standard embedding.
\item"(ii)" Suppose $A_n=L(\freeF_{t_n})$, $1<t_n\le\infty$
and $\pi_n:A_n\rightarrow A_{n+1}$ is
a standard embedding ($n\ge1$).
Consider the inductive limit von Neumann algebra (taken with respect to the
unique traces on $A_n$), $A=\varinjlim(A_n,\pi_n)$ with
the associated inclusions $\psi_n:A_n\hookrightarrow A$.
Then $A=L(\freeF_t)$ where $t=\lim_{n\rightarrow\infty}t_n$ and each $\psi_n$
is a standard embedding.
\endroster
\endproclaim

We now turn to an extension of the notion of a freely complemented subalgebra.

\proclaim{Definition \FreelySubcompl}\rm
Suppose $(\MvN,\phi)$ is a von Neumann algebra with normal state, and let
$1\in A\subseteq\MvN$ be a von Neumann subalgebra.
Let $A_\phi$ denote the centralizer of $\phi\restrict_A$.
We say that $A$ is {\it freely subcomplemented} in $(\MvN,\phi)$ if, for some
index set $I$, there is for each $\iota\in I$ a
self--adjoint projection $h_\iota\in A_\phi$ and a von Neumann subalgebra
$h_\iota\in B_\iota\subseteq h_\iota\MvN h_\iota$ such that
\roster
\item"(i)" $\MvN=W^*(A\cup\bigcup_{\iota\in I}B_\iota)$
\item"(ii)" $B_\iota$ and $h_\iota\MvN_{I\backslash\{\iota\}}h_\iota$ are free
in
$(h_\iota\MvN h_\iota,\phi(h_\iota)^{-1}\phi\restrict_{h_\iota\MvN h_\iota})$,
where
$\MvN_{I\backslash\{\iota_0\}}\eqdef
 W^*(A\cup\bigcup_{\iota\in I\backslash\{\iota_0\}}B_\iota)$.
\endroster
To be more precise, we will say that $A$ is freely subcomplemented in
$(\MvN,\phi)$ by $(B_\iota)_{\iota\in I}$ and with associated projections
$(h_\iota)_{\iota\in I}$.
\endproclaim

The condition~(ii)' of the following lemma will be the most useful one for
proving properties of free subcomplementation.

\proclaim{Lemma \FreelySubcomplWords}
Suppose $(\MvN,\phi)$ , $A$, $h_\iota$ and $B_\iota$ are as in
Definition~\FreelySubcompl, except without
assuming~(ii) holds and suppose also
that the GNS representation of $\MvN$ associated to $\phi$ is faithful.
Let $\Omega$ be the set of reduced words
$a_1a_2\cdots a_n\in\Lambdao(A,\bigcup_{\iota\in I}B_\iota\oup)$ such that
whenever $2\le j\le n-1$, $a_j\in A$, and $a_{j-1},a_{j+1}\in B_\iota\oup$ for
some $\iota\in I$,
then $a_j\in(h_\iota Ah_\iota)\oup$.
Then $\lspan(\Omega\cup\{1\})$ is a s.o.--dense $*$--subalgebra of $\MvN$.
Consider the condition
\roster
\item"(ii)'" $\phi(z)=0$ $\forall z\in\Omega\backslash A$.
\endroster
Then (ii)'$\implies$(ii).
Moreover, (ii)'~implies the existence of a s.o.--continuous projection of
norm~$1$, $E:\MvN\rightarrow A$, such that $E(z)=0$ $\forall
z\in\Omega\backslash A$.
Conversely, if each $h_\iota$ has central support in $A_\phi$ equal to~$1$
then (ii)$\implies$(ii)'.
\endproclaim
\demo{Proof}
One easily sees that $\lspan(\Omega\cup\{1\})$ is the $*$--algebra generated by
$A\cup\bigcup_{\iota\in I}B_\iota$, hence is s.o.--dense in $\MvN$.
For $F\subseteq I$, let $\Omega_F$ be the set of all words $z\in\Omega$ such
that whenever $z$ has a letter coming from $B_\iota\oup$, $\iota\in I$, then
$\iota\in F$.
Let us show that~(ii)' implies the existence of the projection
$E:\MvN\rightarrow A$.
We will realize this projection by compressing to $L^2(A,\phi)$ in the Hilbert
space $L^2(\MvN,\phi)$.
The defining mapping $\MvN\rightarrow L^2(\MvN,\phi)$ will be denoted
$x\mapsto\xh$ and we let $\xi\eqdef\onehat$.
\proclaim{Lemma \TwoProjections}
If $\PvN$ is a von Neumann algebra with normal state $\phi$ and if
$h,k\in\PvN_\phi$ are self--adjoint projections in the centralizer of $\phi$
such that $\phi(hk)=0$, then $\phi(hak)=0$ $\forall a\in\PvN$.
\endproclaim
\demo{Proof}
$0=\phi(hk)=\phi(khk).$
The left support of $kh$ is the support of $khk$ which is $k-k\wedge(1-h)$.
Since $khk\ge0$ and $\phi(khk)=0$ we have $\phi(k-k\wedge(1-h))=0$.
Thus
$$ \phi(hak)=\phi(khak)=\phi((k-k\wedge(1-h))ha(k-k\wedge(1-h)))=0. $$
Hence Lemma~\TwoProjections{} is proved.
\QED

\noindent Continuing with the proof of Lemma~\FreelySubcomplWords,
using Lemma~\TwoProjections{} we see for $\iota_1,\iota_2\in I$ that
$$ \alignedat 2
h_{\iota_1}Ah_{\iota_2}&\subseteq\ker\phi
 &&\text{ if }\phi(h_{\iota_1}h_{\iota_2})=0 \\
h_{\iota_1}Ah_{\iota_2}&=(h_{\iota_1}Ah_{\iota_2})\cap\ker\phi
  +\Cpx h_{\iota_1}h_{\iota_2}
 &&\text{ if }\phi(h_{\iota_1}h_{\iota_2})\neq0
\endalignedat \tag{\twoprojections} $$
Using~(ii)', equation~(\twoprojections) and that $\lspan(\Omega\cup\{1\})$ is
dense in
$\MvN$, we find that  $L^2(\MvN,\phi)$ can be characterized as follows.
Let $\alpha\not\in I$,
$\Hilo_\alpha\eqdef L^2(A,\phi\restrict_A)\ominus\Cpx\onehat$, for $\iota\in I$
let
$\Hilo_\iota\eqdef L^2(B_\iota,\phi\restrict_{B_\iota})\ominus\Cpx\hh_\iota$,
and let
$$ \Hil\eqdef\Cpx\xi\oplus\bigoplus
 \Sb n\ge1\\j_k\in I\cup\{\alpha\}\\j_1\neq j_2\neq\cdots\neq j_n\endSb
 \Hilo_{j_1}\otimes\Hilo_{j_2}\otimes\cdots\otimes\Hilo_{j_n}. $$
Then $L^2(\MvN,\phi)$ is the closed subspace of $\Hil$ that is spanned by $\xi$
together with all vectors
$$ \zeta_1\otimes\cdots\otimes\zeta_n
 \in\Hilo_{j_1}\otimes\cdots\otimes\Hilo_{j_n},\quad j_1\neq\cdots\neq j_n $$
such that
$$ \align
\text{if }2\le k\le n-1,\,j_k=\alpha,
 \text{ (hence }j_{k-1},j_{k+1}\in I\text{),}
 &\text{ then }\zeta_k\in\overline{\{\ah\mid a\in h_{j_{k-1}}Ah_{j_{k+1}}
  \cap\ker\phi\}} \\
\text{if }n\ge2\text{ and }j_1=\alpha,
 \text{ (hence }j_2\in I\text{),}
 &\text{ then }\zeta_1\in\overline{\{\ah\mid a\in Ah_{j_2}
  \cap\ker\phi\}} \\
\text{if }n\ge2\text{ and }j_n=\alpha,
 \text{ (hence }j_{n-1}\in I\text{),}
 &\text{ then }\zeta_n\in\overline{\{\ah\mid a\in h_{j_{n-1}}A
  \cap\ker\phi\}}.
\endalign $$
Let $P_A$ be the self--adjoint projection of $L^2(\MvN,\phi)$ onto
$L^2(A,\phi\restrict_A)=\Cpx\xi\oplus\Hilo_\alpha\subseteq L^2(\MvN,\phi)$.
Let $E:\Bof(L^2(\MvN,\phi))\rightarrow\Bof(L^2(A,\phi\restrict_A))$ be defined
by $E(x)=P_AxP_A$.
Then $E$ is a projection of norm $1$ and is s.o--continuous.
We easily see that
$$ \align
E(a)=a&\quad\forall a\in A \\
E(z)=0&\quad\forall z\in\Omega\backslash A.
\endalign $$
Thus, by s.o.--continuity of $E$, we see that $E$ projects $\MvN$ onto $A$.

Now using the projection $E$, we will show (ii)'$\implies$(ii).
Let $\iota\in I$.
We will show that $B_\iota$ and $h_\iota\MvN_{I\backslash\{\iota\}}h_\iota$ are
free in $h_\iota\MvN h_\iota$.
Clearly $\lspan(\Omega_{I\backslash\{\iota\}}\cup\{1\})$
is a dense $*$--subalgebra of
$\MvN_{I\backslash\{\iota\}}$.
Using the projection $E:\MvN\rightarrow A$, we see that each element of
$(h_\iota\MvN_{I\backslash\{\iota\}}h_\iota)\oup$ is the s.o.--limit of a
bounded net in
$\lspan(h_\iota(\Omega_{I\backslash\{\iota\}}\backslash A)h_\iota
\cup(h_\iota Ah_\iota)\oup)$,
hence it will suffice to show $\phi(y)=0$ whenever
$$ y\in\Lambdao(h_\iota(\Omega_{I\backslash\{\iota\}}\backslash A)h_\iota
\cup(h_\iota Ah_\iota)\oup,B_\iota\oup). $$
But by regrouping we see that
$y\in(\Omega\backslash A)\cup(h_\iota A h_\iota)\oup$, hence $\phi(y)=0$
by~(ii)'.

Assume that each $h_\iota$ has central support in $A_\phi$ equal to~$1$
and let us show (ii)$\implies$(ii)'.
Take any $z=a_1a_2\cdots a_n\in\Omega_F\backslash A$, where $F$ is a finite
subset of $I$.
It will suffice to show $\phi(z)=0$, and we will proceed by induction on $|F|$,
proving the $|F|=1$ case and the inductive step simultaneously.
Let $\iota\in F$ be such that $z$ has a letter from $B_\iota\oup$.
Now we may assume $n\ge3$ and $a_1,a_n\in A$.
By hypothesis, there are partial isometries $v_0,\ldots,v_k\in A_\phi$ such
that $v_j^*v_j\le h_\iota$ $\forall1\le j\le k$ and $\sum_{j=1}^kv_jv_j^*=1$.
Now $\phi(z)=\sum_{j_1,j_2=0}^k\phi(v_{j_1}v_{j_1}^*zv_{j_2}v_{j_2}^*)
=\sum_{j=0}^k\phi(v^*_jzv_j)$, so we may assume without loss of
generality that $a_1\in h_\iota A$ and $a_n\in Ah_\iota$.
Regrouping in $z$ by singling out all letters from $B_\iota\oup$ and grouping
all other strings of neighboring letters together, we see that
$z\in\Lambdao(h_\iota\Omega_{F\backslash\{\iota\}}h_\iota,B_\iota\oup)$ and
if a letter of $z$ is in $h_\iota\Omega_{F\backslash\{\iota\}}h_\iota\cap A$
then it
is in $(h_\iota Ah_\iota)\oup$, unless it is the first or the last letter, in
which case it is in $h_\iota Ah_\iota$.
Using if necessary $h_\iota Ah_\iota=(h_\iota Ah_\iota)\oup+\Cpx h_\iota$ to
expand the first and/or the last letter, we
see that $z$ is a linear combination of up to four elements of
$$ \Lambdao((h_\iota\Omega_{F\backslash\{\iota\}}h_\iota\backslash A)
 \cup(h_\iota Ah_\iota)\oup,B_\iota\oup). $$
By inductive hypothesis in the case $|F|>1$ and by the fact that
$\Omega_{F\backslash\{\iota\}}\subseteq A\oup$ if $|F|=1$, we have that
$$ z\in\lspan\Lambdao((h_\iota\MvN_{I\backslash\{\iota\}}h_\iota)\oup,
 B_\iota\oup), $$
hence by~(ii), $\phi(z)=0$.
\QED

Now we show that the composition of (certain) freely subcomplemented embeddings
is a freely subcomplemented embedding.

\proclaim{Corollary \FreelySubcomplChain}
Suppose
$$ 1\in A_1\subseteq A_2\subseteq\cdots\subseteq\MvN
 \tag{\freelysubcomplchain} $$
is a finite or infinite chain of inclusions of von Neumann algebras, $\phi$ is
a normal state on $\MvN$ and
$$ \MvN=\cases A_n&\text{if the chain~(\freelysubcomplchain) is finite of
length }n \\
\overline{\bigcup_{n\ge1}A_n}
 &\text{if the chain~(\freelysubcomplchain) is infinite}
\endcases $$
Suppose for all $j$ that $A_j$ is freely subcomplemented in
$(A_{j+1},\phi\restrict_{A_{j+1}})$ by $(B_\iota)_{\iota\in I_j}$ with
associated projections $(h_\iota)_{\iota\in I_j}$ such that each $h_\iota$ has
central support in $(A_j)_\phi$ equal to~$1$ and is
Murray--von~Neumann equivalent in $(A_j)_\phi$ to a projection in $A_1$.
Then $A_1$ is freely subcomplemented in $\MvN$ by (algebras isomorphic to)
$(B_\iota)_{\iota\in I}$, with $I=\bigcup_j I_j$.
\endproclaim

The above result is a special case of the next one.

\proclaim{Corollary \FreelySubcomplWell}
Let $J$ be a set and $\le$ a well--ordering on $J$.
Suppose $\forall j\in J$ $A_j$ is a von Neumann algebra with normal state
$\phi_j$ and suppose there are compatible, $\phi_j$--preserving inclusions
$A_{j_1}\subseteq A_{j_2}$ $\forall j_1\le j_2$.
Let $(\MvN,\phi)=\varinjlim_{j\in J}(A_j,\phi_j)$ be the inductive limit von
Neumann algebra.
Let $j_0\in J$ denote the smallest element and suppose
$\forall j\in J\backslash\{j_0\}$ that $\MvN_{\lneqq j}\eqdef
W^*(\{1\}\cup\bigcup_{j_1\lneqq j}A_{j_1})$ is 
freely subcomplemented in $W^*(\MvN_{\lneqq j}\cup A_j)$ (with respect to
the restriction of $\phi$ to $W^*(\MvN_{\lneqq j}\cup A_j)$)
by $(B_\iota)_{\iota\in I_j}$ with associated
projections $(h_\iota)_{\iota\in I_j}$, such that each $h_\iota$ has central
support in $(A_j)_\phi$  equal to $1$ and is Murray--von Neumann equivalent to
a projection in $A_{j_0}$.
Then $A_{j_0}$ is freely subcomplemented in $(\MvN,\phi)$.
by (algebras isomorphic to) $(B_\iota)_{\iota\in I}$ where
$I=\bigcup_{j\in J}I_j$.
\endproclaim
\demo{Proof}
Let $(h_\iota)_{\iota\in I_j}$ be the projections associated to
$(B_\iota)_{\iota\in I_j}$.
We may suppose
$h_\iota\in (A_{j_0})_\phi$ $\forall\iota\in I$.
Now it will suffice to show~(ii)' of Lemma~\FreelySubcomplWords, (with
$A=A_{j_0}$ here).
Let $z\in\Omega_F\backslash A$, where $F\subseteq I$ is any finite subset and
let $j\in J$ be the largest element
of $J$ such that $F\cap I_j$ is nonempty.
We show $\phi(z)=0$ by transfinite induction on $j$.
If $j=j_0$ then $\Omega_F\subseteq A$ and there is nothing to prove.
Otherwise, regrouping in $z$ by grouping letters from all $B_\iota\oup$ for
$\iota\not\in I_j$ and $A$ together, we see that $z$ equals an element of
$\Lambdao(\Omega_{F\backslash I_j},\bigcup_{\iota\in I_j}B_\iota\oup)$
whose every letter lying between letters from $B_\iota\oup$, $\iota\in I_j$
belongs to
$(h_\iota\Omega_{F\backslash I_j}h_\iota\backslash A)
\cup(h_\iota Ah_\iota)\oup$.
This last set is, by inductive hypothesis, a subset of
$(h_\iota\MvN_{\lneqq j}h_\iota)\oup$.
Applying Lemma~\FreelySubcomplWords(ii)' to the situation of $\MvN_{\lneqq j}$
freely subcomplemented in $A_{j}$, we have that $\phi(z)=0$.
\QED

\head \FreeAmalg.  Freeness with amalgamation over a subalgebra. \endhead

Freeness of C$^*$--algebras over a subalgebras was defined by
Voiculescu, and this notion of freeness
corresponds to his construction of free products with amalgamation
over common subalgebras.

\proclaim{Definition \FreenessAmalg (\cite{\VoiculescuZZSymetries})}\rm
Suppose $A$ is an algebra, $B\subseteq A$ a subalgebra,
$E:A\rightarrow B$ a $B$--$B$ bimodule map.
Given subalgebras $B\subseteq A_\iota\subseteq A$
($\iota\in I$),
we say that the family $(A_\iota)_{\iota\in I}$ is {\it free over} $B$ if
$E(a_1a_2\cdots a_n)=0$ whenever $a_j\in A_{\iota_j}$,
$\iota_k\neq\iota_{k+1}$ ($1\le k\le n-1$) and $E(a_j)=0$ $\forall j$.
To be emphatic, we will sometimes say the family $(A_\iota)_{\iota\in I}$ is
{\it free with amalgamation over} $B$.
\endproclaim

\proclaim{Proposition \FreenessAmalgCondExp}
Suppose $A$ is a C$^*$--algebra with state $\phi$,
$B$ is a 
C$^*$--subalgebra of $A$ such that the G.N.S\. representation associated to
$\phi\restrict_B$ is faithful on $B$,
$E:A\rightarrow B$ is a projection of norm one
(hence a conditional expectation) from $A$ onto $B$ such that
$\phi\circ E=\phi$.
Suppose $B\subseteq A_\iota\subseteq A$ are subalgebras.
Then to show that $(A_\iota)_{\iota\in I}$ is free with amalgamation over $B$,
it suffices to show that $\phi(a_1a_2\cdots a_n)=0$ whenever
$a_j\in A_{\iota_j}$, $\iota_1\neq\iota_2\neq\cdots\neq\iota_n$ and $E(a_j)=0$
$\forall j$.
\endproclaim
\demo{Proof}
{}From the fact that $E(a_j)=0$ we have, for any $b_1,b_2\in B$, that
$E(a_nb_2)=E(a_n)b=0=E(b_1a_1)$ and hence, if the condition
of the proposition holds, we have
$0=\phi((b_1a_1)a_2\cdots(a_nb_2))
=\phi\circ E((b_1a_1)a_2\cdots(a_nb_2))
=\phi(b_1E(a_1\cdots a_n)b_2)$
Letting $(\pi,\Hil,\xi)=\text{GNS}(B,\phi_B)$, it follows that
$0=\phi(b_1E(a_1\cdots a_n)b_2)
=\langle\pi(E(a_1\cdots a_n))\bh_2,(b_1^*)\hat{\;}\rangle$ so
$\pi(E(a_1\cdots a_n)=0$ and hence $E(a_1\cdots a_n)=0$.
Thus, the condition of the proposition implies freeness over $B$.
\QED

\proclaim{Proposition \FreelySubcomplAmalg}
Suppose $\MvN$ is a von Neumann algebra with normal state $\phi$,
$1\in A\subseteq\MvN$ a von Neumann subalgebra, $E=E^\MvN_A:\MvN\rightarrow A$
a s.o.--continuous projection of norm~$1$ such that $\phi\circ E=\phi$.
Suppose $A\subseteq\NvN_j\subseteq\MvN$ is a von Neumann subalgebra ($j\in J$)
such that $\bigcup_{j\in J}\NvN_j$ generates $\MvN$ and $(\NvN_j)_{j\in J}$ is
free over $A$ in $(\MvN,E)$.
Suppose for each $j\in J$ $A$ is freely subcomplemented in $\NvN_j$ by
$(B_\iota)_{\iota\in I_j}$ with associated projections
$(h_\iota)_{\iota\in I_j}$ such that the central support of each $h_\iota$ in
the centralizer of $\phi\restrict_A$ is~$1$.
Let $I=\bigcup_{j\in J}I_j$ be a disjoint union.
Then $A$ is freely subcomplemented in $\MvN$ by $(B_\iota)_{\iota\in I}$
with associated projections $(h_\iota)_{\iota\in I}$
and 
$\NvN_j$ is freely subcomplemented in $\MvN$ by
$(B_\iota)_{\iota\in I\backslash I_j}$ with associated projections
$(h_\iota)_{\iota\in I\backslash I_j}$.
\endproclaim
\demo{Proof}
To prove that both $A$ and $\NvN_j$ are freely subcomplemented in $\MvN$,
using Definition~\FreelySubcompl(ii) and
Lemma~\FreelySubcomplWords(ii)' we
see that it will suffice to show $\phi(z)=0$ $\forall z\in\Omega\backslash A$,
with $\Omega$ as in Lemma~\FreelySubcomplWords.
However, by regrouping, we see that
$z\in\Lambdao((\Omega_{I_j}\backslash A)_{j\in J})$.
Using that $\phi(y)=0$ $\forall y\in\Omega_{I_j}\backslash A$, we see that
$\Omega_{I_j}\backslash A\subseteq\NvN_j\cap\ker E^\MvN_A$, hence $\phi(z)=0$.
\QED

\head \CondExp.  Conditional Expectations. \endhead

Recall Tomiyama's famous result, \cite{\Tomiyama}, that a projection of norm
$1$, $E$, 
from a C$^*$--algebra $A$ onto a C$^*$--subalgebra $B$ is necessarily positive
and a
conditional expectation, {\it i.e\.} $E(b_1ab_2)=b_1E(a)b_2$ for all $a\in A$,
$b_1,b_2\in B$.
The following two lemmas, (which we do not claim are new), form a
characterization of conditional expectations in terms of projections on the
Hilbert space of a G.N.S\. representation.

\proclaim{Lemma \CondExpFromDense}
Let $A$ be a von Neumann algebra, $B\subseteq A$ a von Neumann subalgebra and
$\phi$ a normal, faithful state on $A$.
Let $(\pi,\Hil,\xi)=\text{GNS}(A,\phi)$, so $\Hil$ contains
$\{\ah\mid a\in A\}$ as a dense subspace with inner product
$\langle\ah_1,\ah_2\rangle=\phi(a_2^*a_1)$.
Let $\Hil_B=\overline{\{\bh\mid b\in B\}}\subseteq\Hil$, and let
$P:\Hil\rightarrow\Hil_B$ be the orthogonal projection.
Suppose that $\Afr$ is a $\sigma$--weakly dense subset of $A$, and that
$P\pi(a)\restrict_{\Hil_B}\in P\pi(B)\restrict_{\Hil_B}$ for every $a\in \Afr$.
Then there is a projection of norm $1$, $E:A\rightarrow B$, such that
$\phi\circ E=\phi$.
Moreover, $E$ is continuous from the s.o.--topology on $A$ (acting on $\Hil$)
to the s.o.--topology on $B$ (acting on $\Hil_B$).
\endproclaim
\demo{Proof}
We have that $\pi$ is a normal, faithful representation of $A$ on $\Hil$ and
$\pi_B:B\rightarrow\Bof(\Hil_B)$ given by $\pi_B(b)=P\pi(b)\restrict_{\Hil_B}$
is a normal, faithful representation of $B$ on $\Hil_B$, and in fact
$(\pi_B,\Hil_B,\xi)=\text{GNS}(B,\phi\restrict_B)$.
Hence, taking limits in the weak--operator topology we have that
$P\pi(A)\restrict_{\Hil_B}\subseteq \pi_B(B)$ and we may
define $E:A\rightarrow B$ by 
$$ P\pi(a)\restrict_{\Hil_B}=P\pi(E(a))\restrict_{\Hil_B}. $$
Clearly $E(b)=b$ if $b\in B$ and $\nm{E(a)}\le\nm a$ $\forall a\in A$.
Also, $\phi(E(a))=\langle\pi_B(E(a))\xi,\xi\rangle=\langle
P\pi(a)P\xi,xi\rangle=\langle\pi(a)\xi,\xi\rangle=\phi(a)$,
so $\phi\circ E=\phi$.
The s.o.--continuity of $E$ is clear.
\QED

\proclaim{Lemma \CondExpImplByProj}
Suppose $(A,\phi)$ is a C$^*$--algebra with faithful state.
Suppose $B\subseteq A$ is a C$^*$--subalgebra with projection of norm $1$,
$E:A\rightarrow B$, such that $\phi\circ E=\phi$.
Let $(\pi,\Hil,\xi)=\text{GNS}(A,\phi)$, let
$\Hil_B=\overline{\{\bh\mid b\in B\}}$ and let $P:\Hil\rightarrow\Hil_B$ be the
orthogonal projection, with notation as in Lemma~\CondExpFromDense.
Then 
$$ \align
(E(a))\hat{\,}&=P\ah, \tag{\EhatIsP} \\
P\pi(a)\restrict_{\Hil_B}&=P\pi(E(a))\restrict_{\Hil_B}\;\forall a\in A.
\tag{\EImplProj}
\endalign $$
\endproclaim
\demo{Proof}
The identity~(\EhatIsP) follows because
$$ \langle E(a)\hat{\,},\bh\rangle=\phi(b^*E(a))=\phi(E(b^*a))=\phi(b^*a)
=\langle\ah,\bh\rangle=\langle\ah,P\bh\rangle=\langle P\ah,\bh\rangle
\;\forall b\in B $$
and then~(\EImplProj) holds because
$$ P(\pi(a)\bh)=P((ab)\hat{\,})=E(ab)\hat{\,}=(E(a)b)\hat{\,}=E(a)\bh. $$
\QED

\proclaim{Proposition \FreeProdOfCondExp}
Suppose $(A_\iota,\phi_\iota)$ are von Neumann algebras with normal, faithful
states, for $\iota$ in some index set $I$.
Suppose $1\in B_\iota\subseteq A_\iota$ are von Neumann subalgebras with
projections of norm $1$, $E_\iota:A_\iota\rightarrow B_\iota$, such that
$\phi_\iota\circ E_\iota=\phi_\iota$.
Let $(\MvN,\phi)=\freeprodi(A_\iota,\phi_\iota)$ and let
$\NvN=W^*(\bigcup_{\iota\in I}B_\iota)\subseteq\MvN$.
Then there is a projection of norm $1$, $E:\MvN\rightarrow\NvN$ such that
$\phi\circ E=\phi$
and $E(a)=E_\iota(a)$ whenever $a\in A_\iota$.
\endproclaim
\demo{Proof}
Let $(\pi_\iota,\Hil_\iota,\xi_\iota)=\text{GNS}(A_\iota,\phi_\iota)$,
$(\pi,\Hil,\xi)=\text{GNS}(\MvN,\phi)$.
{}From the construction of the free product we have that
$$ \Hil=\Cpx\xi\oplus
\bigoplus\Sb n\ge1 \\ \iota_1\neq\iota_2\neq\cdots\neq\iota_n \endSb
\Hilo_{\iota_1}\otimes\cdots\otimes\Hilo_{\iota_n}, $$
where $\Hilo_\iota=\Hil_\iota\ominus\Cpx\xi_\iota$.
By Lemma~\CondExpImplByProj{} there is a projection
$P_\iota:\Hil_\iota\rightarrow\Hil_{B_\iota}$ implementing the conditional
expectation $E_\iota$.
If $\Omega=\Lambdao((B_\iota\oup)_{\iota\in I})$ then
$\Afr=\lspan(\Omega\cup\{1\})$
is the span of words in the $B_\iota$ and is a s.o.--dense $*$--subalgebra of
$\NvN$. 
Thus it is clear that
$$ \overline{\{\ah\mid a\in\NvN\}}=\overline{\{\ah\mid a\in\Afr\}}
=\Hil_\NvN=\Cpx\xi\oplus
\bigoplus\Sb n\ge1 \\ \iota_1\neq\iota_2\neq\cdots\neq\iota_n \endSb
\Hilo_{B_{\iota_1}}\otimes\cdots\otimes\Hilo_{B_{\iota_n}}\subseteq\Hil. $$
Let $P:\Hil\rightarrow\Hil_\NvN$ be the orthogonal projection.
Then for a word $x=a_1a_2\cdots a_n\in\Omega$, $a_j\in A_{\iota_j}$,
$\iota_j\neq\iota_{j+1}$, we have that
$$ P\pi(x)\restrict_{\Hil_\NvN}
=P\pi(E_{\iota_1}(a_1)E_{\iota_2}(a_2)\cdots
E_{\iota_n}(a_n))\restrict_{\Hil_\NvN}, $$
hence the same holds for all $x\in\Omega$.
Thus, by Lemma~\CondExpFromDense, there is a conditional expectation $E$ with
the desired properties.
\QED

\proclaim{Corollary \FreeFactorCondExp}
Suppose $(A_\iota,\phi_\iota)$ are von Neumann algebras with normal, faithful
states, $(\iota=1,2)$, and let $(\MvN,\phi)=(A_1,\phi_1)*(A_2,\phi_2)$.
The there is a projection, $E$, of norm $1$ from $\MvN$ onto $A_1$ such that
$\phi_1\circ E=\phi$.
Moreover, $E$ is continuous from the strong--operator topology on $\MvN$ (in
its 
GNS representation corresponding to $\phi$) to the strong--operator topology on
$A_1$ (in its GNS representation corresponding to $\phi_1$).
\endproclaim

\head \ExtrAlmPerSt.  Extremal almost periodic states. \endhead

Let us first recall a few definitions and theorems of Connes'.
\proclaim{Definition \AlmPer{} (\cite{\ConnesZZAlmPer})} \rm
Let $(\MvN,\phi)$ consist of a von Neumann algebra, $\MvN$, with normal,
faithful (n.f\.) state or semi--finite weight, $\phi$.
Then $\phi$ is said to be {\it almost periodic} if the modular operator,
$\Delta_\phi$, has pure point spectrum.
If, moreover, $\MvN$ is full, one defines $\Sd(\MvN)$ to be the intersection of
the point spectra of the modular operators $\Delta_\phi$ as $\phi$ ranges over
all n.f\. almost periodic weights.
\endproclaim

Almost periodic states behave well with respect to free products.
We will be able to determine explicitly many free products that are taken with
respect to almost periodic states.
The following basic observation follows easily from Voiculsecu's description of
the conjugation
operator (see Lemma~1.8 of~\cite{\VoiculescuZZSymetries} or the discussion
in~\S1 of~\cite{\DykemaZZTM}).

\proclaim{Proposition \AlmPerFreeProd}
For each $i\in I$, given a von Neumann algebra, $A_i$, with n.f\.
almost periodic state, $\phi_i$, whose modular operator has point spectrum
$\Gamma_i$, consider the free product
$(\MvN,\phi)=\freeprodi(A_i,\phi_i)$.
Then $\phi$ is an almost periodic state and the point spectrum of its modular
operator equals the subgroup of $\Real_+^*$ generated by
$\bigcup_{i\in I}\Gamma_i$.
\endproclaim

\proclaim{Remark \SpectralTh}\rm
It follows from Arveson's and Connes' spectral
theory (\cite{\Arveson}, \cite{\ConnesZZThesis})
that whenever $\phi$ is
an almost periodic n.f\. state or weight on any von Neumann algebra $\MvN$ one
has a decomposition
$$ \MvN=\bigoplus_{\gamma\in\Gamma}\MvN_\phi(\gamma), \tag{\spectraldecomp} $$
where $\Gamma$ is the point spectrum of $\Delta_\phi$ and where
$$
\MvN_\phi(\gamma)=\{a\in\MvN\mid\phi(xa)=\gamma\phi(ax)\;\forall\,x\in\MvN\}.
$$
(See~\S1 of~\cite{\DykemaZZAlmPer}.)
One also has that each $\MvN_\phi(\gamma)$ is w.o.--closed, that
$\MvN_\phi(\gamma_1)\MvN_\phi(\gamma_2)\subseteq\MvN_\phi(\gamma_1\gamma_2)$
and that $\MvN(\gamma)^*=\MvN(\gamma^{-1})$.
Thus if $\gamma>1$ and if $\MvN_\phi(\gamma)$ contains an isometry $v$, then
$\phi(vv^*)=\gamma$ and $\MvN_\phi(\gamma)=\MvN_\phi v$, where
$\MvN_\phi$ denotes the {\it centralizer} of $\phi$, which is by definition
$\MvN_\phi=\MvN_\phi(1)$.
\endproclaim

\proclaim{Definition \Msg}\rm
Suppose $\phi$ is an almost periodic weight on a von Neumann algebra $\MvN$.
Let $\Gamma$ be the subgroup of $\Real_+^*$ generated by the point spectrum of
the modular operator $\Delta_\phi$.
Suppose $\Gamma'$ is a subgroup of $\Gamma$.
Then define $\MvN_\phi(\Gamma')$ to be the weak closure of the linear span of
$\bigcup_{\gamma\in\Gamma'}\MvN_\phi(\gamma)$.
\endproclaim

\proclaim{Lemma \CondExpAlmPer}
In the situation of Definition~\Msg,
there is a projection of norm $1$, $E$, from $\MvN$ onto $\MvN_\phi(\Gamma')$
such that $\phi\circ E=\phi$ and
that is continuous from the strong--operator topology on $\MvN$ (in its GNS
representation from $\phi$) to the strong--operator topology on
$\MvN_\phi(\Gamma')$ (in its GNS representation from $\phi$).
We have for $a\in\MvN_\phi(\gamma)$ that
$$ E(a)
=\cases a&\text{ if }\gamma\in\Gamma' \\ 0&\text{ otherwise.} \endcases $$
\endproclaim
\demo{Proof}
Let $\NvN$ denote $\MvN_\phi(\Gamma')$.
Let $(\pi,\Hil,\xi)=\text{GNS}(\MvN,\phi)$.
Since $\lspan(\bigcup_{\gamma\in\Gamma}\MvN_\phi(\gamma))$ is s.o.--dense in
$\MvN$, we get an orthogonal decomposition of the Hilbert space,
$\Hil=\bigoplus_{\gamma\in\Gamma}\Hil_\phi(\gamma)$, where
$\Hil_\phi(\gamma)=\overline{\{\ah\mid a\in\MvN_\phi(a)\}}$, and
one sees that
$\Hil_{\NvN}=\bigoplus_{\gamma\in\Gamma'}\Hil_\phi(\gamma)$.
If $P:\Hil\rightarrow\Hil_{\NvN}$ is the orthogonal projection
then for $a\in\MvN_\phi(\gamma)$
$$ P\pi(a)\restrict_{\Hil_\NvN}
=\cases P\pi(a)\restrict_{\Hil_\NvN}&\text{ if }\gamma\in\Gamma' \\
0&\text{ otherwise.} \endcases $$
Hence Lemma~\CondExpFromDense{} applies to give $E$ with the desired
properties.
\QED

\proclaim{Proposition \SdProp{} (\cite{\ConnesZZAlmPer})}
Suppose $\MvN$ is a full factor.
Then $\Sd(\MvN)$ is a multiplicative subgroup of the positive reals,
$\Real_+^*$.
If $Sd(\MvN)\neq\Real_+^*$ and if $\phi$ is an almost periodic weight on
$\MvN$, then $\text{(the point spectrum of }\Delta_\phi)=\Sd(\MvN)$ if and only
if $\MvN_\phi$ is a factor.
\endproclaim

The following is in analogy with Connes' result.

\proclaim{Definition \DefExtrAPSt}\rm
An {\it extremal almost periodic} state (or weight) on a von Neumann algebra
$\MvN$ is a normal, faithful, almost periodic state (or weight), $\phi$, whose
centralizer, $\MvN_\phi$, is a factor.
We define $\Sd(\phi)$ to be the point spectrum of the modular operator
$\Delta_\phi$.
\endproclaim

\proclaim{Lemma \MustBeFactor}
Suppose a von Neumann algebra $\MvN$ has an extremal almost periodic state
$\phi$.
Then $\MvN$ is a factor and
\roster
\item"(i)" is type II$_1$ if $\Sd(\phi)=\{1\}$;
\item"(ii)" is type III$_\lambda$ if $\Sd(\phi)
=\{\lambda^n\mid n\in\Integers\}$, $0<\lambda<1$;
\item"(iii)" is type III$_1$ if $\Sd(\phi)$ is dense in $\Real_+^*$.
\endroster
\endproclaim
\demo{Proof}
The center of $\MvN$ lies in the centralizer of $\phi$, which is a factor, so
$\MvN$ is a factor.
The type classification~(i), (ii), and~(iii) of Connes' follows
from Corollary~3.2.7. of~\cite{\ConnesZZThesis}.
\QED

\proclaim{Lemma \ExtrAPStDecomp}
Suppose $\MvN$ is a factor with extremal almost periodic state $\phi$.
Then
\roster
\item"(i)" for every $\gamma\in\Sd(\phi)$, $\gamma>1$, there is an isometry
$v_\gamma$ in the spectral subspace $\MvN_\phi(\gamma)$;
\item"(ii)" $\Sd(\phi)$ is a subgroup of the multiplicative group $\Real_+^*$
\endroster
\endproclaim
\demo{Proof}
Since $\gamma$ is in the point spectrum of $\Delta_\phi$, one sees that
$\MvN_\phi(\gamma)$ is nonempty.
Taking the polar decomposition of any element of $\MvN_\phi(\gamma)$, one shows
that $\MvN_\phi(\gamma)$ contains a partial isometry, $v'$.
Then $\phi(v'v^{\prime*})=\gamma\phi(v^{\prime*}v')$.
Since $\MvN_\phi$ is a factor and $\gamma>1$, there are partial isometries
$x_j$ and $y_j$ in 
$\MvN_\phi$ such that $v=\sum_{j=1}^nx_jv'y_j\in\MvN_\phi(\gamma)$ is an
isometry. 
This proves~(i).
Part~(ii) now follows easily from~(i) and the fact that
$h_1\MvN_\phi h_2\neq\{0\}$ if $h_1,h_2$ are nonzero projections in
$\MvN_\phi$.
\enddemo

\proclaim{Lemma \ExtrAPStCutDowns}
\roster
\item"(i)" Let $\MvN$ be a type III factor having an extremal almost periodic
state $\phi$
and let $h\in\MvN_\phi$ be a self--adjoint projection.
Then $\phi_h=\phi(h)^{-1}\phi\restrict_{h\MvN h}$ is an extremal almost
periodic state on $h\MvN h$ whose centralizer is equal to $h\MvN_\phi h$, and
$\Sd(\phi_h)=\Sd(\phi)$. 
\item"(ii)" Conversely, if $\MvN$ is a von Neumann algebra with almost periodic
state $\phi$ and if $h\in\MvN_\phi$ is a self--adjoint projection whose
central support in $\MvN_\phi$ is $1$ and such that $h\MvN_\phi h$ is a factor,
then $\MvN$ is a factor and $\phi$ is extremal almost periodic.
\endroster
\endproclaim
\demo{Proof}
Part~(i) is clear from Lemma~\ExtrAPStDecomp{} and Remark~\SpectralTh{} and
part~(ii) is tautologically true, because $\MvN_\phi$ must be a factor.
\QED

In the remainder of this section, we prove that an extremal almost
periodic state, $\phi$, is truly extremal, in the sense that if $\psi$ is an
almost periodic state that commutes with $\phi$ in the sense
of~\cite{\PedersenTakesaki}, then the point spectrum of $\Delta_\psi$ contains
the point spectrum of $\Delta_\phi$.

Given faithful, normal states, $\phi$ and $\psi$, on a von
Neumann algebra $\MvN$, by Pedersen and Takesaki's Radon-Nikodym Theorem for
von Neumann algebras, \cite{\PedersenTakesaki},
$\psi$ commutes with $\phi$
if and only if there is a self--adjoint, nonnegative
operator, $h$, affiliated with the centralizer of $\phi$, $\MvN_\phi$, such
that $\psi=\phi(\cdot h)$.
\proclaim{Lemma \WhenAlmPer}
Let $\phi$ be a faithful, extremal almost periodic, normal state on $\MvN$ and
let $\psi$ be a faithful, normal state on $\MvN$ commuting with $\phi$, so
$\psi=\phi(\cdot h)$ as described above.
Then $\psi$ is almost periodic if
and only if $h$ has purely atomic spectral measure, {\rm i.e\.} if and only if
$h=\sum_{k\in I} t_kE_k$ for $t_k\in\Reals^*_+$ and $(E_k)_{k\in I}$ an
orthogonal family of projections in $\MvN_\phi$ whose sum is $1$.
\endproclaim
\demo{Proof}
Let $E(\cdot)$ be the spectral measure of $h$, so
$h=\int_{\Reals_+}tE(\dif t)$.
Let
$$ \align
\eta_\phi:&\MvN\rightarrow L^2(\MvN,\phi) \\
\eta_\psi:&\MvN\rightarrow L^2(\MvN,\psi)
\endalign $$
be the usual embeddings, and denote the inner product in $L^2(\MvN,\phi)$ by
$\langle\cdot,\cdot\rangle_\phi$ and in $L^2(\MvN,\psi)$ by
$\langle\cdot,\cdot\rangle_\psi$.
Let $\Gamma_\phi$ denote the point spectrum of $\Delta_\phi$ and let
$$ \mfr=\bigcup
\Sb 0<t_1<t_2<\infty \\ \gamma\in\Gamma_\phi \endSb
E([t_1,t_2])\MvN_\phi(\gamma)E([t_1,t_2]), $$
so $\mfr$ is s.o.--dense in $\MvN$.
Then for $x\in\mfr\cap\MvN_\phi(\gamma)$ and $y\in\mfr$,
$$ \align
\langle\eta_\psi(x),\eta_\psi(y)\rangle_\psi&=\psi(y^*x)=\phi(y^*xh)
=\langle\eta_\phi(xh),\eta_\phi(y)\rangle_\phi \\
\langle\Delta_\psi\eta_\psi(x),\eta_\psi(y)\rangle_\psi
&=\langle\eta_\psi(y^*),\eta_\psi(x^*)\rangle_\psi
=\psi(xy^*)=\phi((hx)y^*)=\gamma\phi(y^*(hx))= \\
&=\gamma\psi(y^*(hxh^{-1}))
=\gamma\langle\eta_\psi(hxh^{-1}),\eta_\psi(y)\rangle_\psi. \tag{\psiphi}
\endalign $$
For $\gamma\in\Gamma_\phi$, let $F_\gamma$ be the
orthogonal projection from $L^2(\MvN,\psi)$ onto
$\clspan\{\eta_\psi(x)\mid x\in\MvN_\gamma(\phi)\}$.
Then $(F_\gamma)_{\gamma\in\Gamma_\phi}$ is an orthogonal family of projections
whose sum is~$1$.
Let $D_\phi=\sum_{\gamma\in\Gamma_\phi}\gamma F_\gamma$.
Let $L_h$ and $R_{h^{-1}}$ be the closures of the operators defined on $\mfr$
by $L_h\eta_\psi(x)=\eta_\psi(hx)$ and
$R_{h^{-1}}\eta_\psi(x)=\eta_\psi(xh^{-1})$.
Clearly, $L_{E(\cdot)}$ is the spectral measure for $L_h$ and
$R_{E(\sigma(\cdot))}$ is the spectral measure for $R_{h^{-1}}$, where
$\sigma(t)=t^{-1}$.
Then $D_\phi$, $L_h$ and $R_{h^{-1}}$ are commuting,
(see~\S VII.5 of~\cite{\ReedSimonZZI}), nonnegative, self--adjoint
operators on $L^2(\MvN,\psi)$.
{}From~(\psiphi) we see that $\Delta_\psi=D_\phi L_hR_{h^{-1}}$.
Since $\MvN_\phi$ is a factor, and using the spectral decomposition of $\MvN$
with respect to $\phi$, (see Remark~\SpectralTh), we see that $L^2(\MvN,\psi)$
has a basis of eigenvectors of $\Delta_\psi$ if and only if $h$ has purely
atomic spectral measure.
\QED

\proclaim{Proposition \InWhatSenseExtr}
Let $\phi$ and $\psi$ be faithful, almost periodic, normal states on $\MvN$
such that $\psi$ commutes with $\phi$ in the sense
of~\cite{\PedersenTakesaki}.
If $\phi$ is extremal almost
periodic
then the point spectrum of $\Delta_\psi$
contains the point spectrum of $\Delta_\phi$.
\endproclaim
\demo{Proof}
Let $h$, affiliated with $\MvN_\phi$, be such that $\psi=\phi(\cdot h)$.
By Lemma~\WhenAlmPer{} and its proof, $h$ has purely atomic spectral measure
and $\Delta_\psi=D_\phi L_hR_{h^{-1}}$.
Let $\Gamma_h$ be the point spectrum of $h$, so $h=\sum_{t\in\Gamma_h}tE_t$,
where $(E_t)_{t\in\Gamma_h}$ is an orthogonal family of projections in
$\MvN_\phi$ summing to~$1$.
Using the spectral decomposition (see Remark~\SpectralTh) of $\MvN$ with
respect to $\phi$ and the factoriality of $\MvN_\phi$, we see that
$E_{t_1}F_\gamma E_{t_2^{-1}}\neq0$ for every $\gamma\in\Gamma_\phi$ and
every $t_1,t_2\in\Gamma_h$, (where $F_\gamma$ is as in the proof
of~\WhenAlmPer).
Thus the point spectrum of $\Delta_\psi$ is equal to
$\{\gamma t_1t_2^{-1}\mid\gamma\in\Gamma_\phi,\;t_1,t_2\in\Gamma_h\}
\supseteq\Gamma_\phi$.
\QED

\head \FreeProdLemmas.  Some technical lemmas. \endhead

The free etymology proof
of Theorem~1.2 of~\cite{\DykemaZZFreeDim}, which was stated for the finite,
separable case, is valid also more generally.
We restate it for convenience in two steps, first the more often used
version, then its generalization.
\proclaim{Proposition \FreeProdDirectSum}
Suppose $(A_\iota,\phi_\iota)$ is a von Neumann algebra with normal, faithful
state ($\iota=1,2$) and let $(\MvN,\phi)=(A_1,\phi_1)*(A_2,\phi_2)$.
Suppose $h\in A_1$ is a  central projection.
Let $B_1=\Cpx h+(1-h)A_1$ and $\NvN=W^*(B_1\cup A_2)$.
Then $h\MvN h=W^*(h\NvN h\cup A_1h)$ and $h\NvN h$ and $A_1h$ are free in
$(h\MvN h,\phi(h)^{-1}\phi\restrict_{h\MvN h})$.
Moreover, the central support, $\Cc_\MvN(h)$, of $h$ in $\MvN$ equals the
central support, $\Cc_\NvN(h)$, of 
$h$ in $\NvN$.
\endproclaim
\proclaim{Proposition \FreeProdDirectSumMatrixAlg}
Suppose $(A_\iota,\phi_\iota)$ is a von Neumann algebra with normal, faithful
state ($\iota=1,2$) and let $(\MvN,\phi)=(A_1,\phi_1)*(A_2,\phi_2)$.
Suppose $A_1=D_1\otimes M_n(\Cpx)\oplus D_2$ and let
$$ \alignat 2
A_1&=\;&D_1&\otimes M_n(\Cpx)\;\oplus\;D_2 \\
\cup\;& \\
B_1&=&1&\otimes M_n(\Cpx)\;\oplus\;D_2
\endalignat $$
and $\NvN=W^*(B_1\cup A_2)$.
Let $h\in1\otimes M_n(\Cpx)\oplus0\subseteq B_1$ be a minimal projection.
Then $h\MvN h=W^*(h\NvN h\cup(D_1\otimes h))$ and $h\NvN h$ and $D_1\otimes h$
are free in $(h\MvN h,\phi(h)^{-1}\phi\restrict_{h\MvN h})$.
\endproclaim

\proclaim{Lemma \FreeProdDirectSumBoth}
Suppose $(A_\iota,\phi_\iota)$ is a von Neumann algebra with normal, faithful
state ($\iota=1,2$) and let $(\MvN,\phi)=(A_1,\phi_1)*(A_2,\phi_2)$.
Suppose $h_\iota\in A_\iota$ is a central projection
and $\phi_1(h_1)=\phi_2(h_2)$.
Let $B_\iota=\Cpx h_\iota+(1-h_\iota)A_\iota$ and $\NvN=W^*(B_1\cup B_2)$.
There is a unitary $w\in W^*(1,h_1,h_2)$ such that $h_2w=wh_1$.
Then $h_1\MvN h_1=W^*(h_1\NvN h_1\cup A_1h_1\cup w^*A_2h_2w)$
and $(h_1\NvN h_1,A_1h_1,w^*A_2h_2w)$ is free in
$(h_1\MvN h_1,\phi(h_1)^{-1}\phi\restrict_{h_1\MvN h_1})$.
Moreover, the central supports satisfy $\Cc_\MvN(h_1)=\Cc_\NvN(h_1)$.
\endproclaim
\demo{Proof}
That the unitary $w$ exists is easily seen from Theorem~1.1
of~\cite{\DykemaZZFreeDim}.
Let $\NvN_1=W^*(B_1\cup A_2)$.
Then by Proposition~\FreeProdDirectSum,
$h_2\NvN_1h_2=W^*(h_2\NvN h_2\cup h_2A_2)$, $h_2\NvN h_2$ and $h_2A_2$ are free
in $h_2\NvN_1 h_2$ and $\Cc_{\NvN_1}(h_2)=\Cc_\NvN(h_2)$.
Thus $h_1\NvN_1h_1=w^*h_2\NvN_1h_2w=W^*(h_1\NvN h_1\cup w^*h_2A_2w)$,
$h_1\NvN h_1$ and $w^*h_2A_2w$ are free in $h_1\NvN_1h_1$ and
$\Cc_{\NvN_1}(h_1)=\Cc_\NvN(h_1)$.
Applying again Proposition~\FreeProdDirectSum, we have
$h_1\MvN h_1=W^*(h_1\NvN_1h_1\cup h_1A_1)$,
$h_1\NvN_1h_1$ and $h_1A_1$ are free in $h_1\MvN h_1$ and
$\Cc_\MvN(h_1)=\Cc_{\NvN_1}(h_1)$.
Now apply Proposition~2.5.5 of~\cite{\VDNbook}.
\QED

The previous two results were about what happens when one adds central summands
in a free product situation ({\it i.e\.} going from $B_\iota$ to $A_\iota$)
and making the center of $A_\iota$ larger (or keeping it the same). 
The next result is about what happens when adding a partial isometry and making
the center of $A_\iota$ smaller.

\proclaim{Lemma \FreeProdMatUnits}
Suppose $(A_\iota,\phi_\iota)$ is a von Neumann algebra with normal, faithful
state ($\iota=1,2$) and let $(\MvN,\phi)=(A_1,\phi_1)*(A_2,\phi_2)$.
Suppose $v\in A_1$ is a partial isometry
such that $vv^*=k$, $v^*v=h\le1-k$, with $h,k\in(A_1)_{\phi_1}$ and
$$ \forall a\in A_1\quad\phi_1(va)=0\Leftrightarrow\phi_1(av)=0. $$
Set $B_1$ to be either
\roster
\item"(i)" $B_1=\Cpx k+(1-k)A_1(1-k)\subseteq A_1$,
\item"(ii)" or, only if $k$ is minimal in $A_1$,
$B_1=\kt A_1\kt+\htld A_1\htld$
where $\kt,\htld\in(A_1)_{\phi_1}$ are projections, $h\le\htld$, $k\le\kt$ and
$\htld+\kt=1$.
\endroster
Important properties of $B_1$ will be that $k$ is minimal in $B_1$ and
that $kB_1h=\{0\}$.
Let $\NvN=W^*(B_1\cup A_2)$,
let $\NvN_\phi$ denote the centralizer of $\phi\restrict_\NvN$ and suppose
there is $x\in\NvN$ such that $x^*x=h$, $xx^*=k$ and
$$ \forall a\in\NvN\quad\phi(xa)=0\Leftrightarrow\phi(ax)=0. $$
Then $x^*v$ is a Haar unitary and $(h\NvN h,\{x^*v\})$ is $*$--free in
$(h\MvN h,\phi(h)^{-1}\phi\restrict_{h\MvN h})$, and
$$ h\MvN h=W^*(h\NvN h\cup\{x^*v\}). \tag{\hpgen} $$
Moreover, the central supports satisfy $\Cc_\MvN(h)=\Cc_\NvN(h+k)$.
\endproclaim
\demo{Proof}
(\hpgen) is clear by the usual algebraic technique of bringing back
generators with partial isometries,
{\it cf\.} Lemma~1 of~\cite{\RadulescuZZFundGp}.
We will show simultaneously that $x^*v$ is a Haar unitary and that
$(h\NvN h,\{x^*v\})$ is $*$--free by showing that $\phi(y)=0$ for every
alternating product $y$ in $(h\NvN h)\oup$ and
$\{(x^*v)^n\mid n\in\Naturals\}\cup\{(v^*x)^n\mid n\in\Naturals\}$.
Regrouping by sticking $x$ and $x^*$ to letters from $(h\NvN h)\oup$ whenever
possible, we see that $y$ is equal to an alternating product in $\{v,v^*\}$ and
$$ (h\NvN h)\oup\cup(h\NvN h)x^*\cup x(h\NvN h)
\cup x(h\NvN h)\oup x^*, $$
where
\baselineskip=1pt
$$ \spreadlines{0.2ex} \align
\text{every letter from $(h\NvN h)\oup$ }&\text{is (unless it is the first
letter of $y$) preceeded by $v$} \\
\text{and }&\text{is (unless it is the last letter of $y$)
followed by $v^*$;} \\ \vspace{2ex}
\text{every letter from $(h\NvN h)x^*$ }&\text{is (unless it is the first
letter of $y$) preceeded by $v$} \\
\text{and }&\text{is followed by $v$;} \\ \vspace{2ex}
\text{every letter from $x(h\NvN h)$ }&\text{is preceeded by $v^*$} \\
\text{and }&\text{is (unless it is the last letter of $y$)
followed by $v^*$;} \\ \vspace{2ex}
\text{every letter from $x(h\NvN h)\oup x^*$ }&\text{is
preceeded by $v^*$} \\ 
\text{and }&\text{is followed by $v$.}
\endalign $$
\baselineskip=18pt

Using Kaplansky's density theorem, we see that every element of $\NvN$ is
the s.o.--limit of a bounded sequence in
$\lspan(\Omega\cup\{1\})$ where $\Omega=\Lambdao(B_1\oup,A_2\oup)$.
Now using the $\phi$--preserving conditional expectation from $\NvN$ onto
$B_1$ (which is gotten from Proposition~\FreeFactorCondExp), and using the fact
that $k$ is minimal in $B_1$ and $k$ and $h$ are perpendicular projections
belonging to the centralizer of $\phi$, we have that every element of
$$ \align
(h\NvN h)\oup&\text{ is the s.o.--limit of a bounded sequence in }
\lspan\;\Omega_0 \\
x(h\NvN h)&\text{ is the s.o.--limit of a bounded sequence in }
\lspan\;\Omega_r \\
(h\NvN h)x^*&\text{ is the s.o.--limit of a bounded sequence in }
\lspan\;\Omega_l \\
x(h\NvN h)\oup x^*&\text{ is the s.o.--limit of a bounded sequence in }
\lspan\;\Omega_{l,r}
\endalign $$
 where
$$ \align
\Omega_0&=(\Omega\backslash B_1\oup)\cup(hB_1h)\oup \\
\Omega_r&\text{ is the set of elements of $\Omega$ whose right--most letter is
from }A_2\oup \\
\Omega_l&\text{ is the set of elements of $\Omega$ whose left--most letter is
from }A_2\oup \\
\Omega_{l,r}&=\Omega_l\cap\Omega_r.
\endalign $$
Hence it will suffice to show that $\phi(z)=0$ whenever $z$ is an alternating
product in $\{v,v^*\}$ and $\Omega_0\cup\Omega_l\cup\Omega_r\cup\Omega_{l,r}$,
such that every letter that is $v$ is preceded by a letter from
$\Omega_r\cup\Omega_{l,r}$ and every letter that is $v^*$ is followed by a
letter from $\Omega_l\cup\Omega_{l,r}$.
Regrouping by sticking each $v$ and $v^*$ to a letter from $B_1\oup$ whenever
possible, and using that $v(hB_1h)\oup v^*\subseteq A_1\oup$ and
$vB_1\subseteq A_1\oup$, we see that $z$ equals an alternating product in
$A_1\oup$ and $A_2\oup$, hence by freeness, $\phi(z)=0$.

We have left to prove only that $\Cc_\MvN(h)=\Cc_\NvN(h+k)$.
Using $\Cc_\Ac(p)=\bigwedge_{x\in\Ac}\text{ left support}(xp)$, and noting
$vh=v$ has left support equal to $k$,
we have that $\Cc_\MvN(h)=\Cc_\MvN(h+k)\ge\Cc_\NvN(h+k)$.
But since $\Cc_\NvN(h+k)$ commutes with $B_1$, with $A_2$ and with $v$, we
have that $\Cc_\NvN(h+k)$ is in the center of $\MvN$, hence
$\Cc_\MvN(h+k)\le\Cc_\NvN(h+k)$.
\QED

Now we prove that (in certain circumstances) if $A_1$ is freely complemented
in $\MvN$ then $hA_1h$ is freely subcomplemented in $h\MvN h$.

\proclaim{Proposition \FreeProdCutDown}
Let $A_\iota$ be a von Neumann algebra and $\phi_\iota$ a normal,
faithful state on $A_\iota$, ($\iota=1,2$).
Let $(\MvN,\phi)=(A_1,\phi_1)*(A_2,\phi_2)$.
Suppose $(A_1)_{\phi_1}$ is a factor and $(A_2)_{\phi_2}$ has no minimal
projections.
If $1\neq h\in(A_1)_{\phi_1}$ is a self--adjoint projection then
$hA_1h$ is freely subcomplemented in
$(h\MvN h,\phi(h)^{-1}\phi\restrict_{h\MvN h})$ by
$(B_\iota)_{0\le\iota\le n}$, some $n\ge1$,
where $B_0$ is isomorphic to the cut--down of the free product of $A_2$ and a
finite dimensional abelian algebra and where $B_\iota\cong L(\Integers)$
$\forall \iota\ge1$.
\endproclaim
\demo{Proof}
There are $n\ge1$ and partial isometries $v_0,v_1,\ldots,v_n\in(A_1)_{\phi_1}$
such that $v_jv_j^*\le h$ $\forall j$ and $\sum_{j=0}^nv_j^*v_j=1$.
We may assume $v_0=h$.
Let
$$ \align
D_{-1}&=\Cpx h+\Cpx v_1^*v_1+\cdots+\Cpx v_n^*v_n \subseteq(A_1)_{\phi_1} \\
D_0&=D_{-1}+hA_1h\subseteq A_1 \\
D_j&=W^*(D_{j-1}\cup\{v_j\})\subseteq A_1\text{ for }1\le j\le n.
\endalign $$
For $-1\le j\le n$ let $\PvN_j=W^*(D_j\cup A_2)$ so
$(\PvN_j,\phi\restrict_{\PvN_j})\cong(D_j,\phi_1\restrict_{D_j})*(A_2,\phi_2)$.
Clearly $D_n=A_1$ and $\PvN_n=\MvN$.
By Proposition~\FreeProdDirectSum, $hA_1h$ is freely complemented in
$h\PvN_0h$ by $h\PvN_{-1}h$.
{}From~\cite{\DykemaZZFreeDim} we see that there are partial isometries
$x_j$ in the centralizer of $\phi$ restricted to $W^*(D_{-1}\cup(A_2)_{\phi})$
such that $x_j^*x_j=v_j^*v_j$ and $x_jx_j^*=v_jv_j^*$, so
by Lemma~\FreeProdMatUnits{} $h\PvN_{j-1}h$ is freely subcomplemented in
$h\PvN_jh$ by an algebra isomorphic to $L(\Integers)$ with associated
projection equal to $v_jv_j^*\in(A_1)_{\phi_1}$ $\forall 1\le j\le n$.
Thus using Corollary~\FreelySubcomplChain{} and the fact
that $h(A_1)_{\phi_1}h$ is a factor the proposition is proved.
\QED

Now we prove that (in certain circumstances) if $A$ is a freely subcomplemented
von Neumann subalgebra of $\MvN$ then $hAh$ is freely subcomplemented in
$h\MvN h$.

\proclaim{Proposition \FreelySubcomplCutDown}
Let $\MvN$ be a von Neumann algebra with normal, faithful state $\phi$ and
$1\in A\in\MvN$ a von Neumann subalgebra.
Let $A_\phi$ denote the centralizer of $\phi\restrict_A$
and suppose $A_\phi$ is a II$_1$--factor.
Suppose also that $A$ is freely subcomplemented in $\MvN$ by
$(B_\iota)_{\iota\in I}$ 
where for each $\iota\in I$ the centralizer of $\phi\restrict_{B_\iota}$ has
no minimal projections.
If $h\in A_\phi$ is a self--adjoint projection then $hAh$ is freely
subcomplemented in $h\MvN h$ by algebras isomorphic to the $B_\iota$, to
cut--downs of free products of the $B_\iota$ with finite dimensional abelian
algebras and to $L(\Integers)$.
\endproclaim
\demo{Proof}
Denote the associated projections of $(B_\iota)_{\iota\in I}$ by
$(h_\iota)_{\iota\in I}$.
Let us first do the case when $|I|=1$, say $I=\{1\}$.
Since $A_\phi$ as a factor we may assume either $h_1\le h$ or $h\le h_1$.
If $h_1\le h$ then clearly $hAh$ is freely subcomplemented in $h\MvN h$ by
$B_1$.
If $h\le h_1$ then $h_1Ah_1$ and $B_1$ freely generate $h_1\MvN h_1$
so applying Proposition~\FreeProdCutDown{} we see that $hAh$ is freely
complemented in $h\MvN h$ by a cut--down of ($B_1$ free product a finite
dimensional abelian algebra) and copies of $L(\Integers)$, so the case $|I|=1$
is proved.
Moreover in the above case, since $A_\phi$ is a II$_1$--factor we may choose
the associated projections for the free 
subcomplimnetation of $hAh$ in $h\MvN h$ to be in  any given II$_1$--subfactor
of $hAh$.
This allows us to conclude the general case from Corollary~\FreelySubcomplWell.
\QED

\head \FPwrtMAP.  Free products of factors with respect to
extremal almost periodic states. \endhead

\proclaim{Lemma \AddingOneComplimented}
Suppose $\phi_\iota$ is an almost periodic state on a von Neumann
algebra $A_\iota$,
($\iota=1,2$).
Let $\Gamma_{\iota,0}$ be the point spectrum of $\phi_\iota$ and let
$\Gamma_\iota$ be the subgroup of $\Real_+^*$ generated by $\Gamma_{\iota,0}$.
Let
$(\MvN,\phi)=(A_1,\phi_1)*(A_2,\phi_2)$.
Then by Proposition~\AlmPerFreeProd, $\phi$ is almost periodic and its point
spectrum, which we denote 
$\Gamma$, is equal to the subgroup of
$\Real_+^*$ generated by $\Gamma_1\cup\Gamma_2$.
Suppose that $\lambda\in\Gamma_{1,0}$, $0<\lambda<1$, that there is an isometry
$v\in(A_1)_{\phi_1}(\lambda^{-1})$ and that
$\lambda^\Integers\eqdef\{\lambda^n\mid n\in\Integers\}$ is complemented in
$\Gamma_1$, say\footnote{
We will mostly use additive notation for direct sums, complements and cosets of
subgroups of $\Reals_+^*$, even though the operation in $\Reals_+^*$ is
multiplication.}
$\Gamma_1=H\oplus\lambda^\Integers$.
Let $B_1=(A_1)_{\phi_1}(H)$.
Note by Remark~\SpectralTh{} that $\MvN$ is generated by $B_1$ and $A_2$
together with $v$. 
Let
$$ \aligned
\NvN_{[0,0]}&=W^*(B_1\cup A_2), \\
\NvN_{[-n,0]}&=W^*(v^*\NvN_{[-n+1,0]}v\cup\NvN_{[0,0]}),\;n\ge1, \\
\NvN_{(-\infty,0]}&=W^*(\tsize\bigcup_{n\ge0}\NvN_{[-n,0]}), \\
\NvN_{(-\infty,n]}&=W^*(\NvN_{(-\infty,n-1]}\cup v^n\NvN_{[0,0]}(v^*)^n),
 \;n\ge1, \\
\NvN_{(-\infty,\infty)}&=W^*(\tsize\bigcup_{n\ge0}\NvN_{(-\infty,n]}).
\endaligned \tag{\allNvNDef} $$
Then
\roster
\item"(i)" $v^*\NvN_{[-n+1,0]}v$ and $\NvN_{[0,0]}$ are free in $(\MvN,\phi)$
with amalgamation
over $B_1$, $\forall n\ge1$;
\item"(ii)" $h_n\NvN_{(-\infty,n-1]}h_n$ and $v^n\NvN_{[0,0]}(v^*)^n$ are free
in $(h_n\MvN h_n,\lambda^{-n}\phi\restrict_{h_n\MvN h_n})$
with amalgamation over $h_nB_1h_n$, $\forall n\ge1$, where
$h_n=v^n(v^*)^n\in\MvN_\phi$;
\item"(iii)" if $\lambda^\Integers\cap\Gamma'=\{1\}$, where
$\Gamma'$ is the subgroup of $\Gamma$ generated by $H$
together with $\Gamma_2$,
then $\NvN_{(-\infty,\infty)}=\MvN_\phi(\Gamma')$.
\endroster
\endproclaim
\demo{Proof}
We have by Remark~\SpectralTh{} that
$\lspan\big(B_1\cup\bigcup_{n\in\Naturals}(B_1v^n\cup(v^*)^nB_1\big))$ is a
s.o.--dense $*$--subalgebra of $A_1$, hence, letting
$$ \Omega=
\Lambdao\big(B_1\oup\cup\bigcup_{n\in\Naturals}(B_1v^n\cup(v^*)^nB_1),
A_2\oup\big), $$
we have that $\lspan(\Omega\cup\{1\})$ is a s.o.--dense $*$--subalgebra of
$\MvN$. 
Consider a word $x=a_1a_2\cdots a_n$ in $\Omega$, written as an alternating
product of the specified form, {\it i.e\.} each
$a_j\in B_1\oup\cup\bigcup_{n\in\Naturals}(B_1v^n\cup(v^*)^nB_1)\subseteq
A_1\oup$ or $a_j\in A_2\oup$,
with $a_j\in A_1\oup\Longleftrightarrow a_{j+1}\in A_2\oup$.
We will want to keep track of the occurrences of $v$ and $v^*$, hence
let $l(x)=n$ be the length of $x$ and for each $1\le i,j\le n$ let
$$ \align
s(a_i)&=\cases
0&\text{if }a_i\in B_1\oup\cup A_2\oup \\
k&\text{if }a_i\in B_1v^k \\
-k&\text{if }a_i\in(v^*)^kB_1, \endcases \\
t_j(x)&=\sum_{i=1}^js(a_i).
\endalign $$
If $0\in I\subseteq\Integers$ let $\Theta_I$ be the set of alternating products
$x\in\Omega$
such that $t_{l(x)}(x)=0$ and $t_j(x)\in I$ $\forall1\le j\le l(x)$.
\proclaim{Claim \AddingOneComplimented a}
If $I$ is of the form $[-n,0]$, $(-\infty,n]$, $n\in\{0\}\cup\Naturals$, or
$(-\infty,\infty)$ then
$\lspan(\Theta_I\cup\{1\})$ is a s.o.--dense $*$--subalgebra of $\NvN_I$.
\endproclaim
\demo{Proof}
First, show that 
$\lspan(\Theta_I\cup\{1\})$ is a $*$--algebra.  We need only show that
multiplication is closed.  We show by induction on $n+m$ that if we have 
elements of
$\Omega$, $x=a_1\cdots a_n$ (or $x=1$, in
which case we take $n=0$) and $y=b_1\cdots b_m$ (or $y=1$, in which case
$m=0$),
such that $t_j(x)\in I$, $t_{l(x)}(x)+t_j(y)\in I$ $\forall j$ and
$t_n(x)+t_m(y)=0$, then 
$xy\in\lspan(\Theta_I\cup\{1\})$.
When $n+m=0\text{ or }1$, then either $x=1$ or $y=1$, so
$xy\in\Theta_I\cup\{1\}$.
If $a_n\in A_{\iota_1}\oup$, $b_1\in A_{\iota_2}\oup$, $\iota_1\neq\iota_2$,
then $xy=a_1\cdots a_nb_1\cdots b_m$ is already in reduced form, so $xy$
belongs to $\Theta_I$.
If $a_n\in A_\iota\oup$ and $b_1\in A_\iota\oup$, let $c=a_nb_1$.
If $\iota=1$ and $s(a_n)\neq-s(b_1)$ then
$a_nb_1\in(v^*)^kB_1\text{ or }b_1v^k$, some $k>0$, so
$xy=a_1\cdots a_{n-1}cb_2\cdots b_m$ is in reduced form and 
$$ t_j(xy)=\cases t_j(x)&\text{if }1\le j\le n-1 \\
t_{l(x)}(x)+t_{j-n+1}(y)&\text{if }n\le j\le n+m-1, \endcases $$
hence $xy$ is an element of $\Theta_I$.
If $\iota=1$ and $s(a_n)=-s(b_1)$ then $c\in B_1$, or if $\iota=2$ then $c\in
A_2$ and $c=c\oup+\phi(c)1$ where $\phi(c\oup)=0$, so
$$ xy=a_1\cdots a_{n-1}c\oup b_2\cdots b_m
+\phi(c)a_1\cdots a_{n-1}b_2\cdots b_m, $$
and, by the same argument as above, the first term is an element of
$\Theta_I$ while the second term is of the correct form to apply the inductive
hypothesis.
Hence we have shown that $\lspan(\Theta_I\cup\{1\})$ is a $*$--algebra.

Now we show by induction on $n\ge0$ that
$$ \NvN_{[-n,0]}=\clspan(\Theta_{[-n,0]}\cup\{1\}). \tag{\NvNminus} $$
First, note that $\NvN_{[-n,0]}\supseteq\NvN_{[-n+1,0]}$ $\forall n\ge1$,
(which is readily shown from~(\allNvNDef) by induction on $n$).
For $n=0$, (\NvNminus) is clearly true.
Suppose $n\ge1$.
By inductive hypothesis, we have
$v^*\NvN_{[-n+1,0]}v\subseteq\clspan(v^*\Theta_{[-n+1,0]}v\cup\{1\})
\subseteq\clspan(\Theta_{[-n,0]}\cup\{1\})$,
and $\NvN_{[0,0]}\subseteq\clspan(\Theta_{[0,0]}\cup\{1\})$, so the inclusion
$\subseteq$ in (\NvNminus) is true.
For the reverse inclusion, suppose $x=a_1\cdots a_m\in\Theta_{[-n,0]}$ is a
word 
of the specified form.
We show by induction on $m$ that $x\in\NvN_{[-n,0]}$.
If $t_j=0$ $\forall1\le j\le m$ then $x\in\NvN_{[0,0]}$.
Otherwise, let $j_1$ be least such that $t_{j_1}(x)<0$ and $j_2$ least such
that $j_2>j_1$ and $t_{j_2}(x)=0$.
If $j_1>1$ then $a_1\cdots a_{j_1-1}\in\NvN_{[0,0]}$.
Now $a_{j_1}=(v^*)^{k_1}b_1$ and $a_{j_2}=b_2v^{k_2}$, some $b_1,b_2\in B_1$,
$0<k_1,k_2\le n$.
Thus
$(v^*)^{k_1-1}b_1a_{j_1+1}\cdots a_{j_2-1}b_2v^{k_2-1}\in\Theta_{[-n+1,0]}$,
so (by inductive hypothesis for the induction on $n$),
$a_{j_1}\cdots a_{j_2}\in v^*\NvN_{[-n+1,0]}v$.
Finally, $a_{j_2+1}\cdots a_m\in\Theta_{[-n,0]}$, and is of shorter length,
hence 
lies in $\NvN_{[-n,0]}$, (by inductive hypothesis for the induction on $m$).
Thus $a_1\cdots a_m\in\NvN_{[-n,0]}$, and (\NvNminus) is true.

Let us now prove by induction on $n\ge0$ that
$$ \NvN_{(-\infty,n]}=\clspan(\Theta_{(-\infty,n]}\cup\{1\}). \tag{\NvNplus} $$
For $n=0$ we have it by taking the inductive limit of (\NvNminus).
The inclusion $\subseteq$ in (\NvNplus) is true because (using inductive
hypothesis)
$\NvN_{(-\infty,n-1]}\subseteq\clspan(\Theta_{(-\infty,n-1]}\cup\{1\})$ and
$v^n\Theta_{[0,0]}(v^*)^n\subseteq\Theta_{(-\infty,n]}$, while every element of
$\NvN_{[0,0]}$ is the s.o--limit of a bounded sequence in
$\lspan(\Theta_{[0,0]}\cup\{1\})$.
To show $\supseteq$, consider $x=a_1\cdots a_m\in\Theta_{(-\infty,n]}$, a word
of 
the specified form, and show by induction on $m$ that
$x\in\NvN_{(-\infty,n]}$.
If $t_j(x)<n$ $\forall j$, then by inductive hypothesis (for the induction on
$n$),
$x\in\NvN_{(-\infty,n-1]}\subseteq\NvN_{(-\infty,n]}$.
Otherwise let $j_1$ be least such that $t_{j_1}(x)=n$ and $j_2$ least such that
$j_2>j_1$ and $t_{j_2}(x)<n$.
Thus $a_{j_1}=b_1v^{k_1}$ and $a_{j_2}=(v^*)^{k_2}b_2$ for some
$b_1,b_2\in B_1$ and $k_1,k_2>0$.
Then
$x=a_1\cdots a_{j_1}(v^*)^nv^na_{j_1+1}\cdots a_{j_2-1}(v^*)^nv^na_{j_2}\cdots
a_m$.
But $a_1\cdots a_{j_1}(v^*)^n\in\lspan\Theta_{(-\infty,n-1]}
\subseteq\NvN_{(-\infty,n-1]}$,
$v^na_{j_1+1}\cdots a_{j_2-1}(v^*)^n\in v^n(\Theta_{[0,0]}\cup\{1\})(v^*)^n
\subseteq v^n\NvN_{[0,0]}(v^*)^n$
and $v^na_{j_2}\cdots a_m\in\lspan\Theta_{(-\infty,n]}$ is a sum of one or two
words of shorter length, so by
inductive hypothesis, is in $\NvN_{(-\infty,n]}$.
This proves~(\NvNplus).

Finally,
$$ \NvN_{(-\infty,\infty)}=\clspan(\Theta_{(-\infty,\infty)}\cup\{1\}) $$
is obtained by taking the inductive limit of~(\NvNplus).
Hence Claim~{\AddingOneComplimented a} is proved.
\enddemo

We continue with the proof of Lemma~\AddingOneComplimented.
Showing~(i) is the same as showing that $\phi(z)=0$ whenever $z$ is an
alternating product in $v^*\NvN_{[-n+1,0]}v\ominus B_1$ and
$\NvN_{[0,0]}\ominus B_1$.
If $x\in\NvN_{[0,0]}$ then (by Claim~{\AddingOneComplimented a},
Kaplansky's density theorem and separability) $x$ is the strong--operator limit
of a bounded sequence,
$(x_n)_{n=1}^\infty$ in $\lspan(\Theta_{[0,0]}\cup\{1\})$.
Combining Corollary~\FreeFactorCondExp{} and Lemma~\CondExpAlmPer, there is a
$\phi$--preserving conditional expectation, $E_{B_1}$ from $\MvN$ onto $B_1$.
Hence if also $x\perp B_1$, then $E_{B_1}(x)=0$ and thus $x_n-E_{B_1}(x_n)$ is
a bounded sequence converging in s.o.--topology to $x$.
So $x$ is the s.o.--limit of a bounded sequence in
$\lspan(\Theta_{[0,0]}\backslash B_1\oup)$.
Doing the same for the other, note that
$v^*\NvN_{[-n+1,0]}v\ominus B_1=v^*(h_1\NvN_{[-n+1,0]}h_1\ominus h_1B_1h_1)v$.
One similarly shows that if $x\in\NvN_{[-n+1,0]}\ominus B_1$ then $x$ is the
s.o.--limit of a bounded sequence in
$\lspan(\Theta_{[-n+1,0]}\backslash B_1\oup)$.
Hence, to show that $\phi(z)=0$ it will suffice to show that $\phi(y)=0$
whenever $y$ is an alternating product in $\Theta_{[0,0]}\backslash B_1\oup$ and
$v^*(\Theta_{[-n+1,0]}\backslash B_1\oup)v$.
If we ``erase parenthesis'' in $y$ and multiply elements from $A_1$ together
whenever possible, because $B_1vB_1\subseteq A_1\oup$
we see that $y$ is equal to an alternating product in $A_1\oup$ and $A_2\oup$,
hence $\phi(y)=0$ by freeness, and~(i) is proved.

Part~(ii) is proved similarly.
It suffices to show that $\phi(y)=0$ whenever $y$ is an alternating product in
$\Theta_{(-\infty,n-1]}\backslash B_1\oup$ and
$v^n(\Theta_{[0,0]}\backslash B_1\oup)(v^*)^n$.
Since every element of $\Theta_{(-\infty,n-1]}\backslash B_1\oup$ begins with
an element of either $A_2\oup$, $(v^*)^jB_1$, ($j\ge1$), $B_1\oup$ or $B_1v^k$,
($1\le k\le n-1$), and ends with the adjoint of one of these, by regrouping and
multiplying neighbors we see that $y$ is 
equal to an alternating product in $A_1\oup$ and $A_2\oup$.
Hence $\phi(y)=0$ and~(ii) is proved.

Let us now show~(iii).
We assume $\lambda^\Integers\cap\Gamma'=\{1\}$.
In light of Claim~{\AddingOneComplimented a}, we have 
$\NvN_{(-\infty,\infty)}=\clspan(\Theta_\Integers\cup\{1\})$.
Note that $\Theta_\Integers$ equals the set of words $x$ belonging to $\Omega$
such that $t_{l(x)}(x)=0$.
Let $E$ denote the projection from $\MvN$ onto $\MvN_\phi(\Gamma')$, obtained
from Lemma~\CondExpAlmPer.
If $x\in\Omega$, then, writing each element of $A_2\oup$ as the
s.o.--limit of a bounded sequence in
$\lspan\bigcup_{\gamma\in\Sd(\phi_2)}(A_2)_{\phi_2}(\gamma)$ and each element
of $B_1$ as the s.o.--limit of a bounded sequence in
$\lspan\bigcup_{\gamma\in H}(A_1)_{\phi_1}(\gamma)$, we see that
$x\in\clspan\bigcup_{\gamma\in\lambda^{k(x)}+\Gamma'}
\MvN_\phi(\gamma)$, where $k(x)=-t_{l(x)}(x)$.
So $t_{l(x)}(x)=0$
implies $x\in\MvN_\phi(\Gamma')$, and
$t_{l(x)}(x)\neq0$ implies $x\perp\MvN_\phi(\Gamma')$ and hence $E(x)=0$.
Thus we easily see $\MvN_\phi(\Gamma')\supseteq\clspan\Theta_\Integers$.
For the reverse inclusion, let $y\in\MvN_\phi(\Gamma')$.
Since $\lspan(\Omega\cup\{1\})$ is s.o.--dense in $\MvN$, by the Kaplansky
density theorem and separability there is a bounded sequence $(y_n)_{n\ge1}$
converging to $y$ in s.o.--topology and such that
$y_n\in\lspan(\Omega\cup\{1\})$ $\forall n$.
Then, since $E(y)=y$, we have that $(E(y_n))_{n\ge1}$ is a bounded sequence
converging strongly to $y$, and $E(y_n)\in\lspan(\Theta_\Integers\cup\{1\})$
$\forall n$.
Thus~(iii) is proved.

\QED

\proclaim{Lemma \AddingOneTorsion}
Let $(\MvN,\phi)=(A_1,\phi_1)*(A_2,\phi_2)$, $\Gamma_\iota$, $\Gamma_{\iota,0}$
and $\Gamma$ be as in Lemma~\AddingOneComplimented.
Suppose $\lambda\in\Gamma_{1,0}$, $0<\lambda<1$, there is an isometry
$v\in(A_1)_{\phi_1}(\lambda^{-1})$ and there is a subgroup
$H\subseteq\Gamma_1$ that together with $\lambda$ generates
$\Gamma_1$.
Suppose
$\lambda^\Integers\cap H
=\lambda^{N\Integers}\eqdef\{\lambda^{Nk}\mid k\in\Integers\}$, some
$N\in\Naturals$, $N\ge2$.
Let $B_1=(A_1)_{\phi_1}(H)$.
Note by Remark~\SpectralTh{} that $\MvN$ is generated by $B_1$ and $A_2$
together with $v$. 
Let
$$ \aligned
\NvN_{[0,0]}&=W^*(B_1\cup A_2), \\
\NvN_{[0,n]}&=W^*(\NvN_{[0,n-1]}\cup v^n\NvN_{[0,0]}(v^*)^n),\;1\le n\le N-1.
\endaligned $$
Then
\roster
\item"(i)" $h_n\NvN_{[0,n-1]}h_n$ and $v^n\NvN_{[0,0]}(v^*)^n$ are free
in $(h_n\MvN h_n,\lambda^{-n}\phi\restrict_{h_n\MvN h_n})$
with amalgamation over $h_nB_1h_n$, $\forall1\le n\le N-1$,
where $h_n=v^n(v^*)^n$;
\item"(ii)" if also $\lambda^\Integers\cap\Gamma'=\lambda^{N\Integers}$, where
$\Gamma'$ is the subgroup of $\Gamma$ generated by $H$ together with
$\Gamma_2$, then $\NvN_{[0,N-1]}=\MvN(\Gamma')$.
\endroster
\endproclaim
\demo{Proof}
This is a modulo $N$ version of the proof of Lemma~\AddingOneComplimented.
Because $v^N\in B_1$ and, {\it e.g\.} for $1\le n\le N-1$,
$$ (v^*)^nB_1=(v^*)^n(v^*)^{N-n}v^{N-n}B_1(v^*)^{N-n}v^{N-n}
=(v^*)^N\big(v^{N-n}B_1(v^*)^{N-n}\big)v^{N-n}\subseteq B_1v^{N-n},
\tag{\vmodN} $$
we have that $\lspan(B_1\cup\bigcup_{1\le n\le N-1}B_1v^n)$ is a s.o.--dense
$*$--subalgebra of $A_1$, hence, letting
$$ \Omega=\Lambdao(B_1\oup\cup\bigcup_{1\le n\le N-1}B_1v^n,A_2\oup), $$
we have
that $\lspan(\Omega\cup\{1\})$ is a s.o.--dense $*$--subalgebra of $\MvN$.
Consider an alternating product $x=a_1a_2\cdots a_n\in\Omega$
with each
$a_j\in B_1\oup\cup\bigcup_{1\le n\le N-1}B_1v^n\subseteq
A_1\oup$ or $a_j\in A_2\oup$ and
with $a_j\in A_1\oup\Longleftrightarrow a_{j+1}\in A_2\oup$.
We will want to keep track of the number of occurrences of $v$ modulo $N$,
hence 
let $l(x)=n$ be the length of $x$ and for each $1\le i,j\le n$ let
$$ \align
s(a_i)&=\cases
0&\text{if }a_i\in B_1\oup\cup A_2\oup \\
k&\text{if }a_i\in B_1v^k \\ \endcases \\
t_j(x)&=\sum_{i=1}^js(a_i)\in\Integers/N\Integers, \tag{\modNsum}
\endalign $$
taking the sum in~(\modNsum) in the group of integers modulo~$N$.
If $0\in I\subseteq\Integers/N\Integers$ let $\Theta_I$ be the set of words
$x\in\Omega$ 
such that $t_{l(x)}(x)=0$ and $t_j(x)\in I$ $\forall1\le j\le l(x)$.
\proclaim{Claim \AddingOneTorsion a}
For each $1\le n\le N-1$, 
$\lspan(\Theta_{[0,n]}\cup\{1\})$ is a s.o.--dense $*$--subalgebra of
$\NvN_{[0,n]}$. 
\endproclaim
\demo{Proof}
One shows that $\Theta_{[0,n]}$ is a $*$--algebra using induction just like in
the proof of
Claim~{\AddingOneComplimented a},
and that $\NvN_{[0,n]}=\clspan(\Theta_{[0,n]}\cup\{1\})$ using double induction
just like in the proof of~(\NvNplus).
\enddemo

The proof of~(i) is similar to the proof of~\AddingOneComplimented(ii).
It will suffice to show that $\phi(y)=0$ whenever $y$ is an alternating product
in $\Theta_{[0,n-1]}\backslash B_1\oup$ and
$v^n(\Theta_{[0,0]}\backslash B_1\oup)(v^*)^n$.
But every element of  $\Theta_{[0,n-1]}\backslash B_1\oup$ begins with an
element of either $A_2\oup$, $B_1\oup$ or $B_1v^k$, ($1\le k\le n-1$), and ends
with an element of either $A_2\oup$, $B_1\oup$ or $B_1v^{N-k}$,
($1\le k\le n-1$).
Since, for $1\le k\le n-1$, we have
$B_1v^{N-k}v^nB_1\subseteq B_1v^{n-k}\subseteq A_1\oup$ and (see~(\vmodN))
$B_1(v^*)^nB_1v^k\subseteq B_1v^{N-n+k}\subseteq A_1\oup$, by regrouping and
multiplying some neighbors we see that $y$ is equal to an alternating product
in $A_1\oup$ and $A_2\oup$, hence $\phi(y)=0$ by freeness.

In light of Claim~{\AddingOneTorsion a}, showing~(ii) is the same as
showing $\MvN_\phi(\Gamma')=\clspan(\Theta_{[0,N-1]}\cup\{1\})$.
Again, note that $\Theta_{[0,N-1]}$ is the set of words, $x\in\Omega$ such
that $t_{l(x)}(x)=0$.
Let $E$ denote the projection from $\MvN$ onto $\MvN_\phi(\Gamma')$, obtained
from Lemma~\CondExpAlmPer.
As in the proof of~\AddingOneComplimented(iii), we see that every
$x\in\Omega$ is the s.o.--limit of a bounded sequence in
$\lspan(\bigcup_{\gamma\in\lambda^{k(x)}+\Gamma'}\MvN_\phi(\gamma))$, where
$k(x)=-t_{l(x)}(x)$.
Hence $t_{l(x)}(x)=0$ implies $x\in\MvN_\phi(\Gamma')$ and
$t_{l(x)}(x)\in\{1,\ldots,N\!-\!1\}$ 
implies $x\perp\MvN_\phi(\Gamma')$, hence $E(x)=0$.
Thus $\MvN_\phi(\Gamma')\supseteq\clspan(\Theta_{[0,N-1]}\cup\{1\})$.
To prove the reverse inclusion, let $y\in\MvN_\phi(\Gamma')$ and let
$(y_n)_{n=1}^\infty$ be 
a bounded sequence in $\lspan(\Omega\cup\{1\})$ whose s.o.--limit is $y$.  Then
$(E(y_n))_{n=1}^\infty$ is a bounded sequence in
$\lspan(\Theta_{[0,N-1]}\cup\{1\})$ whose s.o.--limit is $E(y)=y$.
Hence $y\in\clspan(\Theta_{[0,N-1]}\cup\{1\})$, proving~(ii).

\QED

\proclaim{Lemma \AddingTwoComplimented}
Let $(\MvN,\phi)=(A_1,\phi_1)*(A_2,\phi_2)$, $\Gamma_\iota$, $\Gamma_{\iota,0}$
and $\Gamma$ be as in Lemma~\AddingOneComplimented.
Suppose $\lambda\in\Gamma_{1,0}\cap\Gamma_{2,0}$,
$0<\lambda<1$, for each $\iota=1,2$ there is an isometry 
$v_\iota\in(A_\iota)_{\phi_\iota}(\lambda^{-1})$ and
$\lambda^\Integers$ is complemented in
$\Gamma_\iota$, say $\Gamma_\iota=H_\iota\oplus\lambda^\Integers$.
Let $B_\iota=(A_\iota)_{\phi_\iota}(H_\iota)$.
Note by Remark~\SpectralTh{} that $\MvN$ is generated by $B_1$ and $B_2$
together with $v_1$ and $v_2$.
For each
$n\ge1$, let $h_{\iota,n}=v_\iota^n(v^*_\iota)^n$.
Then there is a partial isometry $y_n\in W^*(\{1,h_{1,n},h_{2,n}\})$ such that
$y_n^*y_n=h_{1,n}$ and $y_ny_n^*=h_{2,n}$.
Let $u_n=y_n^*v_2^n(v^*_1)^n$.
Let
$$ \aligned
\NvN_{[0,0]}&=W^*(B_1\cup B_2), \\
\NvN_{[0,1]}&=W^*(\NvN_{[0,0]}\cup\{u_1\}), \\
\NvN_{[-n,1]}&=W^*(v_1^*\NvN_{[-n+1,1]}v_1\cup\NvN_{[0,1]}),\;n\ge1, \\
\NvN_{(-\infty,1]}&=W^*(\tsize\bigcup_{n\ge0}\NvN_{[-n,1]}), \\
\NvN_{(-\infty,n]}&=W^*(\NvN_{(-\infty,n-1]}\cup\{u_n\}),\;n\ge2, \\
\NvN_{(-\infty,\infty)}&=W^*(\tsize\bigcup_{n\ge1}\NvN_{(\infty,n]}).
\endaligned $$
Then
\roster
\item"(i)"  $u_1$
is a Haar unitary and $(h_{1,1}\NvN_{[0,0]}h_{1,1},\{u_1\})$ is $*$--free in
$(h_{1,1}\MvN h_{1,1},\lambda^{-1}\phi\restrict_{h_{1,1}\MvN h_{1,1}})$;
\item"(ii)" $v_1^*\NvN_{[-n+1,1]}v_1$ and $\NvN_{[0,1]}$ are free
with amalgamation
over $\NvN_{[0,0]}$ in $(\MvN,\phi)$, $\forall n\ge1$;
\item"(iii)" $u_n$ is a Haar unitary and
$(h_{1,n}\NvN_{(-\infty,n-1]}h_{1,n},\{u_n\})$ is $*$--free in
$(h_{1,n}\MvN h_{1,n},\lambda^{-n}\phi\restrict_{h_{1,n}\MvN h_{1,n}})$,
$\forall n\ge2$;
\item"(iv)" if $\lambda^\Integers\cap\Gamma'=\{1\}$, where
$\Gamma'$ is the subgroup of $\Gamma$ generated by $H_1$
together with $H_2$,
then $\NvN_{(-\infty,\infty)}=\MvN_\phi(\Gamma')$.
\endroster
\endproclaim
\demo{Proof}
Since $h_{1,n}$ and $h_{2,n}$ are free projections having the same weight under
$\phi$,
the existence of the partial isometry
$y_n$ is easily seen from Theorem~1.1 of~\cite{\DykemaZZFreeDim}.
Since
$\lspan(B_\iota\cup\bigcup_{n\ge1}(B_\iota v_\iota^n\cup(v^*_\iota)^nB_\iota))$
is a s.o.--dense $*$--subalgebra of $A_\iota$, letting
$$ \Omega=\Lambdao\bigg(B_1\oup\cup\bigcup_{n\ge1}(B_1v_1^n\cup(v^*_1)^nB_1),\;
B_2\oup\cup\bigcup_{n\ge1}(B_2v_2^n\cup(v^*_2)^nB_2)\bigg), $$
we have that
$\lspan(\Omega\cup\{1\})$ is a s.o.--dense $*$--subalgebra of $\MvN$.
Consider a word $x=a_1a_2\cdots a_n$ belonging to $\Omega$, such that
$a_j\in B_{\iota_j}\oup\cup
\bigcup_{n\in\Naturals}(B_{\iota_j} v_{\iota_j}^n
\cup(v_{\iota_j}^*)^nB_{\iota_j})\subseteq
A_{\iota_j}\oup$,
with $\iota_j\neq\iota_{j+1}$.
We will want to keep track of the occurrences of $v_\iota$ and $v_\iota^*$,
hence 
let $l(x)=n$ be the length of $x$ and for each $1\le i,j\le n$ let
$$ \align
s(a_i)&=\cases
0&\text{if }a_i\in B_1\oup\cup B_2\oup \\
k&\text{if }a_i\in B_1v_1^k\cup B_2v_2^k \\
-k&\text{if }a_i\in(v_1^*)^kB_1\cup(v_2^*)^kB_2, \endcases \\
t_j(x)&=\sum_{i=1}^js(a_i).
\endalign $$
If $0\in I\subseteq\Integers$ let $\Theta_I$ be the set of words $x\in\Omega$
such that $t_{l(x)}(x)=0$ and $t_j(x)\in I$ $\forall1\le j\le l(x)$.
\proclaim{Claim \AddingTwoComplimented a}
If $I$ is of the form $[0,0]$, $[-n,1]$ for $n\in\{0\}\cup\Naturals$,
$(-\infty,n]$ for $n\in\Naturals$ or $(-\infty,\infty)$, then
$\lspan(\Theta_I\cup\{1\})$ is a s.o.--dense $*$--subalgebra of $\NvN_I$.
\endproclaim
\demo{Proof}
One first proves that $\lspan(\Theta_I\cup\{1\})$ is a $*$--algebra like in
Claim~{\AddingOneComplimented a}.
It is clear that $\NvN_{[0,0]}=\clspan(\Theta_{[0,0]}\cup\{1\})$.
Let us prove that
$$ \NvN_{[0,1]}=\clspan(\Theta_{[0,1]}\cup\{1\}). \tag{\NvNtZeroOne} $$
Since $u_1\in\clspan\Theta_{[0,1]}$, (which is seen by writing $y_1$ as the
s.o.--limit of a bounded sequence in $\lspan(\Theta_{[0,0]}\cup\{1\})$), the
inclusion $\subseteq$ in~(\NvNtZeroOne) is clear.
For the reverse inclusion, since $v_\iota^*v_\iota=1$, it will suffice to show
that
$v_1av_1^*,\,v_1av_2^*,\,v_2av_2^*\in\NvN_{[0,1]}$ whenever $a\in B_1\cup B_2$.
If $a\in B_1$ then
$$ \aligned
v_1av_1^*&\in B_1\subseteq\NvN_{[0,1]}, \\
v_1av_2^*&=(v_1av_1^*)(v_1v_2^*y_1)y_1^*
 \in B_1u_1^*y_1\subseteq\NvN_{[0,1]}, \\
v_2av_2^*&=y_1(y_1^*v_2v_1^*)(v_1av_1^*)(v_1v_2^*y_1)y_1^*
 \in y_1u_1B_1u_1^*y_1^*\subseteq\NvN_{[0,1]},
\endaligned \tag{\uOnea} $$
and if $a\in B_2$ then
$$ \aligned
v_1av_1^*&=(v_1v_2^*y)y^*(v_2av_2^*)y(y^*v_2v_1^*)
 \in u_1^*y_1^*B_2y_1u_1\subseteq\NvN_{[0,1]} \\
v_1av_2^*&=(v_1v_2^*y)y^*(v_2av_2^*)\in u_1^*y_1^*B_2\subseteq\NvN_{[0,1]} \\
v_2av_2^*&\in B_2\subseteq\NvN_{[0,1]}.
\endaligned \tag{\uOneb} $$
This proves~(\NvNtZeroOne).

Let us prove by induction on $n\ge0$ that
$$ \NvN_{[-n,1]}=\clspan(\Theta_{[-n,1]}\cup\{1\}). \tag{\NvNtminus} $$
The case $n=0$ was just proved.
Suppose $n\ge1$.
{}From the fact that $v_1^*\Theta_{[-n+1,1]}v_1=\Theta_{[-n,0]}$ and from
inductive hypothesis, one sees that $\subseteq$ holds in~(\NvNtminus).
To show the reverse inclusion, note that every element of $\Theta_{[-n,1]}$
equals a product of elements from $\Theta_{[-n,0]}$ and $\Theta_{[0,1]}$.

Now let us prove by induction on $n\ge1$ that
$$ \NvN_{(-\infty,n]}=\clspan(\Theta_{(-\infty,n]}\cup\{1\}).
\tag{\NvNtplus} $$
The case $n=1$ is true by taking the inductive limit of~(\NvNtminus).
For $n\ge2$, the inclusion $\subseteq$ holds in~(\NvNtplus) by inductive
hypothesis and the 
fact (seen by writing $y_n$ as a s.o.--limit of a bounded sequence in
$\lspan(\Theta_{[0,0]}\cup\{1\})$) that
$u_n\subseteq\clspan\Theta_{[0,n]}$.
For the reverse inclusion, analogous to the proof of~(\NvNplus), we show that
every element of $\Theta_{(-\infty,n]}$ is a product of elements of
$\Theta_{(-\infty,n-1]}$ and elements of the form
$v_{\iota_1}^n\Theta_{[0,0]}(v_{\iota_2}^*)^n$, ($\iota_1,\iota_2\in\{1,2\}$).
By inductive hypothesis, the former are in $\NvN_{(-\infty,n-1]}$, and the
latter are shown to be in $\NvN_{(-\infty,n]}$ using analogues of~(\uOnea)
and~(\uOneb) for $u_n$.

Finally, note that
$$ \NvN_{(-\infty,\infty)}=\clspan(\Theta_{(-\infty,\infty)}\cup\{1\}) $$
follows by taking the inductive limit of~(\NvNtplus).
Hence Claim~{\AddingTwoComplimented a} is proved.
\enddemo

Continuing with the proof of Lemma~\AddingTwoComplimented,
we will now prove~(i) and~(iii) simultaneously.
For consistency of notation, let
$\NvN_{(-\infty,0]}=\clspan(\Theta_{(-\infty,0]}\cup\{1\})$, and note that
$\NvN_{[0,0]}\subseteq\NvN_{(-\infty,0]}$.
(We don't need to show that $\NvN_{(-\infty,0]}$ is an algebra; we only care
that it contains $\NvN_{[0,0]}$.)
We will thus show for every $n\ge1$ that $\phi(z)=0$
whenever $z$ is an alternating product in
$(h_{1,n}\NvN_{(-\infty,n-1]}h_{1,n})\oup$ and
$\{u_n^k\mid k\in\Naturals\}\cup\{(u_n^*)^k\mid k\in\Naturals\}$.
Rewriting $u_n=y_n^*v_2^n(v_1^*)^n$ and $u_n^*=v_1^n(v^*_2)^ny_n$ and
regrouping, we see that $z$ equals and alternating product, call it $z'$, in
$$ (h_{1,n}\NvN_{(-\infty,n-1]}h_{1,n})\oup
\cup(h_{1,n}\NvN_{(-\infty,n-1]}h_{2,n})
\cup(h_{2,n}\NvN_{(-\infty,n-1]}h_{1,n})
\cup(h_{2,n}\NvN_{(-\infty,n-1]}h_{2,n})\oup $$
and $\{v_2^n(v^*_1)^n,v_1^n(v^*_2)^n\}$, where every letter that has
$$ \align
v_2^n(v_1^*)^n\text{ on the left and }
v_1^n(v^*_2)^n\text{ on the right}&\text{ belongs to }
 (h_{1,n}\NvN_{(-\infty,n-1]}h_{1,n})\oup \\
v_1^n(v_2^*)^n\text{ on the left and }
v_2^n(v^*_1)^n\text{ on the right}&\text{ belongs to }
 (h_{2,n}\NvN_{(-\infty,n-1]}h_{2,n})\oup.
\endalign $$
{}From Claim~{\AddingTwoComplimented a} and Kaplansky's density theorem, we know
that every element of $\NvN_{(-\infty,n-1]}$ is the s.o--limit of a bounded
sequence in $\lspan(\Theta_{(-\infty,n-1]}\cup\{1\})$.
Using now the conditional expectation from $\MvN$ onto $B_\iota$, gotten from
combining Corollary~\FreeFactorCondExp{} and Lemma~\CondExpAlmPer, we see that
every element of $(h_{\iota,n}\NvN_{(-\infty,n-1]}h_{\iota,n})\oup$ is the
s.o.--limit of a bounded sequence in
$\lspan((\Theta_{(-\infty,n-1]}\backslash B_\iota\oup)
\cup(h_{\iota,n}B_\iota h_{\iota,n})\oup)$.
Hence, to show that $\phi(z')=0$ it will suffice to show that $\phi(z'')=0$
whenever $z''$ is obtained from $z'$ by replacing every letter of $z'$ that
$$ \align
\text{lies between }v_2^n(v_1^*)^n\text{ on the left and }&
v_1^n(v^*_2)^n\text{ on the right} \\
&\text{ with an arbitrary element of }
(\Theta_{(-\infty,n-1]}\backslash B_1\oup)\cup(h_{1,n}B_1 h_{1,n})\oup) \\
\text{lies between }v_1^n(v_2^*)^n\text{ on the left and }&
v_2^n(v^*_1)^n\text{ on the right} \\
&\text{ with an arbitrary element of }
(\Theta_{(-\infty,n-1]}\backslash B_2\oup)\cup(h_{2,n}B_2 h_{2,n})\oup,
\endalign $$
and replacing every letter of $z'$ not in one of the above two cases with an
arbitrary element of $\Theta_{(-\infty,n-1]}\cup\{1\}$.
Now regrouping and multiplying some neighbors in $z''$ and using that
$(v_\iota^*)^n(h_{\iota,n}B_\iota h_{\iota,n})\oup v_\iota^n
\subseteq B_\iota\oup$ and $v_\iota^nB_\iota\subseteq A_\iota\oup$, we see that
$z''$ is equal to an alternating product in $A_1\oup$ and $A_2\oup$, hence by
freeness $\phi(z'')=0$.
Thus~(i) and~(iii) are proved.

To prove part~(ii), we must show that $\phi(z)=0$ whenever $z$ is an
alternating product in 
$v_1^*\NvN_{[-n+1,1]}v_1\ominus\NvN_{[0,0]}$ and
$\NvN_{[0,1]}\ominus\NvN_{[0,0]}$.
Let $E$ be the conditional expectation from $\MvN$ onto $\NvN_{[0,0]}$ gotten
from Lemma~\FreeProdOfCondExp{} and Lemma~\CondExpAlmPer.
If $x\in v_1^*\NvN_{[-n+1,1]}v_1$, then by Claim~{\AddingTwoComplimented a} and
Kaplansky's density theorem, $x$ is the s.o.--limit of a bounded sequence,
$(x_n)_{n=1}^\infty$, in
$\lspan(v_1^*\Theta_{[-n+1,1]}v_1\cup\{1\})=\lspan(\Theta_{[-n,0]}\cup\{1\})$.
If also $x\perp\NvN_{[0,0]}$, then $E(x)=0$, so $x_n-E(x_n)$ converges strongly
to $x$ as $n\rightarrow\infty$, and
$x_n-E(x_n)\in\lspan(\Theta_{[-n,0]}\backslash\Theta_{[0,0]})$.
Hence every element of $v_1^*\NvN_{[-n+1,1]}v_1\ominus\NvN_{[0,0]}$ is the
s.o.--limit of a bounded 
sequence in $\lspan(\Theta_{[-n,0]}\backslash\Theta_{[0,0]})$.
Similarly, every element of $\NvN_{[0,1]}\ominus\NvN_{[0,0]}$ is the
s.o.--limit of a bounded sequence in
$\lspan(\Theta_{[0,1]}\backslash\Theta_{[0,0]})$.
Hence it will suffice to show that $\phi(z')=0$ whenever $z'$ is an alternating
product in $\Theta_{[-n,0]}\backslash\Theta_{[0,0]}$ and
$\Theta_{[0,1]}\backslash\Theta_{[0,0]}$.
However, note that
$\Theta_{[0,0]}(\Theta_{[-n,0]}\backslash\Theta_{[0,0]})\Theta_{[0,0]}
\subseteq\lspan(\Theta_{[-n,0]}\backslash\Theta_{[0,0]})$.
Hence from every letter of $z'$ that belongs to
$\Theta_{[-n,0]}\backslash\Theta_{[0,0]}$, we can peel off any leading or
trailing letters from $B_1\oup$ or $B_2\oup$ and glue them onto the adjoining
letters of $z'$, (except if we're at the beginning or the end of $z'$).
Therefore, letting $\Psi$ be the set of words
$x=a_1\cdots a_n\in\Theta_{[0,1]}$ 
such that $a_1\in B_\iota v_\iota$, some $\iota=1,2$ and
$a_n\in v_\iota^*B_\iota$, some $\iota=1,2$, it will suffice to show that
$\phi(fz''g)=0$ whenever $z''$ is an alternating product in
$\Theta_{[-n,0]}\backslash\Theta_{[0,0]}$ and $\Psi$, and whenever $f=1$
(respectively $g=1$) if the
first letter (respectively last letter) of $z''$ is from
$\Theta_{[-n,0]}\backslash\Theta_{[0,0]}$
and $f$ (respectively $g$) is in $\Theta_{[0,0]}\cup\{1\}$ if the first letter
(respectively last letter)
of $z''$ is in $\Psi$.
By regrouping and multiplying some neighbors, $fz''g$ is seen to be an
alternating product in $A_1\oup$ and $A_2\oup$, hence $\phi(fz''g)=0$ by
freeness and~(ii) is proved.

In light of the Claim~{\AddingTwoComplimented a}, part~(iv) is proved just
like~\AddingOneComplimented(iii) was proved.

\QED

\proclaim{Lemma \AddingTwoTorsion}
Let $(\MvN,\phi)=(A_1,\phi_1)*(A_2,\phi_2)$, $\Gamma_\iota$, $\Gamma_{\iota,0}$
and $\Gamma$ be as in Lemma~\AddingOneComplimented.
Suppose $\lambda\in\Gamma_{1,0}\cap\Gamma_{2,0}$,
$0<\lambda<1$, for each $\iota=1,2$ there is an isometry 
$v_\iota\in(A_\iota)_{\phi_\iota}(\lambda^{-1})$ and
there is a subgroup $H_\iota\subseteq\Gamma_\iota$
that together with $\lambda$ generates
$\Gamma_\iota$.
Suppose
$\lambda^\Integers\cap H_1=\lambda^{N\Integers}=\lambda^\Integers\cap H_2$,
some $N\in\Naturals$, $N\ge2$.
Let $B_\iota=(A_\iota)_{\phi_\iota}(H_\iota)$.
Note by Remark~\SpectralTh{} that $\MvN$ is generated by $B_1$ and $B_2$
together with $v_1$ and $v_2$.
For each
$n\ge1$, let $h_{\iota,n}=v_\iota^n(v^*_\iota)^n$.
Then there is a partial isometry $y_n\in W^*(\{1,h_{1,n},h_{2,n}\})$ such that
$y_n^*y_n=h_{1,n}$ and $y_ny_n^*=h_{2,n}$.
Let $u_n=y_n^*v_2^n(v^*_1)^n$.
Let
$$ \aligned
\NvN_{[0,0]}&=W^*(B_1\cup B_2), \\
\NvN_{[0,n]}&=W^*(\NvN_{[0,n-1]}\cup\{u_n\}),\;1\le n\le N-1. \\
\endaligned $$
Then
\roster
\item"(i)" $\{u_n\}$ is a Haar unitary and
$(h_{1,n}\NvN_{[0,n-1]}h_{1,n},\{u_n\})$ is $*$--free in
$(h_{1,n}\MvN h_{1,n},\lambda^{-n}\phi\restrict_{h_{1,n}\MvN h_{1,n}})$,
$\forall 1\le n\le N-1$;
\item"(ii)" if $\lambda^\Integers\cap\Gamma'=\lambda^{N\Integers}$, where
$\Gamma'$ is the subgroup of $\Gamma$ generated by $H_1$ and $H_2$, then
$\MvN_\phi(\Gamma')=\NvN_{[0,N-1]}$.
\endroster
\endproclaim
\demo{Proof}
This is a modulo $N$ version of the proof of Lemma~\AddingTwoComplimented.
As in the proof of Lemma~\AddingOneTorsion,
we have that $\lspan(B_\iota\cup\bigcup_{1\le n\le N-1}B_\iota v_\iota^n)$ is a
s.o.--dense $*$--subalgebra of $A_\iota$, hence, letting
$$ \Omega=\Lambdao(B_1\oup\cup\bigcup_{1\le n\le N-1}B_1v_1^n,\;
B_2\oup\cup\bigcup_{1\le n\le N-1}B_2v_2^n), $$
we have that
$\lspan(\Omega\cup\{1\})$ is a s.o.--dense $*$--subalgebra of $\MvN$.
Consider a word $x=a_1a_2\cdots a_n\in\Omega$
with each
$a_j\in B_{\iota_j}\oup\cup
\bigcup_{1\le n\le N-1}B_{\iota_j}v_{\iota_j}^n\subseteq
A_{\iota_j}\oup$
and with $\iota_j\neq\iota_{j+1}$.
We will want to keep track of the number of occurrences of $v_1$ and $v_2$
modulo $N$, hence 
let $l(x)=n$ be the length of $x$ and for each $1\le i,j\le n$ let
$$ \align
s(a_i)&=\cases
0&\text{if }a_i\in B_1\oup\cup B_2\oup \\
k&\text{if }a_i\in B_1v_1^k\cup B_2v_2^k \\ \endcases \\
t_j(x)&=\sum_{i=1}^js(a_i)\in\Integers/N\Integers, \tag{\TmodNsum}
\endalign $$
taking the sum in~(\TmodNsum) in the group of integers modulo~$N$.
If $0\in I\subseteq\Integers/N\Integers$ let $\Theta_I$ be the set of words
$x\in\Omega$ 
such that $t_{l(x)}(x)=0$ and $t_j(x)\in I$ $\forall1\le j\le l(x)$.
\proclaim{Claim \AddingTwoTorsion a}
For each $1\le n\le N-1$, 
$\lspan(\Theta_{[0,n]}\cup\{1\})$ is a s.o.--dense $*$--subalgebra of
$\NvN_{[0,n]}$. 
\endproclaim
\demo{Proof}
One shows that $\Theta_{[0,n]}$ is a $*$--algebra using induction just like in
the proof of
Claim~{\AddingOneComplimented a},
and that $\NvN_{[0,n]}=\clspan(\Theta_{[0,n]}\cup\{1\})$
just like in the proofs of~(\NvNtZeroOne) and~(\NvNtplus).
\enddemo

Proving~(i) is like proving~\AddingTwoComplimented(i) and~(iii).

In light of the Claim~{\AddingTwoTorsion a}, part~(ii) is proved just
like~\AddingOneComplimented(iii) and~\AddingOneTorsion(ii) were proved.
We omit the details.
\QED

\proclaim{Proposition \FPTypeIIIFactors}
Let $A_1$ be a separable, type III factor with extremal almost periodic state
$\phi_1$ and suppose the centralizer of $\phi_1$ is either the hyperfinite
II$_1$--factor, $R$, or $L(\freeF_\infty)$.
Let $A_2$ be a separable von Neumann algebra with normal, faithful state
$\phi_2$ and suppose either
\roster
\item"(diff)" $\phi_2$ is a trace, $A_2$ is either $L(\Integers)$ or an
interpolated free group factor,
\item"(fin)" $\phi_2$ is a trace and $A_2\neq\Cpx$ is the direct sum of a
finite
dimensional algebra, a hyperfinite algebra and (possibly infinitely many)
interpolated free group factors
\endroster
or
\roster
\item"(extr)" $A_2$ is a type III factor and $\phi_2$ is an extremal
almost periodic state whose centralizer is either $R$ or $L(\freeF_\infty)$.
\endroster
Let
$$ (\MvN,\phi)=(A_1,\phi_1)*(A_2,\phi_2). $$
Then $\MvN_\phi\cong L(\freeF_\infty)$.
Furthermore, the inclusion $Q\hookrightarrow\MvN_\phi$ is a standard embedding,
where $Q=W^*((A_1)_{\phi_1}\cup(A_2)_{\phi_2})\cong L(\freeF_t)$.
Thus $\MvN$ is a type~III factor and $\phi$ is an extremal almost periodic
state whose point spectrum $\Sd(\phi)$ is the subgroup of $\Real_+^*$ generated
by $\Sd(\phi_1)\cup\Sd(\phi_2)$.
\endproclaim
Note that part of this proposition is Theorem~\ExtrFreeProd.
\demo{Proof of Proposition~\FPTypeIIIFactors}
We need only show that $\MvN_\phi\cong L(\freeF_\infty)$ and
$Q\hookrightarrow\MvN_\phi$ is a standard embedding.
These assertions will be proved by induction in several steps,  corresponding
to several cases.
First, we show that
$$ \text{(diff)}\implies\text{(fin)}, $$
or more precisely that
whenever the proposition has been proved in case~(diff) for a given value
of $\Sd(\phi_1)$, then the conclusions of the proposition hold also in the
case~(fin)
for that same value of $\Sd(\phi_1)$.
Indeed, in case~(fin), by~\cite{\DykemaZZFreeDim} there is $n$ large enough so
that if $h_1,\ldots,h_n\in(A_1)_{\phi_1}$ are orthogonal projections,
$\phi_1(h_j)=1/n$ $\forall j$, then
$$ \PvN_0\eqdef W^*(\{h_1,\ldots,h_n\}\cup A_2) $$
is an interpolated free group factor.
For $2\le j\le n$ let $v_j\in(A_1)_{\phi_1}$ be a partial isometry such that
$v_j^*v_j=h_1$ and $v_jv_j^*=h_j$.
Then $A_1=W^*(h_1A_1h_1\cup\{v_2,\ldots,v_n\})$.
Let
$$ \align
\PvN_1&=W^*(\PvN_0\cup h_1A_1h_1) \\
\PvN_j&=W^*(\PvN_{j-1}\cup\{v_j\})\quad2\le j\le n.
\endalign $$
Thus $\PvN_n=\MvN$.
We have by Proposition~\FreeProdDirectSum{} that $h_1\PvN_1h_1$ is freely
generated by $h_1\PvN_0h_1$ and $h_1A_1h_1$, so by the proposition in
case~(diff),
and since $h_1$ has central support equal to~$1$, $\PvN_1$ is a type~III
factor, $\phi\restrict_{\PvN_1}$ is an extremal almost periodic state with
centralizer isomorphic to $L(\freeF_\infty)$ and
$Sd(\phi\restrict_{\PvN_1})=\Sd(\phi_1)$.
In exactly the same way, using Lemma~\FreeProdMatUnits, we have
by induction on $2\le j\le n$ that $\PvN_j$ is a type~III
factor, $\phi\restrict_{\PvN_j}$ is an extremal almost periodic state with
centralizer isomorphic to $L(\freeF_\infty)$ and
$Sd(\phi\restrict_{\PvN_j})=\Sd(\phi_1)$.
It only remains to show that $\QvN\hookrightarrow\MvN_\phi$ is a standard
embedding.
We have that
$W^*(h_1\PvN_0h_1\cup h_1(A_1)_{\phi_1}h_1)\hookrightarrow(h_1\PvN_1h_1)_\phi$
is a standard embedding, and by choosing $x$ to lie in $\PvN_0$ for the
applications of Lemma~\FreeProdMatUnits{} to each inclusion
$h_1\PvN_{j-1}h_1\hookrightarrow h_1\PvN_jh_1$, we see that
$h_1(W^*(h_1A_1h_1\cup\{v_2,\ldots,v_j\}\cup A_2))h_1
\hookrightarrow h_1(\PvN_j)_\phi h_1$ is a standard embedding $\forall2\le j\le
n$, hence at $j=n$ we have that $\QvN\hookrightarrow\MvN_\phi$ is a standard
embedding.
Thus the conclusions of the proposition hold also in case~(fin).

Henceforth, whenever $D\subseteq\MvN$ is a von Neumann subalgebra, by $D_\phi$
we will denote the centralizer of $\phi\restrict_D$.

We begin by showing that the conclusions of the proposition hold in the
following cases for $1\le m\le\infty$.
\proclaim{Case I(diff)$_{\bold m}$}
$\phi_2$ is a trace, $A_2$ is $L(\Integers)$ or an interpolated free group
factor and $\Sd(\phi_1)$ is a free abelian group of rank $m$.
\endproclaim
\proclaim{Case I(extr)$_{\bold m}$}
$A_2$ is a type~III factor, $\phi_2$ is an extremal almost periodic state,
$\Sd(\phi_2)=\Sd(\phi_1)$ and $\Sd(\phi_1)$ is a free abelian group of rank $m$.
\endproclaim
Our strategy will be to prove
\roster
\item"(a)" I(diff)$_1$ and
\item"" I(diff)$_{m-1}\implies$I(diff)$_m$, ($2\le m<\infty$)
\item"(b)" I(extr)$_1$ and
\item"" I(diff)$_{m-1}+$I(extr)$_{m-1}\implies$I(extr)$_m$, ($2\le m<\infty$)
\item"(c)" I(diff)$_\infty$ and I(extr)$_\infty$.
\endroster
Let us begin by proving~(a), {\it i.e\.} I(diff)$_m$, $1\le m<\infty$.
We will use Lemma~\AddingOneComplimented{} with
$\Sd(\phi_1)=\lambda^\Integers\oplus H$, for $H$ a free abelian group of rank
$m-1$, ($H=\{1\}$ when $m=1$).
Then also part~\AddingOneComplimented(iii) applies, and hence
$\MvN_\phi=(\NvN_{(-\infty,\infty)})_\phi$.
We will show by induction on $n\ge1$ that
\roster
\item"(a1)" $\NvN_{[-n,0]}$ is a factor, $\phi\restrict_{\NvN_{[-n,0]}}$ is an
extremal almost periodic state with $\Sd(\phi\restrict_{\NvN_{[-n,0]}})=H$.
\item"(a2)" $B_1$ is freely subcomplemented in $\NvN_{[-n,0]}$ by a finite
family $(D_\iota)_{\iota\in I}$ of algebras, where each $D_\iota$ is either
$L(\Integers)$
or $L(\freeF_t)$, $1<t\le\infty$ with $\phi\restrict_{D_\iota}$ being a trace.
\item"(a3)" $\NvN_{[-n+1,0]}$ is freely subcomplemented in $\NvN_{[-n,0]}$ by a
finite family of algebras, each of which is of the form described in~(a2).
\item"(a4)" $(\NvN_{[-n,0]})_\phi\cong L(\freeF_{t_n})$ with $t_n=\infty$ (or
if $m=1$ possibly $t_n<\infty$ but $t_n\nearrow\infty$).
\item"(a5)" the inclusion
$(\NvN_{[-n+1,0]})_\phi\hookrightarrow(\NvN_{[-n,0]})_\phi$ is a standard
embedding.
\item"(a6)" the inclusion
$\QvN\hookrightarrow(\NvN_{[-n,0]})_\phi$ is a standard
embedding.
\endroster
By definition, $B_1$ is freely complemented in $\NvN_{[0,0]}$ by $A_2$,
hence by I(diff)$_{m-1}$ (or by~\cite{\DykemaZZFreeDim} in the case $m=1$), we
have that~(a1), (a2), (a4) and~(a6) hold for $n=0$.
Now suppose $n\ge1$.
We use~\AddingOneComplimented(i) and we want to apply
Lemma~\FreelySubcomplAmalg.
We need that $B_1$ is freely subcomplemented in $v^*\NvN_{[-n+1,0]}v$, and
since $v^*B_1v=B_1$, this is equivalent to $h_1B_1h_1$ being freely
subcomplemented in $h_1\NvN_{[-n+1,0]}h_1$.
{}From~(a2)$_{n-1}$ and from Lemma~\FreelySubcomplCutDown, making use
of~I(fin)$_{m-1}$, we see that $h_1B_1h_1$ is freely subcomplemented in
$h_1\NvN_{[-n+1,0]}h_1$ by algebras of the form described in~(a2).
Hence, by~\FreelySubcomplAmalg{}, we have that~(a2)$_n$ holds.
If $n=1$ then also by~\FreelySubcomplAmalg{}, (a3)$_n$ holds.
If $n\ge2$ then by~(a2)$_{n-2}$ and~(a3)$_{n-1}$,
applying~\FreelySubcomplCutDown{} twice, we have that $B_1$ is freely
subcomplemented in $v^*\NvN_{[-n+2,0]}v$, which is in turn freely
subcomplemented in $v^*\NvN_{[-n+1,0]}v$ by $(D_\iota)_{\iota\in I}$, all by
algebras of the form described 
in~(a2).
Thus, using the idea of the proof of~\FreelySubcomplAmalg, we see that
$\NvN_{[-n+1,0]}=W^*(v^*\NvN_{[-n+2,0]}v\cup\NvN_{[0,0]})$ is freely
subcomplemented in $\NvN_{[-n,0]}$ by $(D_\iota)_{\iota\in I}$,
{\it i.e\.}~(a3)$_n$ holds.
Now~(a4)$_n$, (a5)$_n$ and~(a6)$_n$ all follow from~(a3)$_n$
using~I(diff)$_{m-1}$ and~I(extr)$_{m-1}$.

{}From these we conclude that $\NvN_{(-\infty,0]}$ is a factor,
$\phi\restrict_{\NvN_{(-\infty,0]}}$ is an extremal almost periodic state with
\linebreak
$\Sd(\phi\restrict_{\NvN_{(-\infty,0]}})=H$ and
$(\NvN_{(-\infty,0]})_\phi\cong L(\freeF_\infty)$.
Moreover, we have that
$\QvN\hookrightarrow(\NvN_{(-\infty,0]})_\phi$ is a
standard embedding and that $B_1$ is freely subcomplemented in
$\NvN_{(-\infty,0]}$ by a family of algebras of the form described in~(a2).
We can then use~\AddingOneComplimented(ii) to prove by induction on $n\ge1$
that
\roster
\item"(a7)" $\NvN_{(-\infty,n]}$ is a factor,
$\phi\restrict_{\NvN_{(-\infty,n]}}$ is an 
extremal almost periodic state with $\Sd(\phi\restrict_{\NvN_{(-\infty,n]}})=H$
and $(\NvN_{(-\infty,n]})_\phi\cong L(\freeF_\infty)$.
\item"(a8)" $h_nB_1h_n$ is freely subcomplemented in $h_n\NvN_{(-\infty,n]}h_n$
by a countable family of algebras, each of which is of the form described
in~(a2).
\item"(a9)" $h_n\NvN_{(-\infty,n-1]}h_n$ is freely subcomplemented in
$h_n\NvN_{(-\infty,n]}h_n$ by a finite family of algebras, each of which is of
the form described in~(a2).
\item"(a10)" the inclusion
$(\NvN_{(-\infty,n-1]})_\phi\hookrightarrow(\NvN_{(-\infty,n]})_\phi$ is a
standard embedding.
\item"(a11)" the inclusion
$\QvN\hookrightarrow(\NvN_{(-\infty,n]})_\phi$ is a standard
embedding.
\endroster

Hence by taking inductive limits we have that $\NvN_{(-\infty,\infty)}$ is a
factor,
$\phi\restrict_{\NvN_{(-\infty,\infty)}}$ is extremal almost periodic with
$\Sd(\phi\restrict_{\NvN_{(-\infty,\infty)}})=H$ and
$(\NvN_{(-\infty,\infty)})_\phi\cong L(\freeF_\infty)$ and
$\QvN\hookrightarrow(\NvN_{(-\infty,\infty)})_\phi$ is a standard embedding.
Thus part~(a) is proved.

Now we prove~(b), {\it i.e\.} I(extr)$_m$, $1\le m<\infty$.
We will use Lemma~\AddingTwoComplimented{} with
$\Sd(\phi_1)=\Sd(\phi_2)=\lambda^\Integers\oplus H$, and note
that~\AddingTwoComplimented(iv) applies, so we need only show that
$(\NvN_{(-\infty,\infty)})_\phi\cong L(\freeF_\infty)$ and
$\QvN\hookrightarrow(\NvN_{(-\infty,\infty)})_\phi$ is a standard embedding.
We have by~I(extr)$_{m-1}$ (or by~\cite{\DykemaZZInterp} if $m=1$) that
$\NvN_{[0,0]}$ is a factor, $\phi\restrict_{\NvN_{[0,0]}}$ is an extremal
almost periodic state, $\Sd(\phi\restrict_{\NvN_{[0,0]}})=H$,
$(\NvN_{[0,0]})_\phi\cong L(\freeF_\infty)$ (or possibly $\cong L(\freeF_2)$ if
$m=1$) and $\QvN\hookrightarrow(\NvN_{[0,0]})_\phi$ is a standard embedding (or
$\QvN=(\NvN_{[0,0]})_\phi$ if $m=1$).
Now we will show by induction on $n\ge1$ that
\roster
\item"(b1)" $\NvN_{[-n,1]}$ is a factor, $\phi\restrict_{\NvN_{[-n,1]}}$ is an
extremal almost periodic state with $\Sd(\phi\restrict_{\NvN_{[-n,1]}})=H$.
\item"(b2)" $\NvN_{[0,1]}$ is freely subcomplemented in $\NvN_{[-n,1]}$ by a
finite family $(D_\iota)_{\iota\in I}$ of algebras, where each $D_\iota$ is
either $L(\Integers)$,
$L(\freeF_t)$, $1<t\le\infty$ with $\phi\restrict_{D_\iota}$ being a trace,
or is a type~III factor with $\phi\restrict_{D_\iota}$ being an
extremal
almost periodic state having centralizer isomorphic to $L(\freeF_\infty)$ and
$\Sd(\phi\restrict_{D_\iota})=H$.
\item"(b3)" $\NvN_{[-n+1,1]}$ is freely subcomplemented in $\NvN_{[-n,1]}$ by a
finite family of algebras, each of which is of the form described in~(b2).
\item"(b4)" $(\NvN_{[-n,1]})_\phi\cong L(\freeF_{t_n})$ with $t_n=\infty$ (or
if $m=1$ possibly $t_n<\infty$ but $t_n\nearrow\infty$).
\item"(b5)" the inclusion
$(\NvN_{[-n+1,1]})_\phi\hookrightarrow(\NvN_{[-n,1]})_\phi$ is a standard
embedding.
\item"(b6)" the inclusion
$\QvN\hookrightarrow(\NvN_{[-n,1]})_\phi$ is a standard
embedding.
\endroster
By~\AddingTwoComplimented(i) and the proposition in case~I(diff)$_m$, we see
that~(b1), (b4) and~(b6) hold for $n=0$.
Now suppose $n\ge1$.
We want to use~\AddingTwoComplimented(ii) and Proposition~\FreelySubcomplAmalg,
and hence we must show that $\NvN_{[0,0]}$ is freely subcomplemented in
$v_1^*\NvN_{[-n+1,1]}v_1$.
This is equivalent to showing that $v_1\NvN_{[0,0]}v_1^*$ is freely
subcomplemented in $h_1\NvN_{[-n+1,1]}h_1$.
Applying Lemma~\FreelySubcomplCutDown{} and using~(b2)$_{n-1}$
if $n\ge2$, it thus
suffices that $v_1\NvN_{[0,0]}v_1^*$ be freely subcomplemented in
$h_1\NvN_{[0,1]}h_1$.
Now $v_1\NvN_{[0,0]}v_1^*$ is generated by $v_1B_1v_1^*=h_1B_1h_1$ and
$v_1B_2v_1^*=u_1^*y_1^*B_2y_1u_1$ (with $u_1$ and $y_1$ as
in~\AddingTwoComplimented).
Consider
$$ \alignat 2
D=W^*(\{1\}\cup h_1B_1h_1\cup h_2B_2h_2)=&\;
 (\Cpx(1-h_1)\oplus h_1B_1h_1&)*(\Cpx(1-h_2)\oplus&h_2B_2h_2) \\ \vspace{2ex}
D\supseteq
D_1=&(\Cpx(1-h_1)\oplus\;\Cpx h_1\;&)*(\Cpx(1-h_2)\oplus&\;\Cpx h_2\;).
\endalignat $$
Then by Lemma~\FreeProdDirectSumBoth, $h_1B_1h_1$ and $y_1^*B_2y_1$ are free
and $W^*(h_1B_1h_1\cup y_1^*B_2y_1)$ is freely complemented in $h_1Dh_1$ by
$h_1D_1h_1$, which by Theorem~1.1 of~\cite{\DykemaZZFreeDim}
is of the form $L(\Integers)$ or $L(\Integers)\oplus\Cpx$.
Hence, adding partial isometries from $(B_1)_{\phi_1}$ and $(B_2)_{\phi_2}$
like in the proof of Lemma~\FreelySubcomplCutDown, we see that
$W^*(h_1B_1h_1\cup y_1^*B_2y_1)$ is freely subcomplemented in
$h_1\NvN_{[0,0]}h_1$ by a finite family $(D_\iota)_{\iota\in I}$ of algebras of
the form $L(\Integers)$ or $L(\Integers)\oplus\Cpx$, and all associated
projections can be chosen to lie in $h_1(B_1)_\phi h_1$.
Because $h_1\NvN_{[0,0]}h_1$ and $\{u_1\}$ together freely generate
$h_1\NvN_{[0,1]}h_1$, by a free etymology argument it follows that $h_1B_1h_1$,
$u_1^*y_1^*B_2y_1u_1$,
$(D_\iota)_{\iota\in I}$ and $\{u_1\}$ together generate
$h_1\NvN_{[0,0]}h_1$, that $(h_1B_1h_1,u_1^*y_1^*B_2y_1u_1,\{u_1\})$ is
$*$--free in $h_1\NvN_{[0,1]}h_1$ and that
$W^*(h_1B_1h_1\cup u_1^*y_1^*B_2y_1u_1\cup\{u_1\})$ is freely subcomplemented
in $h_1\NvN_{[0,1]}h_1$ by $(D_\iota)_{\iota\in I}$;
therefore
$v_1\NvN_{[0,0]}v_1^*=W^*(h_1B_1h_1\cup u_1^*y_1^*B_2y_1u_1)$ is freely
subcomplemented in $h_1\NvN_{[0,1]}h_1$ by algebras of the form described
in~(a2).
Now Proposition~\FreelySubcomplAmalg{} can be used to show~(b2)$_n$ and, in
conjunction with~(b3)$_{n-1}$ if $n\ge2$, to show~(b3)$_n$.
Then~(b4)$_n$, (b5)$_n$ and~(b6)$_n$ follow as a matter of course.

Thus we see that $\NvN_{(-\infty,1]}$ is a factor,
$\phi\restrict_{\NvN_{(-\infty,1]}}$ is extremal almost periodic and has
centralizer isomorphic to $L(\freeF_\infty)$.
Now by~\AddingTwoComplimented(iii) it is immediate that $\NvN_{(-\infty,1]}$
is freely subcomplemented in $\NvN_{(-\infty,\infty)}$ by a countable family of
algebras isomorphic to $L(\Integers)$ and from this~(b) is easily proved.

Then~(c) follows by taking inductive limits, and hence the proposition is
proved in 
Case~I(diff) and Case~I(extr).

Next, consider
\proclaim{Case II} $\phi_2$ is a trace and $A_2$ is either $L(\Integers)$
or an interpolated free group factor.
\endproclaim

Of course, if $\Sd(\phi_1)$ is a free abelian group, then this has been
proved as Case~I(diff).
Since every countable subgroup of $\Reals_+^*$ is the union of a chain of
subgroups $H_1\subseteq H_2\subseteq\cdots$ such that each $H_j$ is a free
abelian group and such that $H_{j+1}/H_j$ is a finite cyclic group $\forall j$,
one can prove the proposition in Case~II using Case~I(diff) and
Lemma~\AddingOneTorsion.

Thus all cases~(diff) and ~(fin) have been proved.
Consider now
\proclaim{Case III} $A_1$ is a type~III factor, $\phi_2$ is an extremal almost
periodic state and $\Sd(\phi_1)=\Sd(\phi_2)$.
\endproclaim

This is easily proved using Case~I(extr), the chain
$H_1\subseteq H_2\subseteq\cdots$ and Lemma~\AddingTwoTorsion.

It remains to show the case when $A_1$ and $A_2$ are type~III factors but
$\Sd(\phi_1)\neq\Sd(\phi_2)$.
Let $H=\Sd(\phi_1)\cap\Sd(\phi_2)$ and let
$H=K_1^{(\iota)}\subseteq K_2^{(\iota)}\subseteq\cdots\subseteq\Sd(\phi_\iota)$
($\iota=1,2$) be a finite or infinite sequence of subgroups such that
$\Sd(\phi_\iota)=\bigcup_jK_j^{(\iota)}$ and $K^{(\iota)}_{j+1}/K_j^{(\iota)}$
is cyclic $\forall j$.
If $x$ is an element of $K_{j+1}^{(\iota)}$ that together with $K_j^{(\iota)}$
generates $K_{j+1}^{(\iota)}$, then
$\Integers x\cap K_j^{(\iota)}=\Integers
x\cap(K_j^{(\iota)}+\Sd(\phi_{\iota'}))$ where $\iota'\neq\iota$.
Hence by using Case~III to get started, considering in turn
$W^*((A_1)_{\phi_1}(K^{(1)}_j)\cup A_2(H))$ and
$W^*(A_1\cup(A_2)_{\phi_2}(K^{(2)}_j))$ and using
Lemmas~\AddingOneComplimented{} 
and~\AddingOneTorsion, one proves the proposition.
\QED

\head \FPfd.  Free products of finite dimensional algebras. \endhead

\proclaim{Definition \pAphi}\rm
For a finite dimensional von Neumann algebra $A$ and a faithful positive linear
functional $\phi$ on
$A$, we define the nonnegative integer $\ntr(A,\phi)$, which very roughly
measures how far
$\phi$ is from being a trace, as follows.
If $A=M_n(\Cpx)$ and
$$ \phi=\Tr\left(\cdot
 \pmatrix\lambda_1&&0\\&\ddots\\0&&\lambda_n\endpmatrix\right)
\tag{\TracialWt} $$
then $\ntr(M_n(\Cpx),\phi)=0$ if $\phi$ is a trace and $\ntr(M_n(\Cpx),\phi)=n$
otherwise.
For a general finite dimensional $A$, we may write
$A=\bigoplus_{j=1}^N A_j$ where each $A_j$ is a matrix algebra and we define
$\ntr(A,\phi)=\sum_{j=1}^N \ntr(A_j,\phi\restrict_{A_j})$.
\endproclaim

\proclaim{Remark \PtSpecFd}\rm
Let us take this opportunity to observe that the point spectrum of the modular
operator, $\Delta_\phi$, of a faithful state $\phi$ on a finite dimensional
algebra $A$ is easily computed.
Indeed, if $A=M_n(\Cpx)$ and $\phi$ is as in~(\TracialWt), then
$$ \text{point spectrum}(\Delta_\phi)=\biggl\{\frac{\lambda_i}{\lambda_j}
 \biggm|1\le i,j\le n\biggr\} $$
and for a faithful state on a general finite dimensional $A$, the point
spectrum of $\Delta_\phi$ is just the union of the sets of such ratios for all
the simple summands of $A$.
\endproclaim

\proclaim{Proposition \FPExtrAPfd}
Let
$$ (\MvN,\phi)=(A_1,\phi_1)*(A_2,\phi_2) $$
where $A_1$ is a type~III factor, $\phi_1$ is an extremal almost periodic state
on $A_1$ with centralizer isomorphic to $R$ or $L(\freeF_\infty)$, $A_2$ is
finite dimensional, $A_2\neq\Cpx$ and $\phi_2$ is faithful on $A_2$.
Then $\MvN$ is a type~III factor, $\phi$ is an extremal almost periodic state
with centralizer isomorphic to $L(\freeF_\infty)$ and $\Sd(\phi)$ is the
subgroup of $\Reals_+^*$ generated by
$$ \Sd(\phi_1)\cup\bigl(\text{point spectrum}(\Delta_{\phi_2})\bigr). $$
Furthermore, if $(A_1)_{\phi_1}\cong L(\freeF_\infty)$ then the inclusion
$(A_1)_{\phi_1}\hookrightarrow\MvN_\phi$ is a standard embedding.
\endproclaim
\demo{Proof}
We proceed by induction on $\ntr(A_2,\phi_2)$.
If $\ntr(A_2,\phi_2)=0$ then the conclusions of the proposition hold by
Proposition~\FPTypeIIIFactors(fin).
For $\ntr(A_2,\phi_2)\ge1$ there is a $*$--subalgebra $B_2\subseteq A_2$ and a
partial isometry $v\in A_2$ that together with $B_2$ generates $A_2$, such that
$\ntr(B_2,\phi_2\restrict_{B_2})<\ntr(A_2,\phi_2)$, $v^*v$ and $vv^*$
are orthogonal
minimal projections in $(A_2)_{\phi_2}$ and $(v^*v)B_2(vv^*)=\{0\}$.
Indeed, just let $v\in (A_2)_{\phi_2}(\lambda^{-1})$ be a partial isometry
between appropriate orthogonal minimal projections
of $A_2$ and let $B_2=(1-vv^*)A_2(1-vv^*)+\Cpx v^*v$.
Let $\QvN=W^*(A_1\cup B_2)$.
Then
$$ (\QvN,\phi\restrict_\QvN)\cong(A_1,\phi_1)*(B_2,\phi_2\restrict_{B_2}),
\tag{\Qfp} $$
and by inductive hypothesis, the conclusion of the proposition applies to the
free product in~(\Qfp).
We have $\MvN=W^*(\QvN\cup\{v\})$.
Let $p=v^*v$, $q=vv^*$ and assume without loss of generality that
$0<\lambda\eqdef\frac{\phi(q)}{\phi(p)}<1$.
The proof will be for three separate cases:
\roster
\item"(A)" $\lambda\in\Sd(\phi\restrict_\QvN)$;
\item"(B)" $\lambda^\Integers\cap\Sd(\phi\restrict_\QvN)=\{1\}$;
\item"(C)" $\lambda^\Integers\cap\Sd(\phi\restrict_\QvN)=
\lambda^{N\Integers},$ some $N\ge2$.
\endroster
\demo{Proof in Case A}
There is $y\in\QvN_\phi(\lambda^{-1})$ such that $y^*y=p$ and $yy^*=q$.
Then by Lemma~\FreeProdMatUnits, it follows that
$$ (p\MvN p,\phi(p)^{-1}\phi\restrict_{p\MvN p})
\cong(p\QvN p,\phi(p)^{-1}\phi\restrict_{p\QvN p})
*(L(\Integers),\tau_\Integers), $$
so by Proposition~\FPTypeIIIFactors{} we have that $p\MvN p$ is a type~III
factor, $p\MvN_\phi p\cong L(\freeF_\infty)$ and
$\Sd(\phi\restrict_{p\MvN p})=\Sd(\phi\restrict_{p\QvN p})
=\Sd(\phi\restrict_\QvN)$.
Hence the conclusions of the proposition hold in Case~A.
\enddemo

\demo{Proof in Case B}
There is a partial isometry $y\in\QvN_\phi$ such that  $yy^*=q$,
$h\eqdef y^*y\le p$.
Then $p\MvN p=W^*(p\QvN p\cup\{y^*v\})$ and $w\eqdef y^*v$ is a nonunitary
isometry belonging to $\MvN_\phi(\lambda^{-1})$.
Let
$$ \align
\QvN_{[0,0]}&=p\QvN p \\
\QvN_{[-n,0]}&=W^*(w^*\QvN_{[-n+1,0]}w\cup\QvN_{[0,0]}),\;(n\ge1) \\
\QvN_{(-\infty,0]}&=W^*(\tsize\bigcup_{n\ge1}\QvN_{[-n,0]}) \\
\QvN_{[0,n]}&=W^*(w\QvN_{[0,n-1]}w^*\cup\QvN_{[0,0]}),\;(n\ge1) \\
\QvN_{(0,\infty]}&=W^*(\tsize\bigcup_{n\ge1}\QvN_{[0,n]}) \\
\QvN_{(-\infty,\infty)}&=W^*(w^*\QvN_{(-\infty,0]}w\cup\QvN_{[0,\infty)}).
\endalign $$
Then
\roster
\item"(i)" $w^*\QvN_{[-n+1,0]}w$ and $\QvN_{[0,0]}$ are free in $p\MvN p$
$\forall n\ge1$.
\item"(ii)" $w\QvN_{[0,n-1]}w^*$ and $h\QvN_{[0,0]}h$ are free in
$h\MvN h$ $\forall n\ge1$.
\item"(iii)" $w^*\QvN_{(-\infty,0]}w$ and $\QvN_{[0,\infty)}$ are free in
$p\MvN p$.
\item"(iv)"
$\QvN_{(-\infty,\infty)}=p\MvN_\phi(\Gamma)p$ where
$\Gamma=\Sd(\phi\restrict_\QvN)$.
\endroster
In order to prove~(i)--(iv), we will want to keep track of occurrences of $w$
and $w^*$.
Hence for
$$ x\in\Lambdao((p\QvN p)\oup,\{w^n,(w^*)^n\mid n\ge1\}), \tag{\xAltProd} $$
let $l(x)$ denote the length of the word $x$ and define $t_j(x)$ for
$1\le j\le l(x)$ as follows.
Setting $t_0(x)=0$ we let
$$ t_j(x)=\cases t_{j-1}(x)&\text{if the $j$th letter of $x$ is from
 $(p\QvN p)\oup$} \\
t_{j-1}(x)+k&\text{if the $j$th letter of $x$ is $w^k$} \\
t_{j-1}(x)-k&\text{if the $j$th letter of $x$ is $(w^*)^k$}.
\endcases $$
Let $\chi$ be the set of words $x$ as in~(\xAltProd) such that any letter of
$x$ coming from $(p\QvN p)\oup$ and lying between $w^*$ on the left and $w$ on
the right belongs to $(h\QvN h)\oup$.
\proclaim{Claim  \FPExtrAPfd a}
Let $I\subseteq\Integers$ be an interval of the form $[-n,0]$, $[0,n]$ (for
$n\ge1$), $(-\infty,0]$, $[0,\infty)$ or $(-\infty,\infty)$.
Then
$$ \QvN_I=\lspan(\{p\}\cup\chi_I)\text{ where }
\chi_I=\{x\in\chi\mid t_{l(x)}(x)=0,\,t_j(x)\in I\,\forall1\le j\le l(x)\}.
\tag{\chiI} $$
\endproclaim
This claim is easily proved, ({\it cf}\/ the proof of
Claim~\AddingOneComplimented a, {\it etcetera}).

Let us now show that $\chi_{(-\infty,\infty)}\subseteq\ker\phi$.
Indeed, writing out $w=y^*v$, we see that every element
$x\in\chi_{(-\infty,\infty)}$ is equal to an
$x'\in\Lambdao(\QvN\oup,\{v,v^*\})$ of which each letter coming from $\QvN\oup$
and
$$ \align
\text{lying between }v\text{ on the left and }v\text{ on the right}
&\text{ belongs to }p\QvN q, \\ 
\text{lying between }v\text{ on the left and }v^*\text{ on the right}
&\text{ belongs to }(p\QvN p)\oup, \\
\text{lying between }v^*\text{ on the left and }v\text{ on the right}
&\text{ belongs to }(q\QvN q)\oup, \\
\text{lying between }v^*\text{ on the left and }v^*\text{ on the right}
&\text{ belongs to }q\QvN p.
\endalign $$
Now using the conditional expectation of $\QvN$ onto $B_2$, using that $p$ and
$q$ are minimal projections in $B_2$ and that $pB_2q=\{0\}$, we see that every
element of
$p\QvN q\cup(p\QvN p)\oup\cup(q\QvN q)\oup\cup q\QvN p$
is the
s.o.--limit of a bounded net in $\Lambdao(A_1\oup,B_2\oup)\backslash B_2\oup$.
Thus, regrouping $x'$ by sticking $v$ and $v^*$ to $B_2\oup$ whenever possible,
using that $p,q\in(A_2)_{\phi_2}$, we see that $x'$ equals an element of
$\Lambdao(A_1\oup,A_2\oup)$, hence by freeness $\phi(x')=0$.
Thus we have proved $\chi_{(-\infty,\infty)}\subseteq\ker\phi$.

To show~(i) we need that $\phi(a)=0$ whenever
$a\in\Lambdao((w^*\QvN_{[-n+1,0]}w)\oup,(p\QvN p)\oup)$.
Because $w$ is in a spectral subspace $\MvN_\phi(\lambda^{-1})$ of $\MvN$,
$(w^*\QvN_{[-n+1,0]}w)\oup=w^*(h\QvN_{[-n+1,0]}h)\oup w$.
Using~(\chiI) and the conditional expectation of $\QvN_{[-n+1,0]}$ onto
$\QvN_{[0,0]}$ (gotten from~\FreeFactorCondExp{} and~(i) for~$n-1$), we have
that every element of $(h\QvN_{[-n+1,0]}h)\oup$ is the s.o.--limit of a bounded
net in $\lspan((\chi_{[-n+1,0]}\backslash(p\QvN p)\oup)\cup(h\QvN h)\oup)$.
Thus it will suffice to show $\phi(a')=0$ whenever
$$ a'\in\Lambdao\Bigl(w^*\bigl((\chi_{[-n+1,0]}\backslash(p\QvN p)\oup)
\cup(h\QvN h)\oup\bigr)w,(p\QvN p)\oup\Bigr). $$
It is easily seen that $a'\in\chi_{[-n,0]}\subseteq\chi_{(-\infty,\infty)}$,
which we know is in $\ker\phi$, so~(i) is proved.

To show~(ii) is even easier.
Substitute $\chi_{[0,n-1]}$ for $(\QvN_{[0,n-1]})\oup$ and get an element of
$\chi_{[0,n]}\subseteq\chi_{(-\infty,\infty)}$.

The proof of~(iii) is similar to that of~(i).
We want, however, a $\phi$--preserving conditional expectation from
$\QvN_{(-\infty,0]}$ onto
$\QvN_{[0,0]}$, which follows by taking inductive limits of conditional
expectations from $\QvN_{[-n,0]}$ onto $\QvN_{[0,0]}$, which exists
by~\S\CondExp.

Putting together~(i), (ii) and~(iii), applying Proposition~\FPTypeIIIFactors{}
repeatedly and taking inductive limits, we get that $\QvN_{(-\infty,\infty)}$
is a type~III factor, $\phi\restrict_{\QvN_{(-\infty,\infty)}}$ is extremal
almost periodic with centralizer isomorphic to $L(\freeF_\infty)$ and with
$\Sd\bigl(\phi\restrict_{\QvN_{(-\infty,\infty)}}\bigr)
=\Sd(\phi\restrict_\QvN)=\Gamma$.
Since $w\in\MvN_\phi(\lambda^{-1})$ and
$w^*\QvN_{(-\infty,\infty)}w=\QvN_{(-\infty,\infty)}$, using the hypothesis of
Case~B we see that~(iv) holds and the conclusion of the proposition holds in
Case~B.
\enddemo

\demo{Proof in Case C}
Note that the proof of~(i), (ii) and~(iii) above did not depend on the
assumption $\lambda^{\Integers}\cap\Sd(\phi\restrict_\QvN)=\{1\}$, but remains
valid also in Case~C.
Using the same notation as in the above proof of~(i), (ii) and~(iii), we have
that $\QvN_{[0,N-1]}$ is a type~III factor, $\phi\restrict_{\QvN_{[0,N-1]}}$
is an extremal almost periodic state with centralizer isomorphic to
$L(\freeF_\infty)$ and with
$\Sd\bigl(\phi\restrict_{\QvN_{[0,N-1]}}\bigr)=\Sd(\phi\restrict_\QvN)$.
We have $ww^*=y^*y\in\QvN_{[0,0]}$ and then after $N-1$ steps of~(ii) one has
that $w^N(w^*)^N\in\QvN_{[0,N-1]}$.
Thus there exists a partial isometry $z$ in the spectral subspace
$(\QvN_{[0,N-1]})_\phi(\lambda^{-N})$ such that $z^*z=p$ and
$zz^*=h'\eqdef w^N(w^*)^N$.
Then $z^*w^N\in p\MvN_\phi p$ is a unitary.
We will show that in $p\MvN p$, $z^*w^N$ is a Haar unitary and $\{z^*w^N\}$ and
$\QvN_{[0,N-1]}$ are $*$--free.
Both will be proved at once by showing $\phi(x)=0$ whenever
$$ x\in\Lambdao((\QvN_{[0,N-1]})\oup,\{(z^*w^N)^n\mid n\ge1\}
\cup\{((w^*)^Nz)^n\mid n\ge1\}). $$
Now $x$ is equal to
$$ x'\in\Lambdao(\QvN_{[0,N-1]},\{w^N,(w^*)^N\}) $$
where each letter of $x'$ that
$$ \align
\text{lies between $w^N$ on the left and $(w^*)^N$ on the right }
&\text{belongs to }(\QvN_{[0,N-1]})\oup \\
\text{lies between $(w^*)^N$ on the left and $w^N$ on the right }
&\text{belongs to }(h'\QvN_{[0,N-1]}h')\oup.
\endalign $$
Since $h'\le h$ and using~(\chiI) as well as the $\phi$--preserving conditional
expectation from $\QvN_{[0,N-1]}$ onto $\QvN_{[0,0]}$,
it will suffice to show
$\phi(x'')=0$ whenever $x''$ is obtained from $x'$ by replacing
$$ \align
\text{each letter from }(h'\QvN_{[0,N-1]}h')\oup&\text{ with an arbitrary
element of }\bigl(\chi_{[0,N-1]}\backslash(p\QvN p)\oup\bigr)
\cup(h\QvN h)\oup, \\
\text{each remaining letter from }(\QvN_{[0,N-1]})\oup&\text{ with an
arbitrary element of }\chi_{[0,N-1]}, \\
\text{each remaining letter from }\QvN_{[0,N-1]}&\text{ with an
arbitrary element of }\chi_{[0,N-1]}\cup\{p\}.
\endalign $$
Cancelling $w^*w=p$ whenever such neighbors occur, we see that $x''$ is equal
to an element of $\chi$, hence $\phi(x'')=0$.
Let $\PvN=W^*(\QvN_{[0,N-1]}\cup\{z^*w^N\})$.
Then using Proposition~\FPTypeIIIFactors, we see that $\PvN$ is a type~III
factor, $\phi\restrict_\PvN$ is extremal almost periodic with centralizer
isomorphic to $L(\freeF_\infty)$ and with
$\Sd(\phi\restrict_\PvN)=\Sd(\phi\restrict_\QvN)$.
Moreover, $\MvN=W^*(\PvN\cup\{w\})$ and $w^N\in\PvN$.
{}From this and the assumption of Case~C, one easily proves
$\MvN_\phi=\PvN_\phi$ and that the conclusions of the proposition hold in
Case~C.
\enddemo
\QED

\proclaim{Definition \NotationFD}\rm
We define the notation we shall use for faithful states on finite dimensional
algebras (and other algebras).
Consider a pair $(\MvN,\phi)$ where $\MvN$ is a von Neumann algebra of the form
$$ \MvN=D\quad\text{or}\quad\MvN=\QvN\oplus D, $$
with $D$ finite dimensional and where $\phi$ is a normal, faithful state on
$\MvN$.
To describe $D$ and the functional $\phi\restrict_D$, we will use the following
notation.
For $\MvN=\QvN\oplus D$, we write
$$ (D,\phi\restrict_D)=\bigoplus_{j=1}^K
\smdp{M_{n_j}(\Cpx)\hfil}{\alpha_{j,1},\ldots,\alpha_{j,n_j}}
{p_{j,1},\ldots,p_{j,n_j}} $$
which indicates $D=\bigoplus_{j=1}^K M_{n_j}(\Cpx)$ and
$$ \phi\restrict_D=\oplus_{j=1}^K\phi_j\text{ with }
\phi_j=\Tr_{n_j}\left(\cdot
\left(\smallmatrix\alpha_j,1&&0\\&\ddots\\
0&&\alpha_{j,n_j}\endsmallmatrix\right)\right), $$
where $\Tr_{n_j}$ is the unnormalized trace on $M_{n_j}(\Cpx)$.
By $p_{j,i}$ we denote the projection $\diag(0,\cdots,0,1,0,\cdots,0)$ with $1$
in the $i$th place.
We must have that $\alpha_{j,i}>0$ $\forall j,i$ and, if $\MvN=D$, that
$\sum_{j=1}^K\sum_{i=1}^{n_j}\alpha_{j,i}=1$.
If $\MvN=\QvN\oplus D$ then we may write
$$ (\MvN,\phi)=\smdp\QvN{\alpha_0}{p_0}\oplus\bigoplus_{j=1}^K
\smdp{M_{n_j}(\Cpx)\hfil}{\alpha_{j,1},\ldots,\alpha_{j,n_j}}
{p_{j,1},\ldots,p_{j,n_j}} $$
where $p_0=1_\QvN\oplus0$ and $\phi(p_0)=\alpha_0$.
Since $\phi$ is a state, we have
$\alpha_0+\sum_{j=1}^K\sum_{i=1}^{n_j}\alpha_{j,i}=1$.
\endproclaim

The following proposition contains Theorem~\MainOne{} of the introduction.

\proclaim{Proposition \FPFD}
Let
$$ (\MvN,\phi)=(A_1,\phi_1)*(A_2,\phi_2) $$
where for each $\iota=1,2$, $A_\iota\neq\Cpx$ is a separable von Neumann
algebra of the form 
$D_\iota$, $\QvN_\iota$ or $\QvN_\iota\oplus D_\iota$, where $D_\iota$ is
finite dimensional,
$\QvN_\iota$ is a type~III factor, $\phi_\iota$ is a normal, faithful state on
$A_\iota$
and, whenever $A_\iota=\QvN_\iota$ or $\QvN_\iota\oplus D_\iota$,
$\phi_\iota\restrict_{\QvN_\iota}$ is
extremal almost periodic, with centralizer isomorphic to $L(\freeF_\infty)$ or
the hyperfinite II$_1$ factor $R$.

Then
$$ \MvN=\MvN_0\quad\text{or}\quad\MvN=\MvN_0\oplus D, \tag{\DorNot} $$
where $D$ is finite dimensional.
If $A_1$ and $A_2$ both have linear dimension two, then
$\MvN_0=L(\Integers)\otimes M_2(\Cpx)$.
If $\phi_1$ and $\phi_2$ are traces, (which implies $A_1=D_1$ and $A_2=D_2$),
and at least one of $A_1$ and $A_2$ has dimension $\ge3$, then
$\MvN_0=L(\freeF_t)$, some $1<t<\infty$.
If at least one of $\phi_1$ and $\phi_2$ is not a trace, then $\MvN_0$ is a
type~III factor, $\phi\restrict_{\MvN_0}$ is extremal almost periodic with
centralizer isomorphic to $L(\freeF_\infty)$ and with
$\Sd(\phi\restrict_{\MvN_0})$
equal to the subgroup of $\Reals_+^*$ generated by
$$ \text{point spectrum}(\Delta_{\phi_1})
\;\cup\text{ point spectrum}(\Delta_{\phi_2}). $$
It remains to decide whether $D$ is present in~(\DorNot), and in that case to
describe it.
If $A_1=\QvN_1$ or $A_2=\QvN_2$, then $D$ is not present, {\it i.e\.}
$\MvN=\MvN_0$.
Otherwise, with notation as in Definition~\NotationFD, suppose
$$\aligned
D_1&=\bigoplus_{j=1}^{K_1}\smdp{M_{n_j}(\Cpx)}
{\alpha_{j,1},\cdots,\alpha_{j,n_j}}{p_{j,1},\cdots,p_{j,n_j}} \\
D_2&=\bigoplus_{j=1}^{K_2}\smdp{M_{m_j}(\Cpx).}
{\beta_{j,1},\cdots,\beta_{j,m_j}}{q_{j,1},\cdots,q_{j,m_j}}
\endaligned \tag{\DoneDtwo} $$
If $n_j\ge2$ $\forall1\le j\le K_1$ and $m_j\ge2$ $\forall1\le j\le K_2$, then
$D$ is not present, {\rm i.e\.} $\MvN=\MvN_0$.
If some $n_j=1$ or some $m_j=1$, we may suppose without loss of generality that
$n_1=1$ and that $p_{1,1}$ is the largest (with respect to $\phi$) of the
projections in $D_1$ and $D_2$ that are both minimal and central, {\rm i.e\.}
we suppose that
$$ \alpha_{1,1}=\sup\bigl(\{\alpha_{j,1}\mid1\le j\le K_1,\,n_j=1\}\cup
\{\beta_{j,1}\mid1\le j\le K_2,\,m_j=1\}\bigr). $$
Let
$$ J=\biggl\{1\le j\le K_2\biggm|
\sum_{i=1}^{m_j}\frac1{\beta_{j,i}}<\frac1{1-\alpha_{1,1}}\biggr\}. $$
If $J$ is empty, then $\MvN=\MvN_0$.
If $J$ is nonempty then
$$ \MvN=\MvN_0\oplus\bigoplus_{j\in J}\smdp{M_{m_j}(\Cpx)}
{\gamma_{j,1},\cdots,\gamma_{j,m_j}}{r_{j,1},\cdots,r_{j,m_j}}
\tag{\Mresult} $$
where
$$ \align
\gamma_{j,i}&=\beta_{j,i}
 \biggl(1-(1-\alpha_{1,1})\sum_{t=1}^{m_j}\frac1{\beta_{j,t}}\biggr) \\
r_{j,i}&=\bigwedge_{t=1}^{m_j}
v_{i,t}^{(j)}(p_{1,1}\wedge q_{j,t})v_{t,i}^{(j)},
\endalign $$
where $(v^{(j)}_{s,t})_{1\le s,t\le m_j}$ is a system of matrix units for the
$j$th summand, $\MvN_{m_j}(\Cpx)$, of $D_2$ such that $v_{t,t}^{(j)}=q_{j,t}$.
\endproclaim
\proclaim{Remark \MatrixUnitsMeet}\rm
We have that $(v^{(j)}_{s,1}r_{j,1}v^{(j)}_{1,t})_{1\le s,t\le m_j}$ is a
system of matrix units for summand $j$ of~(\Mresult), which is equal to the
meet (see~\cite{\DykemaZZAlmPer}),
$\{p_{1,1}\}\owedge\{v^{(j)}_{s,t}\mid1\le s,t\le m_j\}$, of matrix
units of the summand of $D_1$ and the summand of $D_2$ that were mated to
produce the $j$th summand of~(\Mresult).
\endproclaim
\demo{Proof of Theorem~\FPFD}
The cases
$$ \aligned
A_1=\QvN_1&\text{ and }A_2=\QvN_2 \\
A_1=\QvN_1&\text{ and }A_2=D_2 \\
A_1=\QvN_1&\text{ and }A_2=\QvN_2\oplus D_2
\endaligned $$ %\tag{\QdominateCases} $$
can be proved by variously applying Propositions~\FPTypeIIIFactors,
\FPExtrAPfd{} and~\FreeProdDirectSum.
In the remaining cases, {\it i.e\.}
$$ A_1=\QvN_1\oplus D_1\text{ or }A_1=D_1\qquad
\text{and}\qquad A_2=\QvN_2\oplus D_2\text{ or }A_2=D_2, \tag{\QplusDcases} $$
the proposition will be proved by induction on
$\ntr(D_1,\phi_1\restrict_{D_1})+\ntr(D_2,\phi_2\restrict_{D_2})$.
If the proposition has been proved in the case of $A_1$ and $A_2$ both finite
dimensional, for a given value of
$\ntr(D_1,\phi_1\restrict_{D_1})+\ntr(D_2,\phi_2\restrict_{D_2})$, then it can
be shown to hold also in the remaining cases of~(\QplusDcases) by using
Propositions~\FreeProdDirectSum, \FPTypeIIIFactors{} and~\FPExtrAPfd{}.

Hence we set about proving the case when $A_1=D_1$ and $A_2=D_2$ by induction.
If
$\ntr(A_1,\phi_1)+\ntr(A_2,\phi_2)=0$ then
$\phi_1$ and $\phi_2$ are traces and the proposition was in this case proved
in~\cite{\DykemaZZFreeDim}.
Otherwise, taking $\ntr(A_2,\phi_2)\ge1$ and arguing as
in the beginning of the proof
of~(\FPExtrAPfd), there are
$0<\lambda<1$, $v\in(A_2)_{\phi_2}(\lambda^{-1})$ such that $h=v^*v$ and
$k=vv^*$ are orthogonal minimal projections of $A_2$, and a $*$--subalgebra
$B_2\subseteq A_2$
which together with $v$ generates $A_2$ and such that $h$ is central in $B_2$
and
$\ntr(B_2,\phi_2\restrict_{B_2})<\ntr(A_2,\phi_2)$.
Let $\NvN=W^*(A_1\cup B_2)$.
By inductive hypothesis,
$$ \NvN=\smp{\NvN_0}{\vphantom{\bigm|}r_0}\quad\text{or}\quad
\NvN=\smp{\NvN_0}{\vphantom{\bigm|}r_0}\oplus\bigoplus_{j\in J}
\smdp{M_{l_j}(\Cpx)}{\gamma'_{j,1},\ldots,\gamma'_{j,l_j}}
{r'_{j,1},\ldots,r'_{j,l_j}} $$
with $\NvN_0$ nonatomic and either $\phi\restrict_{\NvN_0}$ a trace or
$\phi\restrict_{\NvN_0}$ extremal almost periodic with centralizer isomorphic
to $L(\freeF_\infty)$ and  with
$\Sd(\phi\restrict_{\NvN_0})$ equal to
the subgroup of $\Real_+^*$ generated by
point spectrum$(\Delta_{\phi_1})
\cup\text{ point spectrum}(\Delta_{\phi_2\restrict_{B_2}})$.
Clearly $(h+k)\MvN(h+k)=W^*((h+k)\NvN(h+k)\cup\{v\})$.
Thus $\Cc_\MvN(h+k)=\Cc_\NvN(h+k)$ and once we have determined
$(h+k)\MvN(h+k)$, we will be able to find $\MvN$.

One of the following four cases must hold:
\roster
\item"(I)" $k\le r_0$ and $h\le r_0$;
\item"(II)" $k\le r_0$ and $h\not\le r_0$;
\item"(III)" $k\not\le r_0$ and $h\le r_0$;
\item"(IV)" $k\not\le r_0$ and $h\not\le r_0$.
\endroster
By the fact that $h$ is central in $B_2$ and by inductive hypothesis, we see
that case~III cannot occur.
We will show in cases~I and~II that $(h+k)\MvN(h+k)$ is a type~III
factor, meaning that $h$ and $k$ contribute only to the type~III part of
$\MvN$, whereas in case~IV, if $h$ and $k$ are large enough then they can
contribute to a type~I direct summand of $\MvN$.

Case~I is treated in three separate subcases:
\roster
\item"(Ia)" $\NvN_0$ is a type III factor and
  $\lambda\in\Sd(\phi\restrict_{\NvN_0})$;
\item"(Ib)" $\phi\restrict_{\NvN_0}$ is a trace or $\NvN_0$ is a type~III
  factor and $\lambda^\Integers\cap\Sd(\phi\restrict_{\NvN_0})=\{1\}$;
\item"(Ic)"  $\NvN_0$ is a type~III factor and
$\lambda^\Integers\cap\Sd(\phi\restrict_{\NvN_0})=\lambda^{N\Integers}$, some
$N\ge2$.
\endroster
The proposition in each of these cases is proved like the corresponding case in
the proof of Proposition~\FPExtrAPfd.

In Case~II,
$$ (h+k)\NvN(h+k)=\smp{(h+k)\NvN_0(h+k)}{r_0(h+k)}\oplus
\bigoplus_{j\in J'}\smp{M_{n_j}(\Cpx)}{r'_{j,1},\ldots,r'_{j,n_j}}, $$
with $\emptyset\neq J'\subseteq\{1,\ldots,K_1\}$ and $r'_{j,i}\le h$
$\forall(j,i)$.
We will now show that $\phi(r_0h)\ge\phi(k)$.
If $h$ is the largest central and minimal projection in $A_1$ and $B_2$,
then
$$ \align
\phi(r_{j,i}')&=\alpha_{j,i}\bigg(1-(1-\phi(h))
\sum_{t=1}^{n_j}\frac1{\alpha_{j,t}}\bigg) \\
&\le\alpha_{j,i}\bigg(1-(1-\phi(h))\frac1{\alpha_{j,i}}\bigg) \\
&=\phi(h)+\alpha_{j,i}-1.
\endalign $$
Let $N=\sum_{j\in J'}n_j$.
If $N=1$ then $J'=\{j\}$ and
$\phi(r_0h)=\phi(h)-\phi(r_{j,1}')\ge1-\alpha_{j,1}\ge1-\phi(h)\ge\phi(k)$.
If $N\ge2$ then
$$ \align
\phi(r_0h)&=\phi(h)-\sum_{j\in J',\;1\le i\le n_j}\phi(r_{j,i}') \\
 \vspace{2ex}
&\ge\phi(h)-\sum_{j\in J',\;1\le i\le n_j}(\phi(h)+\alpha_{j,i}-1) \\
 \vspace{2ex}
&\ge(N-1)(1-\phi(h))\ge1-\phi(h)\ge\phi(k).
\endalign $$
If $h$ is not the largest central and minimal projection in $A_1$ and $B_2$,
then, by inductive hypothesis and the assumption on the form of $\NvN_0$, the
largest central and minimal projection in $A_1$ and $B_2$ must be
$p=p_{j,1}\in A_1$.
Hence $J'=\{j\}$, $n_j=1$, $r_{j,1}'=p\wedge h$,
$\phi(r_{j,1}')=\phi(p)+\phi(h)-1$ and $\phi(p)+\phi(k)<1$.
Thus $\phi(r_0h)=1-\phi(p)>\phi(k)$.

Consider the three subcases:
\roster
\item"(IIa)" $\lambda\in\Sd(\phi\restrict_\NvN)$;
\item"(IIb)" $\lambda^\Integers\cap\Sd(\phi\restrict_\NvN)=\{1\}$;
\item"(IIc)" $\lambda^\Integers\cap\Sd(\phi\restrict_\NvN)=
\lambda^{N\Integers},$ some $N\ge2$.
\endroster
In Case~IIa, there is $y\in\NvN_\phi(\lambda^{-1})$ such that
$h'\eqdef r_0h=y^*y$ and $k'\eqdef yy^*\le k$, thus
$$ h\MvN h=W^*(v^*\NvN v\cup (h-h')\NvN(h-h')\cup\{y^*v\}). $$
Let
$$ \PvN=W^*(v^*\NvN v\cup (h-h')\NvN(h-h')). $$
One easily shows that $v^*\NvN v$ and $\Cpx h'+(h-h')\NvN(h-h')$ are
free in
$(h\MvN h,\phi(h)^{-1}\phi\restrict_{h\MvN h})$, so by
Proposition~\FPExtrAPfd{}
$\PvN$ is a type~III factor and $\phi\restrict_\PvN$ is extremal almost
periodic with centralizer isomorphic to $L(\freeF_\infty)$ and with
$\Sd(\phi\restrict_\PvN)=\Sd(\phi\restrict_\NvN)$.
Now $v^*yy^*v=v^*k'v$ and $y^*vv^*y=h'$, and there is a partial isometry
$z\in\PvN_\phi$ such that $z^*z=h'$ and $zz^*=v^*k'v$.
Let $w=y^*vz$.
Then $h'\MvN h'=W^*(h'\PvN h'\cup\{w\})$ and we claim that $w$ is a
Haar unitary and that $\{w\}$ and $h'\PvN h'$ are $*$--free in
$h'\MvN h'$ with respect to the restriction of $\phi$.
Clearly $w$ is unitary in $h'\MvN h'$.
It will suffice to show that $\phi(x)=0$ whenever
$$ x\in\Lambdao((h'\PvN h')\oup,\{w^n\mid n\in\Naturals\}
\cup\{(w^*)^n\mid n\in\Naturals\}). $$
Writing $w=y^*vz$ and $w^*=z^*v^*y$ and then sticking $z$ and $z^*$ to
letters from $(h'\PvN h')\oup$ whenever possible, we see that
$$ x=x'\in\Lambdao(\{y^*v,v^*y\},(h'\PvN h')\oup\cup h'\PvN(v^*k'v)
\cup(v^*k'v)\PvN h'\cup((v^*k'v)\PvN(v^*k'v))\oup), $$
where each letter of $x'$ that is $y^*v$ \newline
\centerline{is either the first letter of $x'$
or has to the left a letter from $(h'\PvN h')\oup$ or  $(v^*k'v)\PvN h'$}
\centerline{and}
\centerline{ has to the right a letter from $(v^*k'v)\PvN h'$
 or $((v^*k'v)\PvN(v^*k'v))\oup$,}
and where each letter of $x'$ that is $v^*y$ \newline
\centerline{ has to the left a letter from $h'\PvN(v^*k'v)$ or
 $((v^*k'v)\PvN(v^*k'v))\oup$}
\centerline{and}
\centerline{is either the last letter of $x'$
or has to the right a letter from $(h'\PvN h')\oup$ or $h'\PvN(v^*k'v)$.}
Each element of
$$ \align
\PvN&\text{ is the s.o.--limit of a bounded sequence in }
 \lspan(\{h\}\cup\Omega) \\
(h'\PvN h')\oup&\text{ is the s.o.--limit of a bounded sequence in }
 \lspan(h'\Omega_{l,r}h') \\
h'\PvN(v^*k'v)&\text{ is the s.o.--limit of a bounded sequence in }
 \lspan(\{h'\}\cup h'\Omega_l) \\
(v^*k'v)\PvN h'&\text{ is the s.o.--limit of a bounded sequence in }
 \lspan(\{h'\}\cup\Omega_rh') \\
((v^*k'v)\PvN(v^*k'v))\oup&\text{ is the s.o.--limit of a bounded sequence in }
 \lspan((\Omega\backslash(v^*\NvN v)\oup)\cup v^*(k'\NvN k')\oup v),
\endalign $$
where
$$ \align
\Omega&=\Lambdao((v^*\NvN v)\oup,(\Cpx h'+(h-h')\NvN(h-h'))\oup) \\
\Omega_l&\text{ is the set of words belonging to $\Omega$ whose first letter
 belongs to $(v^*\NvN v)\oup$} \\
\Omega_r&\text{ is the set of words belonging to $\Omega$ whose last letter
 belongs to $(v^*\NvN v)\oup$} \\
\Omega_{l,r}&=\Omega_l\cap\Omega_r.
\endalign $$
Hence, by substituting these into $x'$, we see that in order to show that
$\phi(x')=0$ it will suffice to show that $\phi(x'')=0$ whenever
$$ x''\in\Lambdao(\{y^*v,v^*y\},\{h\}\cup\Omega), $$
and when each letter of $x''$ that
$$ \align
\text{has $y^*v$ on the left and $y^*v$ on the right }&
 \text{belongs to }\{h\}\cup\Omega_r \\
\text{has $y^*v$ on the left and $v^*y$ on the right }&
 \text{belongs to }
 (\Omega\backslash(v^*\NvN v)\oup)\cup v^*(k'\NvN k')\oup v \\
\text{has $v^*y$ on the left and $y^*v$ on the right }&
 \text{belongs to }\Omega_{l,r} \\
\text{has $v^*y$ on the left and $v^*y$ on the right }&
 \text{belongs to }\{h\}\cup\Omega_l.
\endalign $$
Writing out all letters from $\Omega$, using that
$(v^*\NvN v)\oup=v^*(k'\NvN k')\oup v$ and that $y\in k'\NvN h'$, then
cancelling $v^*v$ to $h$ and $vv^*$ to $k$ whenever possible, we see that
$$ x''=x^{(3)}\in\Lambdao(\{v,v^*\},
 (h\NvN h)\oup\cup h\NvN k\cup k\NvN h\cup(k\NvN k)\oup). $$
Each element of $\NvN$ is the s.o.--limit of a bounded sequence in
$\lspan(\{1\}\cup\Psi)$, where $\Psi=\Lambdao(A_1\oup,B_2\oup)$.
Using the $\phi$--preserving conditional expectation from $\NvN$ onto $B_2$,
and the facts that $k\NvN h\subseteq\NvN\oup$, that $h$ and $k$ are minimal
projections in $B_2$ and that $kB_2h=\{0\}$, we have that each element of
$(k\NvN k)\oup\cup k\NvN h\cup h\NvN k\cup(h\NvN h)\oup$ is the s.o.--limit of
a bounded sequence in $\lspan(\Psi\backslash B_2\oup)$.
Thus it suffices to show $\phi(x^{(4)})=0$ for every
$$ x^{(4)}\in\Lambdao\bigl(\Psi\backslash B_2\oup,\{v,v^*\}\bigr). $$
Since $B_2vB_2\subseteq A_2\oup$, we have that
$$ x^{(4)}=x^{(5)}\in\Lambdao(A_1\oup,A_2\oup), $$
so by freeness $\phi(x^{(5)})=0$.
This completes the proof that $w$ is a Haar unitary and that $\{w\}$ and
$h'\PvN h'$ are free in $h'\MvN h'$.

Thus, applying Proposition~\FPTypeIIIFactors{} and using that the central
support of $h'$ in $\MvN$ is $1$, we see that $\MvN$ is a
type~III factor, that $\phi$ is extremal almost
periodic with centralizer isomorphic to $L(\freeF_\infty)$ and
$\Sd(\phi)=\Sd(\phi\restrict_\NvN)$.
This completes the proof in Case~IIa.

In Case~IIb, there is a partial isometry, $y\in\NvN_\phi$, such that
$y^*y\le h'\eqdef r_0h$ and $yy^*=k$.
Let $w=y^*v$.
Then
$$ h\MvN h=W^*(h\NvN h\cup\{w\}). $$
Let
$$ \align
\NvN_{[0,0]}&=h\NvN h \\
\NvN_{[-n,0]}&=W^*(w^*\NvN_{[-n+1,0]}w\cup\NvN_{[0,0]})\;(n\ge1) \\
\NvN_{(-\infty,0]}&=W^*(\bigcup_{n\ge0}\NvN_{[-n,0]}) \\
\NvN_{[0,n]}&=W^*(\NvN_{[0,0]}\cup w\NvN_{[0,n-1]}w^*)\;(n\ge1) \\
\NvN_{[0,\infty)}&=W^*(\bigcup_{n\ge0}\NvN_{[0,n]}) \\
\NvN_{(-\infty,\infty)}&=W^*(w^*\NvN_{(-\infty,0]}w\cup\NvN_{[0,\infty)}).
\endalign $$
Just as in the proof of Proposition~\FPExtrAPfd, one can show that
\roster
\item"(i)" $w^*\NvN_{[-n+1,0]}w$ and $\NvN_{[0,0]}$ are free in $h\MvN h$
$\forall n\ge1$.
\item"(ii)" $w\NvN_{[0,n-1]}w^*$ and $h\NvN_{[0,0]}h$ are free in
$h\MvN h$ $\forall n\ge1$.
\item"(iii)" $w^*\NvN_{(-\infty,0]}w$ and $\NvN_{[0,\infty)}$ are free in
$h\MvN h$.
\item"(iv)"
$\NvN_{(-\infty,\infty)}=h\MvN_\phi(\Gamma)h$ where
$\Gamma=\Sd(\phi\restrict_\NvN)$.
\endroster
Now repeatedly using previously (by inductive hypothesis) proven cases of this
theorem, together with the notion of standard embeddings, we prove that $\MvN$
is a type~III factor, that $\phi$ is extremal almost periodic with centralizer
isomorphic to $L(\freeF_\infty)$ and with $\Sd(\phi)$ equal to the subgroup of
$\Reals^*_+$ generated by $\Sd(\phi\restrict_\NvN)$ and $\lambda$.
This finishes the proof in Case~IIb.

In Case~IIc, note that the facts~(i), (ii) and~(iii) that we proved above are
more generally valid.
Hence we prove that $\NvN_{[-N+1,0]}$ is a type~III factor, that
$\phi\restrict_{\NvN_{[-N+1,0]}}$ is extremal almost periodic with centralizer
isomorphic to $L(\freeF_\infty)$ and with
$\Sd(\phi\restrict_{\NvN_{[-N+1,0]}})=\Sd(\phi\restrict_\NvN)$.
There is $z\in\NvN_\phi(\lambda^{-N})$ such that $z^*z=h$ and
$zz^*=w^N(w^*)^N$.
Let $u=z^*w^N$.
Then as in the proof of Case~C of Proposition~\FPExtrAPfd, we see that
$h\MvN(\Sd(\phi\restrict_\NvN))h=W^*(\NvN_{[-N+1,0]}\cup\{u\})$, that $u$ is a
Haar unitary in
$h\MvN h$ and that $\{u\}$ and $\NvN_{[-N+1,0]}$ are free in $h\MvN h$ with
respect to the restriction of $\phi$.
Hence we prove that $\MvN$ is a type~III factor, that $\phi$ is extremal almost
periodic with centralizer isomorphic to $L(\freeF_\infty)$ and that $\Sd(\phi)$
is the subgroup of
$\Reals^*_+$ generated by $\Sd(\phi\restrict_\NvN)$ and $\lambda$.
This finishes the proof in Case~IIc.

In Case~IV, by inductive hypothesis there must be a central and minimal
projection $p\in A_1$ such that
$\phi(p)>\phi(p')$ whenever $p'$ is a central and minimal projection in either
$A_1$ or $B_2$.
Hence,
$$ (h+k)\NvN(h+k)=\smp{(h+k)\NvN_0(h+k)}{\vphantom{\bigm|}(h+k)r_0}
\oplus\smdp\Cpx{\gamma_1/\phi(h+k)}{r_1}
\oplus\smdp\Cpx{\gamma_2/\phi(h+k)}{r_2} $$ 
with $r_1\le h$ and $r_2\le k$.
By the assumption that $h$ is central in $B_2$ we have
$$ \alignat 2
r_1&=p\wedge h,\qquad&\gamma_1=\phi(r_1)&=\phi(p)+\phi(h)-1 \\
r_2&\le p\wedge k,\qquad&\gamma_2=\phi(r_2)&\le\phi(p)+\phi(k)-1.
\endalignat $$
Hence $\phi(hr_0)=1-\phi(p)\le\phi(kr_0)$ and there exists $y\in\NvN_\phi$ such
that $y^*y=hr_0$ and $yy^*\le kr_0$.
Thus $(h+k)\NvN(h+k)=W^*(k\NvN k\cup\{y,r_1\})$ and using $h=y^*y+r_1$ we see
that
$$ k\MvN k=W^*(k\NvN k\cup\{vy^*,vr_1v^*\}). $$
Let $\PvN=W^*(k\NvN k\cup\{vr_1v^*\})$.
We claim that $k\NvN k$ and $\{vr_1v^*\}$ are free in $k\MvN k$.
Indeed, one can show the stronger fact that $k\NvN k$ and $v\NvN v^*$ are free
using that $\lspan(\{1\}\cup\Lambdao(A_1\oup,B_2\oup))$ is dense in $\NvN$ and
that $k$ and $h$ are minimal projections in $B_2$.
Thus we have that
$$ (\PvN,\phi(k)^{-1}\phi\restrict_\PvN)\cong
(k\NvN_0k\oplus\smdp\Cpx{\gamma_2/\phi(k)}{r_2})
*(\smdp\Cpx{1-(\gamma_1/\phi(h))}{k-vr_1v^*}
\oplus\smdp\Cpx{\gamma_1/\phi(h)}{vr_1v^*}). $$
We see that $\frac{\gamma_2}{\phi(k)}<\frac{\gamma_1}{\phi(h)}$ and hence
$$
\PvN=\cases\PvN_0\oplus
\smdp\Cpx{(\gamma_2/\phi(k))+(\gamma_1/\phi(h))-1}{r_2\wedge vr_1v^*}
&\text{ if }\frac{\gamma_2}{\phi(k)}+\frac{\gamma_1}{\phi(h)}>1 \\
\PvN_0&\text{ otherwise,} \endcases $$
where $\PvN_0$ is nonatomic and where either $\phi\restrict_{\PvN_0}$ is a
trace or $\PvN_0$ is a type~III factor, $\phi\restrict_{\PvN_0}$ is extremal
almost periodic with centralizer isomorphic to $L(\freeF_\infty)$ and
$\Sd(\phi\restrict_{\PvN_0})=\Sd(\phi\restrict_{\NvN_0})$.

We have $k\MvN k=W^*(\PvN\cup\{vy^*\})$ with $yv^*vy^*=yy^*\le kr_0$ and
$vy^*yv^*=k-vr_1v^*$.
Thus if $\frac{\gamma_2}{\phi(k)}+\frac{\gamma_1}{\phi(h)}>1$ then
$r_2\wedge vr_1v^*$ is central and minimal in $k\MvN k$.
Finding the remaining part of $k\MvN k$ ({\it i.e\.} by adding $vy^*$ to
$\PvN_0$), is done in three subcases:
\roster
\item"(IVa)" $\lambda\in\Sd(\phi\restrict_{\NvN_0})$;
\item"(IVb)" $\lambda^\Integers\cap\Sd(\phi\restrict_{\NvN_0})=\{1\}$;
\item"(IVc)" $\lambda^\Integers\cap\Sd(\phi\restrict_{\NvN_0})
 =\lambda^{N\Integers}$ for some $N\ge2$.
\endroster
Let $\kt=k-r_2\wedge vr_1v^*$.
We will now concentrate on proving in each of these three cases that
$\kt\MvN\kt$ is a type~III factor, that $\phi\restrict_{\kt\MvN\kt}$ is an
extremal almost periodic state with centralizer isomorphic to
$L(\freeF_\infty)$ and with $\Sd(\phi\restrict_{\kt\MvN\kt})$ equal to the
subgroup of $\Reals_+^*$ generated by
$\bigl(\text{point spectrum}(\Delta_{\phi_1})
\cup\text{ point spectrum}(\Delta_{\phi_2})\bigr)$.
\demo{Proof in Case IVa}
There is $z\in\PvN_\phi(\lambda^{-1})$ such that $z^*z=yy^*$ and
$zz^*=k-vr_1v^*$.
Let $w=vy^*z^*$.
Then $w^*w=k-vr_1v^*=ww^*$.
Let $k'=w^*w$.
Then $k'\MvN k'=W^*(k'\PvN k'\cup\{w\})$.
\proclaim{Claim \FPFD a}
In $(k'\MvN k',\phi(k')^{-1}\phi\restrict_{k'\MvN k'})$,
$k'\PvN k'$ and $\{w\}$ are $*$--free and
$w$ is a Haar unitary.
\endproclaim
\demo{Proof of Claim~\FPFD a}
It will suffice to show $\phi(x)=0$ whenever
$$ x\in\Lambdao\bigl((k'\PvN k')\oup,
\{w^n\mid n\in\Naturals\}\cup\{(w^*)^n\mid n\in\Naturals\}\bigr). $$
Writing $w=(vy^*)z^*$, we see $x=x'\in\Lambdao(\PvN,\{vy^*,yv^*\})$,
where every letter of $x$ that lies between
$$ \align
vy^*\text{ on the left and $yv^*$ on the right}&\text{ belongs to }
(yy^*\PvN yy^*)\oup \\
yv^*\text{ on the left and $vy^*$ on the right}&\text{ belongs to }
(k'\PvN k')\oup.
\endalign $$
Now every element of
$$ \aligned
\PvN\text{ is the s.o.--limit of a bounded sequence in }
&\lspan(\{k\}\cup\Theta), \\
(yy^*\PvN yy^*)\oup\text{ is the s.o.--limit of a bounded sequence in }
&\lspan\Theta, \\
(k'\PvN k')\oup\text{ is the s.o.--limit of a bounded sequence in }
&\lspan k'\Theta_{l,r}k',
\endaligned \tag{\PTheta} $$
where
$$ \Theta=\Lambdao\bigl((k\NvN k)\oup,
\{k'-\tfrac{\phi(k')}{\phi(k)}k\}\bigr) $$
and $\Theta_{l,r}$ is the set of $a\in\Theta$ whose first and last letters come
from $(k\NvN k)\oup$.

Replacing elements of $\PvN$ with elements of $\{k\}\cup\Theta$ and
$k'\Theta_{l,r}k'$ according to~(\PTheta),
sticking $y$ and $y^*$ to elements of $(k\NvN k)\oup$ whenever
possible, and noting $y^*(k\NvN k)\oup y\subseteq(h\NvN h)\oup$, we see that in
order to show $\phi(x')=0$ it will suffice to show that $\phi(x'')=0$ whenever
$$ x''\in\Lambdao\bigl((k\NvN k)\oup\cup k\NvN h\cup h\NvN k\cup (h\NvN h)\oup,
\{v,v^*\}\bigr). $$
Each element of $\NvN$ is the s.o.--limit of a bounded sequence in
$\lspan(\{1\}\cup\Psi)$, where $\Psi=\Lambdao(A_1\oup,B_2\oup)$.
Using the $\phi$--preserving conditional expectation from $\NvN$ onto $B_2$,
and the facts that $k\NvN h\subseteq\NvN\oup$, that $h$ and $k$ are minimal
projections in $B_2$ and that $kB_2h=\{0\}$, we have that each element of
$(k\NvN k)\oup\cup k\NvN h\cup h\NvN k\cup(h\NvN h)\oup$ is the s.o.--limit of
a bounded sequence in $\lspan(\Psi\backslash B_2\oup)$.
Thus it suffices to show $\phi(x^{(3)})=0$ for every
$$ x^{(3)}\in\Lambdao\bigl(\Psi\backslash B_2\oup,\{v,v^*\}\bigr). $$
Since $B_2vB_2\subseteq A_2\oup$, we have that
$$ x^{(3)}=x^{(4)}\in\Lambdao(A_1\oup,A_2\oup), $$
so by freeness $\phi(x^{(4)})=0$.
This completes the proof of Claim~\FPFD a.
\enddemo

Using the above claim and Proposition~\FPExtrAPfd, then expanding out with
partial isometries from $(\PvN_0)_\phi$, allows us to see the proof of the
required facts about $\kt\MvN\kt$ and $\phi\restrict_{\kt\MvN\kt}$ in
Case~IVa.
\enddemo

\demo{Proof in Case IVb}
Since $vy^*\in\MvN_\phi(\lambda^{-1})$, there is a partial isometry
$z\in(\PvN_0)_\phi$ such that $zz^*=k-vr_1v^*$ and $z^*z\le yy^*$.
Let $w=z^*vy^*$, so that $w^*w=yy^*$ and $ww^*=z^*z\le yy^*$.
Let $k'=w^*w$, $k_1=ww^*\le k'$.
Then $k'\MvN k'=W^*(k'\PvN k'\cup\{w\})$.
Let
$$ \align
\PvN_{[0,0]}&=k'\PvN k' \\
\PvN_{[-n,0]}&=W^*(w^*\PvN_{[-n+1,0]}w\cup\PvN_{[0,0]}),\;(n\ge1) \\
\PvN_{(-\infty,0]}&=W^*(\tsize\bigcup_{n\ge1}\PvN_{[-n,0]}) \\
\PvN_{[0,n]}&=W^*(w\PvN_{[0,n-1]}w^*\cup\PvN_{[0,0]}),\;(n\ge1) \\
\PvN_{(0,\infty]}&=W^*(\tsize\bigcup_{n\ge1}\PvN_{[0,n]}) \\
\PvN_{(-\infty,\infty)}&=W^*(w^*\PvN_{(-\infty,0]}w\cup\PvN_{[0,\infty)}).
\endalign $$
Then
\roster
\item"(i)" $w^*\PvN_{[-n+1,0]}w$ and $\PvN_{[0,0]}$ are free in $k'\MvN k'$
$\forall n\ge1$.
\item"(ii)" $w\PvN_{[0,n-1]}w^*$ and $k_1\PvN_{[0,0]}k_1$ are free in
$k_1\MvN k_1$ $\forall n\ge1$.
\item"(iii)" $w^*\PvN_{(-\infty,0]}w$ and $\PvN_{[0,\infty)}$ are free in
$k'\MvN k'$.
\item"(iv)"
$\PvN_{(-\infty,\infty)}=k'\MvN_\phi(\Gamma)k'$ where
$\Gamma=\Sd(\phi\restrict_{\NvN_0})$.
\endroster

Statements~(i)--(iv) are proved in the same spirit as the proof of Case~IVa
above and
similarly to the proof of Case~B of Proposition~\FPExtrAPfd.
Then Proposition~\FPTypeIIIFactors{} is applied repeatedly to conclude the
required facts about $\kt\MvN\kt$ and $\phi\restrict_{\kt\MvN\kt}$.
\enddemo

In Case~IVc, the required facts about $\kt\MvN\kt$ and
$\phi\restrict_{\kt\MvN\kt}$ are proved similarly to the proof of Case~C of
Proposition~\FPExtrAPfd.

After expanding out with partial isometries from the centralizer of
$\phi\restrict_{\NvN_0}$, we see that the same
facts hold for $\MvN_0$ and $\phi\restrict_{\MvN_0}$, namely, that $\MvN_0$ is
a type~III factor, that $\phi\restrict_{\MvN_0}$ is extremal almost periodic
with centralizer isomorphic to $L(\freeF_\infty)$ and that
$\Sd(\phi\restrict_{\MvN_0})$ is equal to the subgroup of $\Reals_+^*$
generated by
$\bigl(\text{ point spectrum}(\Delta_{\phi_1})
\cup\text{ point spectrum}(\Delta_{\phi_2})\bigr)$.
It remains to verify the size (with respect to $\phi$) of the type~I summand of
$\MvN$ to which $h$ and $k$ contribute.
After changing indices, we may assume that $p$, $h$, $k$ and $v$, are, in the
notation of~(\DoneDtwo),
$$ \align
p&=p_{1,1}\text{ with }n_1=1 \\
h&=q_{1,1}\text{ with }m_1>1 \\
k&=q_{1,2} \\
v&=v_{21}^{(1)}.
\endalign $$
Then
$$ \align
\gamma_1&=\beta_{1,1}\left(1-(1-\alpha_{1,1})\frac1{\beta_{1,1}}\right) \\
\gamma_2&=\beta_{1,2}\left(1-(1-\alpha_{1,1})
 \left(\sum_{i=2}^{m_2}\frac1{\beta_{1,i}}\right)\right) \\
\phi(r_2\wedge vr_1v^*)
 &=\phi(k)\left(\frac{\gamma_1}{\phi(k)}+\frac{\gamma_2}{\phi(h)}-1\right) \\
 &=\beta_{1,2}\left(1-(1-\alpha_{1,1})\left(\frac1{\beta_{1,1}}+
 \sum_{i=2}^{m_2}\frac1{\beta_{1,i}}\right)\right),
\endalign $$
which gives the correct formula for the type~I summand of $\MvN$ to which $h$
and $k$ contribute, thus finishing the proof of the proposition in Case~IV.
\QED

The following is essential for investigating free products with $\Bof(\Hil)$.
\proclaim{Proposition \StdEmbFPFD}
Let
$$ \align
(\MvN,\phi)=&(A_1,\phi_1)*(A_2,\phi_2) \\
\cup \\
(\NvN,\phi\restrict_\NvN)
=&(B_1,\phi_1\restrict_{B_1})*(B_2,\phi_2\restrict_{B_2}),
\endalign $$
where $\phi_1$ and $\phi_2$ are faithful states, $A_1$ and $A_2$ are
finite dimensional, $\Cpx\neq B_1\subsetneqq A_1$, $\Cpx\neq B_2\subseteq A_2$
and at least one of
$\phi_1\restrict_{B_1}$ and $\phi_2\restrict_{B_2}$ is not a trace.
Thus by Proposition~\FPFD,
$$ \alignat 2
\MvN&=\smp{\MvN_0}{r_0}\oplus D\quad&\text{or}\quad\MvN&=\MvN_0 \\
\NvN&=\smp{\NvN_0}{s_0}\oplus C\quad&\text{or}\quad\NvN&=\NvN_0
\endalignat $$
with $\MvN_0$ and $\NvN_0$ type~III factors such that the centralizers of
$\phi\restrict_{\MvN_0}$ and $\phi\restrict_{\NvN_0}$ are isomorphic to
$L(\freeF_\infty)$.
Then $s_0\le r_0$ and the inclusion of centralizers,
$$ (\NvN_0)_\phi\hookrightarrow s_0(\MvN_0)_\phi s_0, \tag{\NzeroInMzero}$$
is a standard embedding.
\endproclaim
As the techniques of proof of~\FPFD{} are similar to those of~\FPExtrAPfd, one
proves that~(\NzeroInMzero) is a standard embedding like in~\FPExtrAPfd.

Now Theorems~\AlsoBofH{} and~\FPInfMany{} are proved by applying previously
proved results and taking inductive limits.

\head \Fullness.  Fullness of free product factors. \endhead

\cite{\MvNZZiv} defined property Gamma for factors of
type~II$_1$ in terms of the existence of nontrivial central sequences, and they
showed that the von Neumann algebra of a free group, $L(\freeF_2)$, is
non--Gamma.
Fullness is a property of factors that was invented by
\cite{\ConnesZZAlmPer} and that when restricted to the~II$_1$ situation
is the property of being non--Gamma.
More recently, in Theorem~11 of \cite{\Barnett}, it was shown, for a free
product of von~Neumann algebras, $(\MvN,\phi)=(A,\phi_A)*(B,\phi_B)$, that
$\MvN$ is full if the centralizers of $\phi_A$ and $\phi_B$ contain enough
unitaries.
His result is based on the ``$14\epsilon$ Lemma'' of \cite{\Pukansky}.
{}From these partial results, and by analogy with the situation for free
products with respect to traces, \cite{\DykemaZZFreeDim}, one would expect that
$\MvN_0$ in Theorems~\MainOne,~\AlsoBofH{} and~\FPInfMany{} is always a full
factor.
Although this has not been proven in general, there are several classes of
cases in which $\MvN_0$ of Theorems~\MainOne,~\AlsoBofH{} and~\FPInfMany{}
is known to be a full factor:
\roster
\item"(i)" if $\MvN_0$ is a
type~III$_\lambda$ factor for $0<\lambda<1$ then if follows
from Proposition~2.3 of \cite{\ConnesZZAlmPer} and the fact that the
centralizer
of $\phi\restrict_{\MvN_0}$ is $L(\freeF_\infty)$ and it thus non--Gamma, that
$\MvN_0$ is full;
\item"(ii)" in many cases, the criteria of Theorem~11 of \cite{\Barnett} will
apply, or can be used in conjunction with Propositions~\FreeProdDirectSum{}
and~\FreeProdDirectSumMatrixAlg;
\item"(iii)" in the case of
$(\MvN,\phi)=(M_2(\Cpx),\phi_1)*(M_2(\Cpx),\phi_2)$, for $\phi_1$ and $\phi_2$
arbitrary faithful states, it is proved in \cite{\DykemaZZSNU} that $\MvN$
is a full factor, the proof being based on Pukansky's ``$14\epsilon$ Lemma.''
\endroster

For example, using only Theorem~11 of \cite{\Barnett} and our
Theorem~\FPInfMany, we can prove the existence of a full factor having
arbitrary countable 
$\Sd$~invariant and an almost periodic state with centralizer isomorphic to
$L(\freeF_\infty)$.
\proclaim{Proposition \ExampleSd}
For every countable subgroup, $\Gamma$, of $\Reals_+^*$, there is a full
type~III factor, $\MvN$, having extremal almost periodic state, $\phi$, such
that $\Sd(\phi)=\Gamma$ and whose centralizer, $\MvN_\phi$, is isomorphic to
$L(\freeF_\infty)$.
\endproclaim
\demo{Proof}
For $n\ge1$ let $0<\lambda_n<1$ be such that
$\Gamma$ is generated by $\{\lambda_n\mid n\ge1\}$, and let
$$ (\MvN,\phi)=\operatornamewithlimits{*}_{n\ge1}(M_2(\Cpx),
\Tr\left(\left(\matrix\tfrac1{1+\lambda_n}&0\\0&\tfrac{\lambda_n}{1+\lambda_n}
\endmatrix\right)\cdot\right)). $$
Grouping the free copies of $M_2(\Cpx)$ into two infinite sets and applying
Theorem~\FPInfMany, we see that $(\MvN,\phi)$ is isomorphic to the free product
of $(\MvN_1,\phi_1)$ and $(\MvN_2,\phi_2)$ with the centralizer of each
$\phi_\iota$ being $L(\freeF_\infty)$.
Hence Barnett's theorem applies to show that $\MvN$ is a full factor, while by
Theorem~\ExtrFreeProd, $\phi$ is extremal almost periodic with centralizer
$L(\freeF_\infty)$ and $\Sd(\phi)=\Gamma$.
\QED

\head\S\Questions.  Questions. \endhead

We conclude this paper with a few questions.

\proclaim{Question \Qisomlambda}\rm
For a given $0<\lambda<1$, are the type~III$_\lambda$
factors that are obtainable by taking free products of various algebras
isomorphic to each other?
One apraoch to this question would be to find the trace scaling automorphism of
the core factor $L(\freeF_\infty)\otimes\Bof(\Hil)$.
\cite{\RadulescuZZlambda} has identified this automorphism in
the case of
$$ (L(\Integers),\tau_\Integers)
*\left(M_2(\Cpx),\Tr\left(\left(
\matrix\tfrac1{1+\lambda}&0\\0&\tfrac\lambda{1+\lambda}\endmatrix
\right)\cdot\right)\right) $$
to be the $\lambda$--scaling automorphism that is an element of his
one--parameter trace--scaling
action \cite{\RadulescuZZOneParamGp}, which he constructed via random matrices.
One might hope that all type III$_\lambda$ factors constructed as free products
of finite dimensional or hyperfinite algebras
in this paper are isomorphic to R\u adulescu's.
It may be that the approach of \cite{\DykemaZZAlmPer} would shed light on this.
\endproclaim

\proclaim{Question \Qcoreone}\rm
What are the continuous cores of the full type~III$_1$  factors obtained in
this paper as free products of
finite dimensional or hyperfinite algebras?
Are the continuous cores of those
that have different $\Sd$~invariants isomorphic or non--isomorphic?
If they were isomorphic, then there would be uncountably many non--outer
conjugate trace scaling actions of $\Reals_+^*$ on this II$_\infty$ factor.
By ``continuous core'' of a type III$_1$ factor, $\MvN$, we mean the
type~II$_\infty$ factor, $\NvN$, such that
$\MvN\cong\NvN\rtimes_\alpha\Reals_+^*$, where $\alpha$ is a trace--scaling
action,
(see \cite{\TakesakiZZDualityCrossProd}).
The discrete core of each of these full factors is
$L(\freeF_\infty)\otimes\Bof(\Hil)$, and
by \cite{\DykemaZZAlmPer} the continuous core is in each case isomorphic to the
fixed point subalgebra of
$L^\infty(\Reals_+^*,L(\freeF_\infty)\otimes\Bof(\Hil))$ under an certain
action of the discrete group $\Sd(\MvN)$.
\endproclaim

\proclaim{Question \Qisomone}\rm
Are the full type III$_1$ factors having the same $\Sd$~invariant that are
obtainable by taking free products of various finite dimensional or hyperfinite
algebras isomorphic to each other?
\endproclaim

\proclaim{Question \Qnonfull}\rm
Is there a type~III$_1$ factor, $\MvN$, with extremal almost periodic state,
$\phi$, such that $\MvN_\phi$ is $L(\freeF_\infty)$ (or any other non--Gamma
II$_1$ factor) but $\MvN$ is not full?
\endproclaim

\refstyle{B}

\enddocument